\documentstyle[11pt, amssymb]{article}
 
\newcommand{\text}{\rm}
\renewcommand{\Bbb}{\bf}
\newcommand{\Det}{{\rm Det}}
\newcommand{\Hom}{{\rm Hom}}
\newcommand{\Gal}{{\rm Gal}}
\newcommand{\Oo}{{\cal O}}
\newcommand{\X}{{\cal X}}
\newcommand{\Q}{{\bf Q}}
\newcommand{\R}{{\bf R}}
\newcommand{\Rr}{{\rm R}}
\newcommand{\Pf}{{\rm Pf}}
\newcommand{\Hh}{{\rm H}}
\newcommand{\func}{\rm}
\newcommand{\tbigcap}{\bigcap}
\newcommand{\tbigsqcup}{\bigsqcup}
\newcommand{\tsum}{\sum}
\newcommand{\tprod}{\prod}

\newcommand{\tbigwedge}{\bigwedge}
\newcommand{\Z}{{\bf Z}}
\newcommand{\Spec}{{\rm Spec}}
\newcommand{\F}{{\cal F}}
\newcommand{\Y}{{\cal Y}}
\newcommand{\G}{{\cal G}}
\renewcommand{\Pf}{{\rm Pf}}
\newcommand{\E}{{\cal E}}
\newcommand{\C}{{\bf C}}
\newcommand{\omlogY}{\Omega_{\Y/\Z}^{1}( \log \Y_{S}^{\rm red}/\log S)}
\newcommand{\omlogX}{\Omega_{\X/\Z}^{1}( \log \X_{S}^{\rm red}/\log S)}
\newcommand{\hK}{\widehat {\rm K}}

\newcommand{\no}{\noindent}

\begin{document}
\title{$\varepsilon$-constants and equivariant Arakelov Euler characteristics}
\author{Ted Chinburg\thanks{Supported in part by NSF grant DMS97-01411.}, Georgios Pappas\thanks{Supported in part by NSF grant DMS99-70378 and by a Sloan Research Fellowship.}, Martin J. Taylor\thanks{EPSRC Senior Research Fellow.}}\date{}
\maketitle

\section{Introduction}

Let $R[G]$ be the group ring of a finite group $G$ over
a ring $R$.  In this article, we study Euler characteristics of bounded 
metrised complexes of finitely generated $\Z[G]$-modules, with applications
to Arakelov theory and the determination of $\epsilon$-constants.  

A metric on a bounded complex $K^\bullet$ of finitely generated $\Z[G]$-modules is
specified by giving for each irreducible character $\phi$ of $G$ a metric
on the determinant of the $\phi$-isotypic piece of the complex ${\bf Q} \otimes_{\bf Z} K^\bullet$
of $\Q[G]$-modules.  Ignoring metrics for the moment,
the alternating sum of the terms of $K^\bullet$ yields an Euler
characteristic in the Grothendieck group G$_0(\Z[G])$ of finitely
generated $\Z[G]$-modules.  If $K^\bullet$ is perfect, in the sense
that all its terms are projective, one has an Euler characteristic
in the finer Grothendieck group K$_0(\Z[G])$ of all finitely generated
projective $\Z[G]$-modules.  To take metrics into account, we will use
a metrized version $A(\Z[G])$ of the projective class group of $\Z[G]$.
We will construct in $A(\Z[G])$ an ``Arakelov-Euler characteristic" 
associated to each bounded perfect
metrised complex of $\Z[G]$-modules. 

Our primary interest will be metrised complexes arising in
the following way.  Let $\mathcal{X}$ be a scheme which is 
projective and flat over $\Spec(\Z)$ and which has smooth
generic fibre. We suppose that $\mathcal{X}$ supports an action by a finite
group $G\,\;$and that the action is tame in the sense that, for each closed
point $x$ of $\mathcal{X}$, the order of the inertia group of $x$ is coprime
to the residue characteristic of $x$.  Choose a $G$-invariant K\"{a}hler
metric $h\;$on the tangent bundle of the associated complex manifold $%
\mathcal{X}(\C) $. We are then able to construct an
Arakelov-Euler characteristic for any hermitian $G$-bundle $( {\cal F},j)$ on 
$\mathcal{X}$\ by endowing the equivariant determinant of
cohomology of $\Rr\Gamma  ( \mathcal{X},\mathcal{F} ) $ with
equivariant Quillen metrics $j_{Q,\phi }$ for each irreducible character $%
\phi $ of $G$. This construction can  be extended to give an
Arakelov-Euler characteristic for a bounded complex of hermitian $G$%
-bundles. Roughly speaking, our main results show that the
Arakelov-Euler characteristic of the (logarithmic) de Rham complex of 
$\mathcal{X}$ determines certain $\varepsilon$-constants
which are associated to the L-functions of the Artin motives obtained from $%
\mathcal{X}$ and the symplectic representations of $G$.

Let us now describe these results in more detail. Firstly we need to say  a little more
about arithmetic classgroups.  Let $R_G$ (resp. $R_{G}^{s}$) be the group of virtual
characters (resp. virutal symplectic characters) of
$G$.  In Section 4 we obtain a quotient $A^{s}(\Z[G])$ of
$A(\Z[G]) $ called the
 symplectic arithmetic classgroup  by restricting  
functions on $R_G$ to $R_G^s$.
We show that $A^{s} ( \Z[G] ) $ contains
a subgroup $R(\Z[G])$, called the group of rational
classes, which supports a natural isomorphism 
$\theta :R(\Z[G]) \rightarrow {\rm Hom}_{\Gal} ( R_{G}^{s},\Bbb{Q}^{\times } ) .$
Now let $S$ denote a finite set of primes which includes those primes where 
$\mathcal{X}$ has non-smooth reduction.  Let $\Omega _{\mathcal{X}\text{%
/}\Bbb{Z}}^{1}\left( \log \X^{\text{red}}_S/\log
S\right) $ denote the sheaf of degree one relative logarithmic differentials
of $\mathcal{X\;}$with respect to the morphism $\left( \X,
{\X}_{S}^{\text{red}}\right) \rightarrow \left( \Spec\left( 
\Z\right), S\right) $ of schemes with log-structures. 
This sheaf is locally free if (as we now assume) each special fibre of $\mathcal{X\;}$is a divisor
with strictly normal crossings and the multiplicities of the irreducible
components of each fibre are prime to the residue characteristic. The logarithmic de Rham complex
$\Omega^\bullet_{\X/\Z}(\log \X^{\rm red}_S/\log S)$ is defined to be the complex 
\[
\Oo_{\mathcal{X}}\rightarrow \omlogX  \rightarrow \cdots
  \rightarrow \Omega^d_{\X/\Z}( \log 
\X_{S}^{\text{red}}/\log S) 
\]
We view each term of this complex as carrying the hermitian metric given by the
corresponding exterior power of $h^{D}.$ In Theorem 7.1 we completely
describe the image ${\frak{c}}^{s}\;$of the logarithmic de Rham 
Arakelov-Euler characteristic $%
{\frak{c}}=\chi  ( {\Rr}\Gamma({\cal X}, \Omega _{\X/\Bbb{Z}}^{\bullet
} ( \log \X_{{S}}^{\text{red}}/\log S )) ,\wedge
^{\bullet }h^{D}_Q ) $ in the symplectic arithmetic classgroup $%
A^{s}(\Z[G]) .$ Here we limit ourselves to stating the result
for characters of degree zero, since this result is particularly striking
(see Theorem 7.1 for the result for characters of arbitrary degree): 
\medskip

\no{\bf Theorem.} {\sl The element ${\frak{c}}^{s}$ of $A^s(\Z[G])$
is a rational class and for any
virtual symplectic character $\psi $ of degree zero
\[
\theta \left( {\frak{c}}^{s}\right) \left( \psi \right) =\varepsilon_{0}(\Y,\psi)^{-1} . 
\]}

Whilst this result is of interest in itself, our principal
concern rests with the sheaf of differentials $\Omega _{\mathcal{X}\text{/}%
\Bbb{Z}}^{1}$.  One can associate to $\Omega
_{\mathcal{X}\text{/}%
\Bbb{Z}}^{1}$ and the metric $h^D$ a class $\Omega$ in the arithmetic
Grothendieck group $\hK_0(\X)$, which is a $\lambda$-ring by 
[GS1,\S 7]. In [CPT2] we show 
$$\sum_{i=0}^d \chi(\lambda^i(\Omega)) = -{\rm log}|\varepsilon(\X)|.$$
when $G$ is trivial and $\X$ satisfies the previous fibral hypotheses, where 
$\varepsilon(\X)$ is the constant in the
functional equation of the Hasse-Weil zeta function of $\X$.  Furthermore, under 
much more general circumstances,
we show that such an equality is equivalent to Bloch's Conjecture
characterising the conductor of $\mathcal{X}$/$\Bbb{Z}$ in terms of the top
localised Chern class of the differentials of $\mathcal{X}$/$\Bbb{Z}$. The
results of this present paper give equivariant
refinements to the results of [CPT2]. Our main result, Theorem 8.3,
shows that the equivariant  de Rham Arakelov-Euler 
characteristic ${\frak{d}}\in A(\Z[G])$
of $\mathcal{X}$/$\Bbb{Z}$, together with the arithmetic
ramification class ${\rm AR}(\X)\in A(\Z[G])$ (see 8.2), completely determines the $\varepsilon $-constants
(together with their sign) for
the symplectic representations of the group $G$, which acts on $\mathcal{X}$.
\medskip

\no{\bf Theorem.} {\sl The element ${\frak{d}}^{s}\cdot {\rm AR}^s(\X)^{-1}$
of $A^s(\Z[G])$
is a rational class and for any
 symplectic character $\psi $ 
\[
\theta \left( {\frak{d}}^{s}\cdot {\rm AR}^s(\X)^{-1}\right) \left( \psi \right) =\varepsilon(\Y,\psi)^{-1} . 
\]}

This result can be thought of as a ``converse" to the main Theorems of [CEPT1]
and [CPT1]; there the class of the de Rham Euler characteristic
in ${\rm K}_0(\Z[G])$ is shown to be determined by $\varepsilon$-factors.
Let us point out that the arithmetic ramification class ${\rm AR}(\X)$
is modeled on the ``ramification class" ${\rm R}(\X/\Y)$ of [CEPT1]. It
only depends on the branch locus of the cover $\X\rightarrow\Y$;
under our assumptions this branch locus is contained on a finite set
of fibers of $\Y\rightarrow\Spec(\Z)$.

In the classical case of the ring of integers $\Oo_{N}\;$of a tame Galois
extension of number fields $N/K\;$with $G=\Gal\left( N/K\right) $, which
corresponds to ${\cal X}=\Spec( \Oo_{N}) $, one has
Fr\"{o}hlich's hermitian conjecture (now the theorem of
Cassou-Nogu\`{e}s-Taylor in [CNT]). Roughly speaking, this
shows that the $\Z[G]$-module $\Oo_N$ together with the additional structure
provided by the hermitian pairing of the trace form can be used to
determine the symplectic $\varepsilon$-constants.  A main observation of the present paper is that
the role of the trace form is not so central; we can replace the trace form
by Arakelov hermitian metrics. In fact, the arithmetic classgroup $A(\Z[G])$ may be considered as
generalisation of Fr\"{o}hlich's hermitian classgroup. There are, however,
two crucial differences: firstly, we consider complex valued hermitian
forms, whereas his forms are rational valued; the second key-feature of our
approach is that it brings in all the characters of $G$, and not only the
symplectic characters. From our point of view, the 
Cassou-Nogu\`{e}s-Taylor result may be reformulated as follows. They
show that the signs of the Artin L-functions for symplectic representations
of $G$ can be recovered from the isomorphism class of $\Oo_{N}$ as a metrised $\Z[G]$-module. 
Here a metrised $\Z[G]$-module is a $\Z[G]$-module $M $ 
together with a $G$-invariant metric on$\;\Bbb{C}$ $\otimes _{\mathbf{Z}}
M$.  For the ring of integers $\Oo_{N},$ the metric on $\C
\otimes _{\Z}\Oo_{N}$ is given by $z\otimes a\mapsto \left(
\tsum_{\sigma }\left| z\sigma  ( a\right) \right| ^{2} ) ^{1/2}$,
where the sum extends over the distinct embeddings $\sigma :N\rightarrow 
\Bbb{C}$. The theorems in this article
provide generalizations of this Cassou-Nogu\`{e}s-Taylor
result to higher dimensions. Let us point out that, in higher dimensions, a partial result 
of similar flavour was obtained in [CEPT2]; however, the information, which was
used there to characterise the $\varepsilon_0$-constants, does not come
directly from the de Rham cohomology.

This article is structured as follows: in Section 2 we define our notation
and present a number of preliminary results. Then, in Section 3, we define
the arithmetic classes for suitable bounded $\Z[G]$-complexes and
establish a number of their basic properties. The construction of the arithmetic
class is a
rather delicate matter, since we wish to produce an invariant which reflects
the fact that the terms in the complex are projective, whilst the metrics
are only defined on the determinants of cohomology. The main point here is
to show that our notion of arithmetic class is invariant under
quasi-isomorphisms which preserve metrics in an appropriate sense.
The formation of arithmetic classes may also be seen
to be closely related to the refined Euler characteristics with values in
relative K-groups introduced by\ D. Burns in [Bu1] and [Bu2].

Our arithmetic classes take values in the arithmetic classgroup $A(\Z[G]))$.
In
practice, it is often convenient to work with quotient groups of $A(\Z[G])$ such as
those in Section 4.

In Section 5 we consider an arithmetic variety $\X$ which carries a
tame action by a finite group $G$ and we define an arithmetic class for a
hermitian $G$-bundle on $\X$ which supports a set of metrics on the
equivariant determinant of cohomology; we then carry out a number of
calculations in the case when $\X$ is the spectrum of a ring of
integers. In Section 6 we fix a choice of K\"{a}hler metric on the
tangent bundle of $\X(\C) $. For a complex $({\cal G}^{\bullet },h_{\bullet })$ of hermitian $G$-bundles on 
$\X$,  we use the equivariant Quillen metrics on the equivariant
determinants of hypercohomology (see [B]) to construct an arithmetic class $
\chi ( {\rm R}\Gamma {\cal G}^{\bullet },h_{ Q\bullet } ) $; we
call this class the Arakelov-Euler characteristic of $ ( {\cal G}%
^{\bullet }, h_{\bullet } ) $. We then briefly detail the functorial
properties of such Euler characteristics and calculate such Euler
characteristics when $\X$ has dimension one. In Section 7 we
describe the arithmetic class associated to the de Rham complex of
logarithmic differentials; this then enables us to derive our
characterisation of the symplectic $\varepsilon _{0}$-constants of $\X$. 
Finally, in Section 8, we consider the de Rham Arakelov-Euler
characteristic associated to the (regular) differentials of $\X/\Z$, 
and we show how this arithmetic class determines the symplectic $\varepsilon $-constants of $\X$.

The third author wishes to express his gratitude to Christophe Soul\'{e} for
helpful conversations and exchanges of letters concerning the work presented
in the latter part of this article.

\newpage

\section{Preliminary results}

\noindent {\bf 2.A. \ Hermitian complexes.}
\medskip

Let $R$ denote a commutative ring which is endowed with a fixed embedding
into the field of complex numbers $\C$; in applications $R$ will be
either $\Z$, ${\bf R}$ or $\C$. We consider bounded (cochain)
complexes $K^{\bullet }$ of finitely generated left $R[G]$-modules 
\[
K^{\bullet }\ :\ \cdots\rightarrow K^{i}\stackrel{d_{K}^{i}}{\rightarrow }K^{i+1}\rightarrow \cdots  
\]
so that the boundary maps $d_{K}^{i}$ are all $R[G]$-maps. Thus the $i$th
cohomology group, denoted $\Hh^{i}=\Hh^i(K^{\bullet })$, is an $R[G]$-module.
Recall that the complex $K^{\bullet }$\ is called {\sl perfect} if in
addition all the modules $K^{i}$ are $R[G]$-projective. 

\medskip 
\no {\bf Definition 2.1}  Let $x\mapsto \overline{x}$ denote
the complex conjugation automorphism of $\Bbb{C}$; we extend complex
conjugation to an involution of the complex group algebra $\C[G]$ by the
rule $\overline{\sum a_{g}g}=\sum \overline{a_{g}}g^{-1}.$ An element $x\in 
\C[G]$ is called \textit{symmetric} if $\overline{x}=x.$ A  
{\sl hermitian $R[G]$-complex} is a pair ($K^{\bullet }, k^{\bullet })$
where $K^{\bullet }$ is an $R[G]$-complex, as above, and where each $K_{\C
}^{i}=\C\bigotimes_{R}K^{i}$ is endowed with a non-degenerate
positive-definite $G$-invariant hermitian form 
\[
k^{i}:K_{\C}^{i}\times K_{\C}^{i}\rightarrow \C. 
\]
Thus, in particular, each $k^{i}$ is left $\Bbb{C}$-linear and $\overline{%
k^{i}(x,y)}=k^{i}(y,x).$ The metric associated to $k^{i}$ is defined by $%
\left\| x\right\| ^{i}=\sqrt{k^{i}\left( x,x\right) }$ for $x\in K_{\Bbb{C%
}}^{i}$, where $k^{i}\left( x,x\right) \geq 0$ because $k^{i}$ is a positive
definite hermitian form.
\medskip

Equivalently (as per p 164 in [F2]) we may work with the $\C[G]$-valued
hermitian forms 
\[
\widehat{k}^{i}:K_{\Bbb{C}}^{i}\times K_{\Bbb{C}}^{i}\rightarrow \C[G] 
\]
given by the rule that for $x,y\in K_{\Bbb{C}}^{i}$
\[
\widehat{k}^{i}(x,y)=\sum_{g\in G}k^{i}(x,gy)g. 
\]
\smallskip Thus $\widehat{k}^{i}$ is $\C[G]$-left linear and is reflexive
in the sense that $\overline{\widehat{k}^{i}(x,y)}=\widehat{k}^{i}(y,x).$
Conversely, given $\widehat{k}^{i}$, we may of course recoup $k^{i}$ by
reading off the coefficient of $1_{G}$ in $\C[G]$.
\medskip

\no{\bf Example 2.2}\ The module $\C[G]$ carries the so-called \textit{standard}
positive $G$-invariant hermitian form 
\[
\mu :\C[G]\times \C[G]\rightarrow \C
\]
given by the rule $\mu \left( \sum x_{g}g,\sum y_{h}h\right) =\tsum x_{g}%
\overline{y}_{g}$. Then the associated $\C[G]$-valued hermitian form $%
\widehat{\mu }$ 
\[
\widehat{\mu }:\C[G]\times \C[G]\rightarrow \C[G] 
\]
is the so-called multiplication form $\widehat{\mu }\left( x,y\right) =x\cdot
\overline{y}$.
\bigskip

\no{\bf 2.B. \ Metrised complexes.}
\medskip

Let $G$ again denote a finite group, let $\widehat{G}$ denote the set of
irreducible complex characters of $G$, and once and for all for each $\phi
\in \widehat{G}$ we let $W=W_{\phi }$ denote the simple 2-sided $\C[G]$%
-ideal with character $\phi \left( 1\right) \overline{\phi }$, where $%
\overline{\phi }$ is the contragredient character of $\phi .\;$For a
finitely generated $\C[G]$-module $M$ we define $M_{\phi }=\left(
M\otimes _{\Bbb{C}}W\right) ^{G}$, where $G$ acts diagonally and on the left
of each term; more generally, for a bounded complex $P^{\bullet }$ of
finitely generated $\C[G]$-modules, we put $\Hh^i=\Hh^i(
P^{\bullet }) $ and we write 
\[
P_{\phi }^{\bullet }=\left( P^{\bullet }\otimes _{\Bbb{C}}W\right) ^{G}\;\;\;%
\hbox{\rm and \ \ } \Hh_{\phi }^{i}=\left( {\Hh}^{i}\otimes _{\Bbb{C}}W\right)
^{G}. 
\]
We then construct the complex lines $\det \left( P^{\bullet }\left( \phi
\right) \right) $ and $\det \left( \Hh^{\bullet }\left( \phi \right)
\right) $ such that 
\[
\det \left( P^{\bullet }\left( \phi \right) \right) ^{\otimes \phi \left(
1\right) }=\det \left( P_{\phi }^{\bullet }\right) =\otimes _{i}\left(
\wedge ^{{\rm top}}P_{\phi }^{i}\right) ^{(-1)^{i}}\;\;\; 
\]
\[
\hbox{\rm and \ \ }\det \left( \Hh^{\bullet }\left( \phi \right) \right)
^{\otimes \phi \left( 1\right) }=\det \left( \Hh_{\phi }^{\bullet
}\right) =\otimes _{i}\left( \wedge ^{{\rm top}}\Hh_{\phi }^{i}\right)
^{(-1)^{i}} 
\]
where for a complex vector space $V$ of dimension $d,$ $\wedge ^{\text{top}%
}V$ denotes $\wedge ^{d}V$ and where for a complex line $L$ we write $L^{-1%
\text{ }}$ for the dual line Hom$\left( L,\mathbf{C}\right) .\;$ Thus $%
P_{\phi }^{\bullet }$ should be thought of as corresponding to the character 
$\phi \left( 1\right) \overline{\phi }$, whereas $P^{\bullet }\left( \phi
\right) $ corresponds to the character $\overline{\phi }$. Observe that a
metric on $\det  ( P_{\phi }^{\bullet } ) $ determines a unique
metric on $\det  ( P^{\bullet } ( \phi )  ) $ and of
course conversely a metric on $\det  ( P^{\bullet } ( \phi  )
 ) $ determines a unique metric on $\det  ( P_{\phi }^{\bullet
} ) .\;$Note also that here and in the sequel for two finite
dimensional vector spaces $V_{i}$ of dimension $d_{i}$ we normalise the
standard isomorphism $\wedge ^{d_{1}d_{2}}\left( V_{1}\otimes V_{2}\right)
\cong \wedge ^{d_{1}d_{2}}\left( V_{2}\otimes V_{1}\right) \;$by multiplying
by $\left( -1\right) ^{d_{1}d_{2}}$, in order to avoid subsequent sign
complications. We refer to the set of lines $\det ( \Hh_{\phi
}^{\bullet }) $ as the \textit{equivariant determinants of cohomology }%
of $P^{\bullet }$. From Theorem 2 in [KM] we have a canonical isomorphism 
\begin{equation}
\xi _{\phi }:\det \left( P_{\phi }^{\bullet }\right) \cong \det \left(  
\Hh_{\phi }^{\bullet }\right) .
\end{equation}

\medskip For ease of computation we use the above definition of $P_{\phi
}^{\bullet }$; however, alternatively one can also work with the isotypical
components $\overline{W}P^{\bullet },$ as shown in the following lemma. Here
and in further applications we shall often need the renormalised form $\nu :
\C[G]\times \C[G]\rightarrow \C$ of the hermitian form $\mu $ of
(2.2) given by 
\[
\nu \left( x,y\right) =\left| G\right|\cdot  \mu \left( x,y\right) \;\;\hbox{\rm for\ }%
x,y\in \C[G]. 
\]

\medskip

\no{\bf Lemma 2.3} \ {\sl For a $\C[G]$  module $V$ 
with a $G$-invariant metric $\left\| -\right\| ,$ the
natural isomorphism  
\[
\alpha :\left( V\otimes _{\Bbb{C}}W\right) ^{G}\cong \overline{W}V 
\]
given by $\alpha \left( \tsum_{i}v_{i}\otimes w_{i}\right)
=\tsum_{i}\overline{w_{i}}v_{i}$ is an isometry, where both terms
carry the natural metrics induced by $\nu $ and $\left\| -\right\| $; 
that is to say $\overline{W}V$ carries the metric given by the
restriction of $\left\| {-}\right\| ,$  and $\left( V\otimes _{\Bbb{%
C}}W\right) ^{G}$  carries the metric given by the restriction of
the tensor metric associated to $\left\| {-}\right\| $ and $\nu $
on $V\otimes \C[G]$.}
\medskip

\no{\sc Proof.} Let\ $\left\| -\right\| _{1}$ resp. $\left\| -\right\| _{2}$ denote
the given metric on $\left( V\otimes _{\Bbb{C}}W\right) ^{G}$ resp. 
$\overline{W}V$. If $e=|G|^{-1}\cdot \tsum_{g}\phi \left(
g\right) g$ is the central idempotent associated to $W,$ then\ for $x\in 
\overline{W}V,$ we have $\alpha ^{-1}\left( x\right) = |G|^{-1}\cdot \tsum_{g\in G}gx\otimes ge$ and so 
\[
\alpha ^{-1}\left( x\right) =\frac{1}{\left| G\right| ^{2}}\tsum_{g,h\in
G}gx\otimes g\phi \left( h\right) h 
\]
\[
=\frac{1}{\left| G\right| ^{2}}\tsum_{g,h\in G}ghh^{-1}\phi \left( h\right)
x\otimes gh 
\]
\[
=\frac{1}{\left| G\right| }\tsum_{f\in G}f\overline{e}x\otimes f=\frac{1}{%
\left| G\right| }\tsum_{f\in G}fx\otimes f. 
\]
Thus 
\[
\left\| \alpha ^{-1}\left( x\right) \right\| _{1}^{2}=\frac{1}{\left|
G\right| }\tsum_{g}\left\| x\right\| _{2}^{2}=\left\| x\right\|
_{2}^{2}.\;\;\;\Box 
\]

\medskip

\no{\bf Definition 2.4} \ Let $R$ again denote a subring of $\Bbb{C}$. A%
\textit{\ metrised }$R[G]$-\textit{complex} is a pair $( P^{\bullet
},p_{\bullet }) $ , where $P^{\bullet }$ is a bounded complex of
finitely generated (not necessarily projective) $R[G]$-modules and the $%
p_{\phi }$ are a set of metrics given by positive definite hermitian forms
on the complex lines $\det ( \Hh_{\phi }^{\bullet }) $, one
for each $\phi \in \widehat{G}$.
\medskip

Let $R_{G}$ denote the group of complex virtual characters of $G$. For each $%
\phi \in \widehat{G},$ let $p( \phi ) $ denote the metric on $%
\det P^{\bullet }\left( \phi \right) $ induced by $p_{\phi }$. For a virtual
character $\chi =\tsum_{\phi }n_{\phi }\phi \in R_{G},$ we write 
\[
\det \left( P^{\bullet }\left( \chi \right) \right) =\otimes _{\phi }\det
\left( P^{\bullet }\left( \phi \right) \right) ^{\otimes n_{\phi }} 
\]
and we endow this complex line with the\ product metric $p\left( \chi
\right) =\otimes _{\phi }p\left( \phi \right) ^{n_{\phi }}.$

\medskip

\no{\bf Example 2.5} \ A hermitian complex $( K^{\bullet }, k^{\bullet
}) $ affords a metrised complex in the following way: endow $(
K^{i}\otimes _{\Bbb{C}}W) ^{G}$ with the form induced by $k^{i}$ on $K^{i}$ 
and by the restriction of the standard form on $W$, which is given by
the restriction of $\nu$. The alternating tensor product of the top
exterior products of these forms is then a positive definite hermitian form
on the complex line $\det  ( K_{\phi }^{\bullet } ) $ and so
induces a positive definite hermitian form on the complex line $\det  ( 
 \Hh_{\phi }^{\bullet } ) $ via $ ( 1 ) $; this then\
determines a unique metric on the complex line $\det  (  \Hh^{\bullet } ( \phi  )  )$.
Each such form determines a
metric which we denote $\det \left( k\left( \phi \right) \right) .$

\bigskip

\section{ Arithmetic Classes}

\no {\bf 3.A. \ The arithmetic classgroup.}
\medskip

In this sub-section we shall define the arithmetic classgroup in which our
arithmetic classes take their values.

The notation is that of [CEPT2] and so we recall it only briefly: $R_{G}$
denotes the group of complex characters of $G$;\, $\overline{\Bbb{Q}}$ is the
algebraic closure of $\Bbb{Q}$ in $\C$, so that we have the
inclusion map $\overline{\Bbb{Q}}\Bbb{\hookrightarrow }\,\C$. We
set $\Omega =\Gal ( \overline{\Bbb{Q}}/\Bbb{Q} )$; $J_{f}$
is the group of finite ideles in $\overline{\Bbb{Q}}$, that is to say the
direct limit of the finite idele groups of all algebraic number fields $E$
in $\overline{\Bbb{Q}}$.

Let $\widehat{\Z}=\prod_{p}\Z_{p}$ denote the ring of integral
finite ideles of $\Z$. For $x\in $ $\widehat{\Z}G^{\times }$, 
the element $\Det (x) \in  \Hom_{\Omega } 
R_{G},J_{f} ) $ is defined by the rule that for a representation $T$
with character $\psi $%
\[
\Det ( x )  ( \psi  ) =\det  ( T(x)) ; 
\]
the group of all such homomorphisms is denoted 
\[
\Det( \widehat{\Bbb{Z}}G^{\times } ) \subseteq \Hom
_{\Omega } ( R_{G},J_{f} ) . 
\]
More generally, for $n>1$ we can form the group $\Det\left( GL_{n} ( 
\widehat{\Bbb{Z}}G )  \right) ;$ as each group ring $\Z_p[G]$ is
semi-local we have the equality $\Det\left( GL_{n} ( \widehat{\Bbb{Z}}%
G ) \right) =\Det ( \widehat{\Bbb{Z}}G^{\times } ) $ (see
1.2.6 in [T2]).

For an $n\times n$ invertible matrix $A$ with coefficients in $\C[G]$, 
$\left| \Det\left( A\right) \right| 
 \in \; \Hom\left( R_{G},{\bf R}%
_{>0}\right) $ is defined by the rule 
\[
\left| \Det( A) \right| \left( \psi \right) =\left|  
\Det\left( A\right) \left( \psi \right) \right| . 
\]

\no{\bf Lemma 3.1} \ {\sl Extending the involution $x\mapsto \overline{x%
}$ on $\C[G]$  to matrices over $\C[G]$ by
transposition, for $\psi \in R_{G}$}  
\[
\left| \Det( \overline{A}) \right| \left( \psi \right)
=\left| \Det( A) \right| \left( \psi \right) . 
\]
\medskip

Replacing the ring $\widehat{\Bbb{Z}}$ by $\Bbb{Q}$, in the same way we
construct 
\[
\Det( \Q[G[^{\times } ) \subseteq \Hom_{\Omega
} ( R_{G},\overline{\Bbb{Q}}^{\times } ) . 
\]

The product of the natural maps $\overline{\Bbb{Q}}^{\times }\rightarrow
J_{f}$ and $\left| -\right| :\overline{\Bbb{Q}}^{\times }\rightarrow {\bf R}
_{>0}$ yields an injection 
\[
\Delta : \Det (\Q[G]^{\times })\rightarrow \Hom_{\Omega
}(R_{G},J_{f})\times \text{Hom}(R_{G},{\bf R}_{>0}). 
\]
\medskip
\no{\bf  Definition 3.2 } \ The \textit{arithmetic classgroup} 
$A(\Z[G])$ is defined to be the quotient group 
\begin{equation}
A(\Z[G])=\left( \frac{\Hom_{\Omega }(R_{G},J_{f})\times \Hom
(R_{G},\R_{>0})}{\left( \Det(\widehat{\Z}[G]^{\times })\times
1\right) \ {\rm Im}(\Delta )}\right) .
\end{equation}
\smallskip

\no{\bf Remarks.} (1) Note that, in the case when $G=\left\{ 1\right\} ,$\ $A(\Z)$ coincides with the usual 
Arakelov divisor class group of $\Spec(\Z)$ (see 7.7 for further details).

(2) As indicated in the Introduction, there are two crucial differences
between this arithmetic classgroup and the hermitian classgroup of
Fr\"{o}hlich (see II.5\ in\ [F2]): firstly, we work with positive definite
complex hermitian forms; secondly, as a consequence of this, we are able to
work in a uniform manner with all characters of $G.$\medskip

\medskip

\no {\bf 3.B. \ The arithmetic class of a complex.}
\medskip

Let $\left( P^{\bullet },\;p_{\bullet }\right)$ be a perfect metrised $\Z[G]$-complex; 
that is to say $P^{\bullet }$ is a bounded metrised
complex all of whose terms are finitely generated projective (and therefore
locally free) $\Z[G]$-modules. For each $i$, suppose that $d_{i}$ is the
rank of $P^{i}$ as a $\Z[G]$-module, and choose bases $ \{
a^{ij} \} ,\;$resp. $ \{ \alpha _{p}^{ij} \} $ of 
\[
\Q\otimes P^{i}=\tsum_{j}\Q[G]\cdot a^{ij}, \ \hbox{\rm resp.\ }\Z
_{p}\otimes P^{i}=\tsum_{j}\Z_{p}[G]\cdot \alpha _{p}^{ij} 
\]
over $\Q[G]$ resp. $\Z_p[G]$. As both $\{ a^{ij}\}$ and 
$\{ \alpha _{p}^{ij}\} $ are $\Q_{p}[G]$-bases of $\Q_{p}\otimes P^{i}$, we can find 
$\lambda _{p}^{i}\in GL_{d_{i}} ( \Q_{p}G ) $ such that $ ( a^{ij} ) _{j}=\lambda _{p}^{i} (
\alpha _{p}^{ij} ) _{j}$, where $ ( a^{ij} ) _{j}$ denotes the
column vector with $j$-th entry $a^{ij}.$

For $a\in \Q\otimes P^{i}$ we put 
\begin{equation}
r\left( a\right) =\tsum ga\otimes g\in P^{i}\otimes \Q[G].
\end{equation}
Note that for $h\in G$
\begin{equation}
r\left( ha\right) =r\left( a\right) \left( 1\otimes h^{-1}\right)
\end{equation}
and for $w\in W,$ the action of $r\left( a\right) $ on $1\otimes w$ is
defined to be 
\begin{equation}
r\left( a\right) \left( 1\otimes w\right) =\tsum_{g}ga\otimes gw\in \left(
P^{i}\otimes W\right) ^{G}.
\end{equation}
For each $\phi \in \widehat{G}$ we choose an orthonormal basis $\left\{
w_{\phi ,k}\right\} $ of $W$ $=W_{\phi }\;$with respect to the standard form 
$\nu $\ on $\C[G]$,\ then the $\left\{ r\left( a^{ij}\right) \left(
1\otimes w_{\phi ,k}\right) \right\} $ form a $\Bbb{C}$-basis of $\left(
P^{i}\otimes W\right) ^{G}.$ By $\left( 4\right) $ and by linearity we have
that for $\eta =\tsum_{h\in G}\eta _{h}h\in \Q[G]$
\begin{equation}
r\left( \eta a\right) \left( 1\otimes w\right) =\tsum_{h,g}\eta
_{h}gha\otimes gw=r\left( a\right) \left( 1\otimes \overline{\eta }w\right) .
\end{equation}
In the sequel for given $i$ we shall write $\bigwedge r ( a^{i} )
 ( 1\otimes w_{\phi } ) $ for the wedge product 
\[
\wedge _{j,k}\left( r\left( a^{ij}\right) \left( 1\otimes w_{\phi ,k}\right)
\right) \in \det \left( P^{i}\otimes W\right) ^{G}. 
\]
We again adopt the notation of 2.B and let $\left( P^{\bullet },p_{\bullet
}\right) $ be a perfect metrised $\Z[G]$-complex; recall from (1) that
for each $\phi \in \widehat{G}$ we have an isomorphism 
\[
\xi _{\phi }:\det  ( P_{\phi }^{\bullet } ) \cong \det  ( \Hh_{\phi }^{\bullet } ) . 
\]

\medskip
\no{\bf Definition 3.3} \ With the above notation and hypotheses $%
\chi  ( P^{\bullet },p_{\bullet } ) $, the arithmetic class of $%
 ( P^{\bullet },p_{\bullet } ) ,$ is defined to be that class in $A(\Z[G])$ 
represented by the homomorphism on $R_{G}$ which maps
each $\phi \in \widehat{G}$ to the value in $J_{f}\times \R_{>0}$%
\[
\tprod_{p}\left( \tprod_{i}\Det(\lambda _{p}^{i})(\phi
)^{(-1)^{i}}\right)\ \times \ p_{\phi }\left( \xi _{\phi }\left( \otimes
_{i}\left( \bigwedge r( a^{i}) ( 1\otimes w_{\phi })
\right) ^{(-1)^{i}}\right) \right) ^{\frac{1}{\phi (1)}}. 
\]
\medskip

In the sequel we shall refer to the first coordinate as the \textit{finite}
coordinate and the second coordinate as the \textit{archimedean} coordinate.
In order to verify that this class is well-defined, we now show that it is
independent of choices:

(i) If $\{ \widetilde{\alpha }_{p}^{ij}\} $ is a further set of
$\Z_p[G]$-bases for the $\Z_p\otimes P^{i}$, then we can find
$z_{p}^{i}\in GL_{d_{i}} ( \Z_p[G] ) $ such that 
\[
 ( \widetilde{\alpha }_{p}^{ij} ) _{j}=z_{p}^{i} ( \alpha
_{p}^{ij} ) _{j} 
\]
and so the $p$-component\thinspace of the finite coordinate of the
homomorphism representing the class only changes by 
\[
\prod_{i}\Det(z_{p}^{i})^{(-1)^{i}}\in \Det(\Z_p[G]^{\times
}). 
\]

(ii) If $ \{ \widetilde{a}^{ij} \} $ is a further set of $\Q[G]$%
-bases for the $\Q\otimes P^{i}$, then we can find $\eta ^{i}\in
GL_{d_{i}} ( \Q[G]) $ such that 
\[
 ( \widetilde{a}^{ij} ) _{j}=\eta ^{i} ( a^{ij} ) _{j}. 
\]
Now for each pair $i,j,$ we have the equality $\widetilde{a}
^{ij}=\tsum_{l}\eta _{jl}^{i}a^{il}$ and so by (6) we get 
\[
r ( \widetilde{a}^{ij} )  ( 1\otimes w_{\phi ,k} )
=\tsum_{l}r ( a^{il} )  ( 1\otimes \overline{\eta }
_{jl}^{i}w_{\phi ,k} ) ; 
\]
hence 
\[
\tbigwedge r ( \widetilde{a}^{i} )  ( 1\otimes w_{\phi } )
=\Det(\overline{\eta }^{i})(\overline{\phi })^{\phi (1)}\tbigwedge
r ( a^{i} )  ( 1\otimes w_{\phi } ) . 
\]
As the $\eta _{jl}^{i}$ have rational coefficients $\Det ( \overline{\eta 
}^{i} ) (\overline{\phi })=\Det(\eta ^{i})(\phi )$, and so the
homomorphism representing the class only changes by the homomorphism which
maps $\phi $ to 
\[
\tprod_{i}{\text{Det}}(\eta ^{i})(\phi )^{(-1)^{i}}\times \tprod_{i}\left| 
{\text{Det}}(\eta ^{i})(\phi )^{(-1)^{i}}\right| ; 
\]
and again this comes from an element of the denominator of (2).\medskip

(iii) If $\left\{ \widetilde{w}_{\phi ,k}\right\} $ is a further orthonormal
basis of $W$, then the wedge product $\tbigwedge r\left( a^{i}\right) \left(
1\otimes \widetilde{w}_{\phi }\right) $ differs from $\tbigwedge r\left(
a^{i}\right) \left( 1\otimes w_{\phi }\right) $ by a power of the
determinant of a unitary base-change, which therefore has absolute value
1.
\bigskip

The following two properties of arithmetic classes follow readily from the
definition.
\medskip

\no{\bf Lemma 3.4} \ {\sl $(P^{\bullet },p_{\bullet }),\;(Q^{\bullet
},q_{\bullet })$  be perfect metrised $\Z[G]$-complexes
and endow the complex $P^{\bullet }\oplus Q^{\bullet }$ with
metrics $p_{\phi }\otimes q_{\phi }$  on the equivariant
determinants of cohomology via the identification 
\[
\det \left( \Hh^{\bullet }\left( P_{\phi }^{\bullet }\oplus Q_{\phi
}^{\bullet }\right) \right) =\det \left( \Hh^{\bullet }\left( P_{\phi
}^{\bullet }\right) \right) \otimes \det \left( \Hh^{\bullet }\left(
Q_{\phi }^{\bullet }\right) \right) . 
\]
Then}
\[
\chi \left( P^{\bullet }\oplus Q^{\bullet },p_{\phi }\otimes q_{\phi
}\right) =\chi (P^{\bullet },p_{\bullet })\chi (Q^{\bullet },q_{\bullet }). 
\]
\smallskip

\no{\sc Proof.}  This follows on choosing bases for $P^{\bullet }$ and $Q^{\bullet }\;$%
and then using these bases to form a basis of $P^{\bullet }\oplus Q^{\bullet
}$.\ \ \ \ $\Box \medskip $

Recall that $\left| -\right|$ denotes the standard metric on $\Bbb{C}$.
\medskip

\no{\bf Lemma 3.5} \ {\sl If $P^{\bullet }$  is an acyclic 
perfect metrised $\Z[G]$-complex and if we endow each complex
line
\[
\det \left( \Hh^{\bullet }\left( P_{\phi }^{\bullet }\right) \right)
=\det \left( \left\{ 0\right\} \right) =\Bbb{C} 
\]
with the metric $\left| -\right| $, then $\chi (P^{\bullet },\left|
-\right| _{\bullet })=1$.}
\medskip

\no{\sc Proof.} As $P^{\bullet }$ is acyclic and its terms are projective, it is
isomorphic to a complex 
\[
\cdots \rightarrow W^{i-1}\oplus W^{i}\rightarrow W^{i}\oplus W^{i+1}\rightarrow
\cdots
\]
where the $W^{i\text{ }}$ are all projective and where the boundary maps are
projection to the second factor. Using bases of the $W^{i}$ to form bases of
the $P^{i}$, together with the standard properties of $\det ,$ we see that
the products in 3.3 all telescope to 1.\ \ \ $\Box \medskip $
\medskip

\no{\bf Lemma 3.6} \ {\sl If $p_{\bullet }$ and $q_{\bullet }$
are two sets of metrics on the equivariant determinants of
cohomology of $P^{\bullet }$, then for each $\phi \in \widehat{G}
,\;p_{\phi }=\alpha \left( \phi \right) ^{\phi \left( 1\right) }q_{\phi }$%
 for a unique positive real number $\alpha \left( \phi \right) .\,$ The class $\chi \left( P^{\bullet },p_{\bullet }\right) \chi \left(
P^{\bullet },q_{\bullet }\right) ^{-1}$ in $A(\Z[G])$ 
is represented by the homomorphism which maps each $\phi \in \widehat{G}$
 to the value $1\times \alpha \left( \phi \right)$.}
\medskip

\no{\sc Proof.} This follows immediately from (3.3).\ \ \ $\Box$

\medskip

\no{\bf 3.C \ Invariance under quasi-isomorphism.}
\medskip

Let $\left( C^{\bullet },c_{\bullet }\right) $ and $\left( D^{\bullet
},d_{\bullet }\right) $ denote bounded (not necessarily perfect) metrised $\Z[G]$-complexes and suppose that there is a $\Z[G]$-cochain map
$\alpha :C^{\bullet }\rightarrow D^{\bullet }$. Recall that  $\alpha $
is called a \textit{quasi-isomorphism} if it induces an isomorphism on the
cohomology of the complexes. Theorem 2 in [KM] implies that if $\alpha $ is
a quasi-isomorphism, then it induces natural isomorphisms 
\[
\det \left( \Hh(\alpha _{\phi })\right) :\det \left( \Hh^{\bullet
}\left( C_{\phi }^{\bullet }\right) \right) \cong \det \left( \Hh
^{\bullet }\left( D_{\phi }^{\bullet }\right) \right) 
\]
so that the following square commutes: 
\[
\begin{array}{lll}
\det ( C_{\phi }^{\bullet } ) & \rightarrow & \det ( D_{\phi
}^{\bullet }) \\ 
\ \ \ \downarrow &  & \ \ \downarrow \\ 
\det  ( \Hh^{\bullet } ( C_{\phi }^{\bullet } )  ) & 
\rightarrow & \det  ( \Hh^{\bullet } ( D_{\phi }^{\bullet
} )  )
\end{array}
\]
where the top horizontal map is $\det (\alpha _{\phi })$ and where the
vertical isomorphisms are $\xi _{C,\phi }$ and $\xi _{D,\phi }$ of
(1).
\medskip

\no{\bf Definition 3.7} \ A quasi-isomorphism $\alpha :C^{\bullet }\rightarrow
D^{\bullet }$ is called a \textit{metric} \textit{quasi-isomorphism} from $%
\left( C^{\bullet },c_{\bullet }\right) $ to $\left( D^{\bullet },d_{\bullet
}\right) $ if $c_{\phi }=d_{\phi }\circ \det  ( \Hh(\alpha _{\phi
}) )$ for each $\phi  \in \widehat{G}.$
\medskip

The following result is an immediate consequence of the definitions:
\medskip

\no{\bf Lemma 3.8} \ {\sl Suppose again that  $\alpha :C^{\bullet }$ $%
\rightarrow D^{\bullet }$  is a quasi-isomorphic cochain map and
that metrics $d_{\phi }$  are given on the $\det ( \Hh^{\bullet }( D_{\phi }^{\bullet }) )$.
 Then there
is a unique set of metrics $c_{\phi }$ on $\det (\Hh
^{\bullet } ( C_{\phi }^{\bullet } )  ) \;$ such that $\alpha :(C^{\bullet },c_{\bullet })\rightarrow (D^{\bullet },d_{\bullet })$
 is a metric quasi-isomorphism; we call the metrics $c_{\bullet }$
the metrics on the equivariant determinants of cohomology induced
from $d_{\bullet }$ via $\alpha $.
If $\beta :C^{\bullet }\rightarrow D^{\bullet }$ is a
further quasi-isomorphic cochain map and if
$\Hh^{i}(\alpha )=\Hh^{i}(\beta )$ for all $i,$ 
then $\det  ( \Hh ( \alpha _{\phi } )  ) =\det
 ( \Hh ( \beta_{\phi } )  ) $  for all $\phi
\in \widehat{G}$ and so $\alpha $  and $\beta $ 
induce the same metrics on the equivariant determinant of cohomology of $C^{\bullet }$.
\medskip

The main result of this sub-section is 
\smallskip 
\medskip

\no{\bf Definition-Theorem 3.9} \  {\sl With the above notation and
hypotheses, let $\alpha :\left( C^{\bullet },c_{\bullet }\right)
\rightarrow \left( D^{\bullet },d_{\bullet }\right) $  be a metric
quasi-isomorphism and suppose further that we can find perfect metrised $\Z[G]$-complexes $ ( P^{\bullet },p_{\bullet } )$, 
  resp. $\left( Q^{\bullet },q_{\bullet }\right) \;$ which
support metric quasi-isomorphisms $f:$ $\left( P^{\bullet
},p_{\bullet }\right) \rightarrow \left( C^{\bullet },c_{\bullet }\right) ,$
resp. $g: (Q^\bullet, q_\bullet)\rightarrow \left( D^{\bullet },d_{\bullet }\right) $. Then $\chi \left(
P^{\bullet },p_{\bullet }\right) =\chi \left( Q^{\bullet },q_{\bullet
}\right)$.
 
In particular: for a metrised }$\Z[G]$-complex $\left(
C^{\bullet },c_{\bullet }\right) $  with the property that $
C^{\bullet }$  is quasi-isomorphic to a perfect complex $P^{\bullet
},$ we let $p_{\bullet }$ denote the metrics on the
equivariant determinant of cohomology of $P^{\bullet }$ induced by 
$c_{\bullet }$; then we can unambiguously define the arithmetic
class of $\left( C^{\bullet },c_{\bullet }\right) $ to be the
class $\chi \left( P^{\bullet },p_{\bullet }\right) ;$ this class
depends only on $\left( C^{\bullet },c_{\bullet }\right) $ and not
on the particular choice of perfect complex $P^{\bullet }$. Thus
with this definition we have the equality $\chi \left( C^{\bullet
},c_{\bullet }\right) =\chi \left( D^{\bullet },d_{\bullet }\right) $.}
\medskip 

Before proving the theorem we first need some preliminary results.
\medskip

\no{\bf Lemma 3.10} \ {\sl Given maps of $\Z[G]$ -complexes $
M^{\bullet }\stackrel{\varphi }{\rightarrow }L^{\bullet }\stackrel{\pi }{%
\leftarrow }N^{\bullet }$\textit{\ with }$\pi $  a surjective
quasi-isomorphism and with $M^{\bullet }\,$ perfect, there exists a 
$\Z[G]$ -cochain map $\psi :M^{\bullet }\rightarrow N^{\bullet
}\; $ such that $\pi \circ \psi =\varphi$.}
\medskip

\smallskip

\no{\sc Proof.} See VI.8.17 in [Mi].\medskip\ \ \ $\Box $

\medskip

\no{\bf Corollary 3.11} \   {\sl If $0\rightarrow A^{\bullet }\stackrel{%
\alpha }{\rightarrow }B^{\bullet }\stackrel{\beta }{\rightarrow }C^{\bullet
}\rightarrow 0$ is an exact sequence of perfect $\Z[G]$ 
-complexes and if $A^{\bullet }$ is acyclic, then there exists a
cochain map  $i:C^{\bullet }\rightarrow B^{\bullet }$
 which is a section of $\beta$.}
\medskip

\no{\sc Proof.} Apply the above lemma to $C^{\bullet }=C^{\bullet }\stackrel{\beta }{%
\leftarrow }B^{\bullet }.\;\;\Box$.
\medskip

\no{\sc Proof of theorem.} First we choose an acyclic perfect complex $L^{\bullet }$
and a map $\lambda :L^{\bullet }\rightarrow D^{\bullet }$ such that $\lambda
\oplus g$ is surjective. We then endow the equivariant determinants of the
cohomology of $L^{\bullet }$ with the trivial metrics $l_{\bullet }$ as
per Lemma 3.5. Then by 3.4 and 3.5  
\[
\chi \left( L^{\bullet }\oplus Q^{\bullet },l_{\bullet }q_{\bullet }\right)
=\chi \left( Q^{\bullet },q_{\bullet }\right) . 
\]
Thus, without loss of generality, we may now assume that $g$ is surjective.

Consider the diagram 
\[
P^{\bullet }\stackrel{f}{\rightarrow }C^{\bullet }\stackrel{\alpha }{%
\rightarrow }D^{\bullet }\stackrel{g}{\leftarrow }Q^{\bullet }. 
\]
By Lemma 3.10 we can find a $\Z[G]$-map $\beta :P^{\bullet }\rightarrow
Q^{\bullet }$ such that $\alpha \circ f=g\circ \beta .$ As $f,g$ and $\alpha 
$ are all quasi-isomorphisms, $\beta $ is also a quasi-isomorphism.

As previously, by adding an acyclic complex with trivial metrics $\left(
L^{\bullet },l_{\bullet }\right) \;$to $\left( P^{\bullet },p_{\bullet
}\right) $ , setting $P^{\prime \bullet }=L^{\bullet }\oplus P^{\bullet }$
and$\;p_{\bullet }^{^{\prime }}=l_{\bullet }p_{\bullet },\;$we obtain a
surjective quasi-isomorphism $\beta ^{\prime }:P^{\prime \bullet
}\rightarrow Q^{\bullet }$ and 
\[
\chi \left( P^{\prime \bullet },p_{\bullet }^{\prime }\right) =\chi \left(
P^{\prime \bullet },l_{\bullet }p_{\bullet }\right) =\chi \left( P^{\bullet
},p_{\bullet }\right) . 
\]
$\;$We let $f^{\prime }:P^{\prime \bullet }\rightarrow C^{\bullet }$ denote
the composition of $f$ with the natural projection map. Then in general of
course it will not be true that $\alpha \circ f^{\prime }=g\circ \beta
^{\prime };$ however, as $L^{\bullet }$ is acyclic,$\;$we do know that $%
\alpha \circ f^{\prime }$\ and $g\circ \beta ^{\prime }$ agree on
cohomology, i.e. $\Hh^{i}\left( \alpha \circ f^{\prime }\right) =\Hh
^{i}\left( g\circ \beta ^{\prime }\right) $ for all $i.$

In order to complete the proof of Theorem 3.9, we apply Corollary 3.11 to
choose a section $\gamma :Q^{\bullet \text{ }}\rightarrow P^{\prime \bullet
} $ of $\beta ^{\prime }$. Again as per Lemma 3.5 we endow the equivariant
determinants of cohomology of $\ker \beta ^{\prime }$ with the trivial
metric $s_{\bullet }$ ; as per Lemma 3.4 we endow $P^{\prime \bullet }$ with
the metric $\widetilde{q}_{\bullet }$ given by $s_{\bullet }.\gamma _{\ast
}q_{\bullet }$. Then 
\[
\chi \left( P^{\prime \bullet },\widetilde{q}^{\bullet }\right) =\chi \left(
\ker \beta ^{\prime },s_{\bullet }\right) \chi \left( \gamma Q^{\bullet
},\gamma _{\ast }q_{\bullet }\right) =\chi \left( \gamma Q^{\bullet },\gamma
_{\ast }q_{\bullet }\right) =\chi \left( Q^{\bullet },q_{\bullet }\right) . 
\]
However, as the metrics $q_{\bullet },\;$on the equivariant determinants of
the cohomology of $Q^{\bullet },$ are induced from $d_{\bullet }$ via $H(g)$%
,\ the metrics \ $\widetilde{q}_{\bullet }$ are the transport to $P^{\prime
\bullet }$ of the metrics $d_{\bullet }$ via $\Hh(g\circ \beta ^{\prime })=\Hh
(\alpha \circ f^{\prime })$. Thus $p_{\bullet }^{\prime }\;$and $\widetilde{q%
}_{\bullet }\;$are both transports of the $d_{\bullet }$ via $\Hh(g\circ \beta
^{\prime })=H(\alpha \circ f^{\prime })$, and so by Lemma 3.8 they are
equal. Therefore we have shown$\;$ 
\[
\chi \left( P^{\bullet },p_{\bullet }\right) =\chi \left( P^{\prime \bullet
},p_{\bullet }^{\prime }\right) =\chi \left( P^{\prime \bullet },\widetilde{q%
}_{\bullet }\right) =\chi \left( Q^{\bullet },q_{\bullet }\right) 
\]
which is the desired result. \ \ $\Box$
\bigskip

\section{Arithmetic Classgroups}

\no{\bf 4.A \ Symplectic arithmetic classes.}
\medskip

The arithmetic classgroup $A(\Z[G])$ carries a great deal of
information. In consequence, it is often advantageous in practice to work
with various image groups. The most important of these is the \textit{%
symplectic} arithmetic classgroup.

Recall that by the Hasse-Schilling norm theorem 
\begin{equation}
\Det( \Q[G]^{\times } ) =\Hom_{\Omega }^{+} (
R_{G},\overline{\Bbb{Q}}^{\times } )
\end{equation}
where the right-hand expression denotes Galois equivariant homomorphisms
whose values on $R_{G}^{s},$ the group of virtual symplectic characters, are
all totally positive. By analogy with the map $\Delta $ of 3.A, we again
have a diagonal map 
\[
\Delta ^{s}:\Hom_{\Omega }^{+} ( R_{G}^{s},\overline{\Bbb{Q}}%
^{\times } ) \rightarrow \Hom_{\Omega } (
R_{G}^{s}, J_{f} ) \times \Hom ( R_{G}^{s},\R_{>0} ) 
\]
where $\Delta ^{s}\left( f\right) =f\times \left| f\right| =f\times
f$.
\medskip

\no{\bf Definition 4.1} \ The group of \textit{symplectic arithmetic classes} $%
A^{s}(\Z[G]) $ is defined to be the quotient group 

\[
A^{s} ( \Z[G])) =\frac{\Hom_{\Omega } (
R_{G}^{s}, J_{f} ) \times \Hom ( R_{G}^{s}, \R_{>0} ) 
}{\func{Im}\Delta ^{s}\cdot\ ( \text{Det}^{s}\ ( \widehat{\Z}[G]^{\times } ) \times 1 ) } 
\]
where $\Det^{s} ( \widehat{\Z}[G]^{\times } ) $ denotes the
restriction of $\Det ( \widehat{\Z}[G]^{\times } ) $ to $%
R_{G}^{s} $. In general, given a homomorphism $f$ on $R_{G},$ we shall write 
$f^{s}$ for the restriction of $f$ to $R_{G}^{s}$. Clearly, restriction from 
$R_{G}$ to $R_{G}^{s}$ induces a homomorphism 
\[
\rho : A( \Z[G]) \rightarrow A^{s} ( \Z[G]) . 
\]
\smallskip

\no{\bf 4.B \ Torsion classes.}
\medskip

Let K$_{0}$T$(\Z[G]) $ denote the Grothendieck group of
finite, cohomologically trivial $\Z[G]$-modules and let K$_{0}$T$(
\Z_{p}[G]) $ denote the Grothendieck group of finite,
cohomologically trivial $\Z_p[G]$-modules. Thus the decomposition of a
finite module into its $p$-primary parts induces the direct sum
decomposition 
\[
{\rm K}_{0}{\rm T} ( \Z[G] ) =\oplus _{p}{\rm K}_{0}{\rm T}
( \Z_{p}[G]) . 
\]
We write K$_{0} ( {\bf F}_{p}[G]) $ for the Grothendieck group of
finitely generated projective ${\Bbb{F}}_{p}[G]$-modules; since each such module
may be considered as a finite, cohomologically trivial $\Z_p[G]$%
-module, we have a natural map 
\[
{\rm K}_{0}\left( {\Bbb{F}}_{p}[G]\right) \rightarrow {\rm K}_{0}{\rm T}
 ( \Z_p[G] ) . 
\]
>From Chapter 1, Theorem 3.3 in [T2] recall that there is the Fr\"{o}hlich
isomorphism 
\[
{\rm K}_0{\rm T}(\Z[G] ) \cong \frac{\Hom_{\Omega
} ( R_{G},J_{f} ) }{\Det( \widehat{\Bbb{Z}}[G]^{\times
}) }; 
\]
thus there is a natural map $\nu :{\rm K}_0{\rm T}(\Z[G])
\rightarrow A ( \Z[G] )$, induced by $f\mapsto f\times 1$ for $%
f\in \;\Hom_{\Omega } ( R_{G},J_{f} ) .$

In order that our invariants agree with the standard invariants in Arakelov
theory, our convention here is that of I.3.2 in [T2]: namely, if 
$M=\Z_p[G]/\alpha \Z_p[G]$ is a 
$\Z_{p}$-torsion $\Z_p[G]$-module,
then the class of $M$ in K$_{0}$T$\left( \Z[G]\right) $ is represented by
$\Det( \alpha  ) ;$ this then is the inverse of the description
given in 4.4 in [C]. It will be important in the sequel to keep this in\
mind when performing various torsion calculations in Sections 7 and 8.
\medskip

\no{\bf 4.C \ Tame arithmetic classes.}
\medskip

Although we shall ultimately always be interested in forming arithmetic
classes over the integral group ring $\Z[G]$, in carrying out
calculations it will often be advantageous to work with more general group
rings, where we allow tame coefficients. With this in mind, we let $T$
denote the maximal abelian tame extension of $\Bbb{Q}$ in $\overline{\Bbb{Q}}$
and we set 
\[
{\Det}^{s} ( \widehat{\Oo_T}[G]^{\times } ) =\lim_{%
\overrightarrow{\scriptsize L}}{\Det}^{s} ( \widehat{\Oo_L}[G]^{\times } ) 
\]
where the direct limit extends over all finite extensions $L$ of $\Bbb{Q}$
in $T$ and where $\widehat{\Oo_L}$ is the ring of integral adeles $\widehat{%
\Z}\otimes \Oo_{L}$.

Next we define the group $\Hom_{\Omega }^{+} ( R_{G}^{s},\overline{\Bbb{Q}%
}^{\times } ) $ to be the subgroup of $\Hom_{\Omega } ( R_{G}^{s},
\overline{\Bbb{Q}}^{\times } ) $ whose values are all totally positive
(note that the values of the homomorphisms in the latter group are of course
all real by Galois equivariance); we then define $\Delta ^{\prime }$ to be
the diagonal homomorphism 
\[
{\Hom}_{\Omega }^{+} ( R_{G}^{s},\overline{\Bbb{Q}}^{\times } )
\rightarrow {\Hom}_{\Omega } ( R_{G}^{s},J_{f} ) \times {\rm
Hom} ( R_{G}^{s},\R_{>0} ) . 
\]
In arithmetic calculations we shall often need to work with the \textit{tame
symplectic arithmetic classgroup }defined as 
\begin{equation}
A_{T}^{s}\left( \Z[G]\right) =\frac{\left({\Det}^{s} ( \widehat{\Oo_T}G^{\times } )\cdot {\Hom}_{\Omega } ( R_{G}^{s},J_{f} )
\right) \times {\Hom} ( R_{G}^{s},\R_{>0} ) }{{\rm Im}%
(\Delta')\cdot \left( {\Det}^{s} ( \widehat{\Oo_T}[G]^{\times
} ) \times 1\right) }.
\end{equation}
Inclusion then induces a homomorphism 
\[
\eta :A^{s}\left( \Z[G]\right) \rightarrow A_{T}^{s}\left(\Z[G]\right) . 
\]
For a perfect metrised $\Z[G]$-complex $ ( P^{\bullet },p_{\bullet
} ) ,$ we write $\chi ^{s} ( P^{\bullet },p_{\bullet } ) $ for
the image of $\chi  ( P^{\bullet },p_{\bullet } ) \;$in $%
A_{T}^{s} ( \Z[G]) $. 
\medskip

\no {\bf 4.D \ Rational classes.}
\medskip

Rational classes are ubiquitous in arithmetic applications. The subgroup of 
\textit{rational symplectic arithmetic classes} is defined to be the
subgroup of $A_{T}^{s}\left( \Z[G]\right) $ generated by Hom$_{\Omega
} ( R_{G}^{s},{\Q}^{\times } ) \times 1$, that is to say 
\[
R^{s} ( \Z[G] ) =\frac{\left( {\Hom}_{\Omega } (
R_{G}^{s},\Q^{\times } )\cdot {\Det}^{s}( \widehat{\Oo_T}%
[G]^{\times } ) \times 1\right) \cdot{\func{Im}}(\Delta')}{\left( \Det^{s} ( \widehat{\Oo_T}[G]^{\times } ) \times 1\right) \cdot \func{Im}(
\Delta')}. 
\]
The natural map $\overline{\Bbb{Q}}^{\times }\hookrightarrow J_{f}$ induces
a map 
\[
\theta ': {\func{Im}}(\Delta')\cdot \left( {\Hom}_{\Omega
} ( R_{G}^{s},J_{f} ) \times 1\right) \rightarrow {\Hom}
_{\Omega } ( R_{G}^{s},J_{f} ) 
\]
which is defined as follows: consider $h\in {\func{Im}}(\Delta')\cdot \left( {\Hom}_{\Omega
} ( R_{G}^{s},J_{f} ) \times 1\right)$
and let $h_{f}$ resp. $h_{\infty }$ denote the finite resp. archimedean
component of $h$. Then $ h_{\infty }$ determines a unique element $%
h_{\infty }'$ of ${\func{Im}}(\Delta');$ we define $\theta' ( h ) =h_{f}h_{\infty }^{\prime -1}$.
Clearly $\theta'$ vanishes on ${\func{Im}}(\Delta ')$ and so induces a
homomorphism 
\[
\theta :R^{s}\left( \Z[G]\right) \rightarrow \frac{{\Hom}_{\Omega
} ( R_{G}^{s},\Q^{\times } )           {\Det}^{s} ( \widehat{\Oo_T}[G]^{\times } ) }
{{\Det}^{s} ( \widehat{\Oo_T}[G]^{\times
} ) }. 
\]
>From [CNT] (see also Corollary 3 to Theorem 17 in [F2]) we know that 
\begin{equation}
{\Hom}_{\Omega } ( R_{G}^{s},\Q^{\times } ) \cap {\Det}
^{s} ( \widehat{\Oo_T}G^{\times } ) = \{ 1 \}
\end{equation}
and so by (9) we see that $\theta $ may be written as an isomorphism 
\[
\theta :R^{s}\left( \Z[G]\right) \rightarrow {\Hom}_{\Omega }\left(
R_{G}^{s},\Q^{\times }\right) . 
\]
\medskip

\no{\bf 4.E \ Passage to degree zero.}
\medskip

In this sub-section we describe a useful procedure for changing arithmetic
classes by passage to characters of degree zero. In practice this will allow
us to disregard various free classes which arise in our calculations. For an
abelian group $A$ and for $f\in \;$Hom$\left( R_{G},A\right) ,$ we write $%
\widetilde{f}\in \;$Hom$\left( R_{G},A\right) $ for the $\hom $omorphism
defined by the rule $\widetilde{f}\left( \chi \right) =f\left( \chi -\chi
\left( 1\right) 1_{G}\right) $, where $1_{G}$ denotes the trivial character
of $G.$ Note that for $z\in \Z_p[G]^{\times }$%
\[
\widetilde{{\Det}} ( z ) ={\Det} ( zd^{-1} )\ \ \ \hbox{where\ \ } d={\Det} ( z )  ( 1_{G} ) , 
\]
and so $\widetilde{\Det( \Z_p[G]^{\times } ) }\subset \Det
\left( \Z_p[G]^{\times }\right) ;$ similarly $\widetilde{{\text{Im}}%
 ( \Delta  ) }\subset {\rm Im} ( \Delta  )$. Thus for a class 
${\frak{c}}\in A\left( \Z[G]\right)$, represented by a homomorphism $f$
under (2), we can unambiguously define a new class $\widetilde{\frak{c}}$,
depending only on $\frak{c}$,\ to be that class represented by the
homomorphism $\widetilde{f}.$

\medskip The map ${\frak{c}}\mapsto \widetilde{\frak{c}}$ can be interpreted
in the following fashion in terms of $G$-fixed points together with the
induction map 
\[
\text{Ind}:A\left( \Z\right) \rightarrow A\left( \Z[G]\right) 
\]
given in terms of character maps by Ind$\left( f\right) \left( \psi \right)
=f ( {\rm Res}_{G}^{\left\{ 1\right\} }\psi  ) =f ( \psi
\left( 1\right) .1_{\left\{ 1\right\} } ) $ for $\psi \in R_{G}.$
\medskip

\no{\bf Lemma 4.2} \  {\sl With the notation of 3.B , let ${\frak{c}
}=\chi \left( P^{\bullet },p_{\bullet }\right) \in A\left( \Z[G]\right) \ 
$ and let ${\frak{c}}_{0}$  be the class in $A\left( \Z
\right) $  of $\left( P^{\bullet G},p_{1}\right)$, where 
$P^{\bullet G}$  denotes the complex obtained from $P^{\bullet }$
by taking $G$-fixed points and where $p_{1}$ 
denotes the metric on the determinant of the cohomology of $P^{\bullet G}$,
  obtained by identifying $P_{\Bbb{C}}^{\bullet G}$  with
the isotypic component of the $P_{\Bbb{C}}^{\bullet }$ for the
trivial character of $G$. There is then an equality $\widetilde{%
\frak{c}}={\frak{c}}\cdot {\rm Ind}({\frak{c}}_{0})^{-1}$  
in $A\left(\Z[G]\right)$.}

\medskip

\no{\sc Proof.} Let $\Sigma =\sum_{g\in G}g$. As each term of $P^{\bullet }$ is
projective, $P^{\bullet G}=\Sigma P^{\bullet }.$ We adopt the notation of
3.B and assume that $f$ is the representative character map for the class $%
{\frak{c}}=\chi  ( P^{\bullet },p_{\bullet } ) $ obtained by using
local bases $ \{ \alpha _{p}^{ij} \} , \{ a^{ij} \}$.
Then we let $h$ denote the representative for the class ${\frak{c}}_{0}=\chi
 ( P^{\bullet G},p_{1} ) $ obtained by using local bases $ \{
\Sigma \alpha _{p}^{ij} \} , \{ \Sigma a^{ij} \} .$ To prove
the lemma it will then suffice to show that $f\left( 1_{G}\right) =h (
1_{\left\{ 1\right\} } ) .$

We start by considering the non-archimedean coordinates. With the notation
of 3.B we have $ ( a^{ij} ) _{j}=\lambda _{p}^{i} ( \alpha
_{p}^{ij} ) _{j},$ and so $ ( \Sigma a^{ij} ) _{j}=e\cdot \lambda
_{p}^{i} ( \Sigma \alpha _{p}^{ij} ) _{j}$ where $e=\Sigma /\left|
G\right| $ is the idempotent associated to the trivial character of $G$.
Since  
$\det ( e\lambda _{p}^{i} ) =\Det( \lambda _{p}^{i} )
 ( 1_{G} )$, we conclude that the non-archimedean coordinates of $%
f $ and $h$ are equal.

To conclude we consider the archimedean coordinates. As $e$ resp. $1$ is a
basis element of length 1 for the trivial isotypic component of $\C[G]$
resp. $\C$ with respect to $\nu _{G}$ resp. $\nu _{1}$ (see 2.3 for
the definition of $\nu$), then as in 3.B we see that the archimedean
coordinate of $f ( 1_{G} ) $ resp. $h ( 1_{ \{ 1 \}
} ) $ is obtained by evaluating $p_{1}$ on the wedge product $\bigwedge
\alpha _{G}\left(  r_{G} ( a^{ij}\otimes e ) \right) ^{ (
-1 ) ^{i}}$ resp.  
$\bigwedge \alpha _{1} ( \Sigma a^{ij}\otimes 1 ) ^{ ( -1 )
^{i}}$ (see 2.3 to recall the definition of $\alpha$).
 Since 
\[
\alpha _{G}\left( r_{G} ( a^{ij}\otimes e ) \right) =\alpha
_{G} ( \Sigma a^{ij}\otimes e ) =\Sigma a^{ij}=\alpha _{1} (
\Sigma a^{ij}\otimes 1) 
\]
it follows that the archimedean coordinates of $f\left( 1_{G}\right) $ and $%
h ( 1_{ \{ 1 \} } ) $ are also equal.\ \ \ $\Box $

\bigskip

\section{Arithmetic applications}

\no{\bf 5.A \ Preliminary results.}
\medskip

Let $\mathcal{X}$ be a projective scheme over $\Spec(\Z)$
with structure morphism 
$f:\X\rightarrow \Spec(\Z)$. Suppose further
that $\mathcal{X}$ is flat over $\Spec(\Z)$ with
equidimensional fibres of dimension $d$ and that the generic fibre of 
$\mathcal{X}$ is smooth. For the sake of brevity, in the sequel we shall
refer to $\mathcal{X}$ simply as an arithmetic variety. Suppose further that 
$\mathcal{X}$ is endowed with an action $(\X, G) $ by a
given finite group $G$; we shall suppose that the action is\textit{\ tame},
in the sense that for each closed point $x$ of $\mathcal{X}$ the inertia
group of $x$ has order coprime to the residue characteristic of $x$. Since $%
\mathcal{X}$ is projective, the quotient scheme ${\mathcal{Y}}=\mathcal{X}$/$%
G$ is defined and we denote the quotient morphism by $\pi :\mathcal{%
X\rightarrow Y}$. Let $b$ denote the branch locus on $\mathcal{Y}$ of the
cover $\mathcal{X}$/$\mathcal{Y}$. From now on we shall suppose that $%
\mathcal{Y}$ is connected and that the branch locus $b$ is a Cartier divisor
on $\mathcal{Y}$ with strictly normal crossings. Since $G$ acts tamely on $%
\mathcal{X}$, we note that by the valuative criterion for properness it
follows that $G$ must act freely on the generic fiber $\X_\Q$ (see
1.2.4(d) in [CEPT1]).

We now consider the construction of arithmetic classes for complexes of
sheaves on $\mathcal{X}$. For a detailed account of the formation of Euler
characteristics (without metrics) associated to a tame action, the reader is
referred to [CEPT4]. Let $\mathcal{F}^{\bullet }\;$denote a bounded complex
of coherent $G$-$\X$ sheaves. Consider a $G$-stable open affine cover 
$\mathcal{U}$ of $\mathcal{X}$ and take the chain complex $C^{\bullet }$
which is the associated simple complex to the double complex $C^{\bullet
} (\mathcal{U},\mathcal{F}^{\bullet } ) $. There is an isomorphism
on the derived category between $C^{\bullet }$ and $\Rr\Gamma  ( \mathcal{X%
},\mathcal{F}^{\bullet } ) $ which induces isomorphisms 
\[
\det  ( {\Hh}^{\bullet } ( \Rr\Gamma ( \X, \F%
^{\bullet } )  ) _{\phi } ) \cong \det  ( \Hh%
^{\bullet } ( C^{\bullet } ) _{\phi } ) \;\hbox{\rm \ for all\  }\phi
\in \widehat{G}. 
\]
\medskip

\no{\bf Lemma 5.1} \ {\sl For $C^{\bullet }$ as above, there is
a perfect $\Z[G]$-complex $P^{\bullet }$ with a
quasi-isomorphism $\gamma :P^{\bullet }\rightarrow C^{\bullet }$.}
\medskip

\no{\sc Proof.} For full details we refer to the proof of Theorem 1.1 in [C]; so we
shall now briefly only sketch the proof for the reader's convenience. From
Lemma III.12.3 in [H] we may construct a quasi-isomorphism $\gamma
_{1}:P_{1}^{\bullet }\rightarrow C^{\bullet }$ where the complex $%
P_{1}^{\bullet }$ is a bounded complex of finitely generated 
$\Z[G]$-modules all of whose terms except the initial term, $P_{1}^{N}$ say, are
free $\Z[G]$-modules. Since the mapping cylinder of $\alpha \circ \gamma
_{1}$ is acyclic with all terms, except possibly $P_{1}^{N}$, being
cohomologically trivial $\Z[G]$-modules, we therefore deduce that $%
P_{1}^{N}$ is a cohomologically trivial $\Z[G]$-module, and it may
therefore be written as the quotient of two projective $\Z[G]$-modules;
replacing $P_{1}^{N}$ by this perfect complex of length 2\ provides $%
P^{\bullet }$ and $\gamma .\;\;$\ $\Box$
\medskip

\no{\bf Definition 5.2} \ Suppose now that we are given metrics $h_{\phi }$ on
the $\det( \Hh^{\bullet } ( \Rr\Gamma  ( \X, \F^{\bullet } )  ) _{\phi } ) \;$for all $%
\phi \in \widehat{G}$. These metrics then induce metrics $p_{\phi }$ on $%
\det  ( \Hh^{\bullet } ( P_{\phi }^{\bullet } )  ) $
and by Theorem 3.9 we know that the arithmetic class $\chi  ( P^{\bullet
},p_{\bullet } )$ is independent of choices; we denote this class 
\[
\chi  \left( \Rr\Gamma  ( \X, \F^{\bullet
} ) ,h_{\bullet } \right) 
\]
and the image of this class in the symplectic arithmetic classgroup $%
A^{s}(\Z[G]) $ will be denoted $\chi ^{s} (\Rr
\Gamma  ( \X, \F^{\bullet } ) ,h_{\bullet
} ) $.  
\medskip

The following results describe some basic properties of such arithmetic
classes. The first two results follow immediately from\ 3.4 and 3.6.
\medskip

\no{\bf Proposition 5.3} \ {\sl Let $\F^\bullet$, $\G^\bullet$ 
be bounded complexes of coherent $G$-$\X$
 sheaves; let $h_{\bullet },\;$ resp. $g_{\bullet }$
 be metrics on the equivariant determinants of cohomology of $
\Rr\Gamma  ( \X, \F^{\bullet } ) $, resp. $\Rr\Gamma (\X, \G^{\bullet } )$. Then
\[
\chi (\Rr\Gamma  (\X, \F^{\bullet
}\oplus \G^{\bullet } ) ,h_{\bullet }g_{\bullet } ) =\chi
 ( \Rr\Gamma ( \X, \F^{\bullet
} ) ,h_{\bullet } ) \cdot \chi  ( \Rr\Gamma  ( \X, \G^{\bullet } ) , g_{\bullet } ) . 
\]}
\smallskip

\no{\bf Proposition 5.4} \ {\sl Let $j_{\bullet }$  denote a
further set of metrics on the equivariant determinants of cohomology of 
${\rm R}\Gamma ( \X, \F^{\bullet } ) $
 and suppose that for each $\phi \in \widehat{G},\;h_{\phi
}=\alpha ( \phi ) ^{\phi ( 1) }j_{\phi }\;$
for $\alpha \left( \phi \right) \in \R_{>0}$. Then the
hermitian class  
\[
\chi  ( \Rr\Gamma  ( \X, \F^{\bullet
} ) ,h )\cdot \chi  ( \Rr\Gamma  ( \X, \F^{\bullet } ) ,j ) ^{-1} 
\]
 is represented by the homomorphism which maps each $\phi \in 
\widehat{G}$ to
\[
\varphi \longmapsto 1\times \alpha \left( \phi \right) . 
\]}
\smallskip

\no{\bf Proposition 5.5} \ {\sl If $0\rightarrow \F
\rightarrow \G\rightarrow {\cal H}\rightarrow 0$  is an
exact sequence of coherent $G$-$\X$  sheaves with metrics $
f_{\bullet }, g_{\bullet }, h_{\bullet }, $ on their equivariant
determinants of cohomology, with the property that $f_{\phi }\otimes
h_{\phi }=g_{\phi }$  under the isomorphisms  
\[
\det  ( \Hh^{\bullet } ( \F ) _{\phi } )
\otimes \det  ( \Hh^{\bullet } ( {\cal H} ) _{\phi
} ) \cong \det  ( \Hh^{\bullet } ( \G )
_{\phi } ) 
\]
 for each $\phi \in \widehat{G},\;$ then there is an equality
of arithmetic classes
\[
\chi  ( \Rr\Gamma  ( \X,\F ) ,f_{\bullet
} ) \cdot \chi  ( {R}\Gamma  ( \X, {\cal H} )
,h_{\bullet } ) =\chi  ( \Rr\Gamma  (\X, \G
 ) ,g_{\bullet } ) . 
\]}
\medskip

\no{\sc Proof.} Let $\mathcal{U}$ denote a $G$-stable affine cover of $\mathcal{X}$.
Then we get the associated exact sequence of Cech complexes 
\[
0\rightarrow {\cal C}^{\bullet }({\cal U},\F)\rightarrow 
{\cal C}^{\bullet }({\cal U},\G)
\rightarrow
{\cal C}^{\bullet }({\cal U},{\cal H})
\rightarrow 0 
\]
For brevity we put ${\cal C}_1^\bullet={\cal C}^{\bullet }({\cal U},\F)$,
${\cal C}_2^\bullet={\cal C}^{\bullet }({\cal U},\G)$, ${\cal C}_3^\bullet=
{\cal C}^{\bullet }({\cal U},{\cal H})$.
 As
mentioned at the start of this section, since the $G$-action is tame, we can
then find perfect $\Z[G]$-complexes with surjective quasi-isomorphisms 
\[
P_{2}^{\bullet \prime }\stackrel{\gamma }{\rightarrow }C_{2}^{\bullet
},\;\;\;P_{3}^{\bullet }\rightarrow C_{3}^{\bullet }\;.\;\;\;\; 
\]
We assert that we can construct a commutative diagram in which the vertical
maps are all surjective quasi-isomorphisms and in which the $P_{i}^{\bullet
} $ are perfect $\Z[G]$-complexes: 
\[
\begin{array}{lllllllll}
0 & \rightarrow & \mathcal{C}^{\bullet }(\mathcal{U},\mathcal{F}) & 
\rightarrow & \mathcal{C}^{\bullet }(\mathcal{U},\mathcal{G}) & \rightarrow
& \mathcal{C}^{\bullet }(\mathcal{U},\mathcal{H}) & \rightarrow & 0 \\ 
&  & \ \uparrow &  & \ \uparrow &  & \ \uparrow &  &  \\ 
0 & \rightarrow & P_{1}^{\bullet } & \rightarrow & P_{2}^{\bullet } & 
\rightarrow & P_{3}^{\bullet } & \rightarrow & 0
\end{array}
\]
The result will then follow on taking bases for the $P_{j}^{i}$ for $j=1,3$
and using these to form bases of the $P_{2}^{i}$ $.$

We now briefly sketch the construction of $P_{1}^{\bullet }$ and$%
\;P_{2}^{\bullet }.$ $\;$By 3.10 we can find a cochain map $\beta $ such
that the following diagram commutes: 
\[
\begin{array}{lll}
C_{2}^{\bullet } & \rightarrow & C_{3}^{\bullet } \\ 
\uparrow \gamma &  & \uparrow \\ 
P_{2}^{\bullet ^{\prime }} & \stackrel{\beta }{\rightarrow } & 
P_{3}^{\bullet }
\end{array}
\]
By adding a free acyclic complex to $P_{2}^{\bullet \prime }$ we may assume
that $\beta $ is surjective; this then implies that $\ker \beta $ is a
perfect complex, and so the restriction of $\gamma $ to $\ker \beta $
provides a quasi-isomorphism to $C_{1}^{\bullet }$. By adding a free acyclic
complex to $\ker \beta $ we obtain a surjective quasi-isomorphism onto $%
C_{1}^{\bullet }$, and the resulting complex is denoted $P_{1}^{\bullet
}.\;\;\;\Box$

\medskip

\no{\bf Proposition 5.6} \ {\sl Suppose that the $G$-$\X$ 
 sheaf $\mathcal{F}$ is fibral, that is to say it is supported over
a finite set of primes $S$ in $\Spec(\Z)$. Then the equivariant determinants of cohomology all
identify with the trivial complex line $\C$, which we endow
with the standard metric $\left| -\right| $ and
\[
\chi  ( \Rr\Gamma  ( \X, \F) ,\left| -\right|
_{\bullet } ) =\nu \circ f_{\ast }^{T} ( \F) 
\]
where $f_{\ast }^{T}$ denotes the composition of 
\[
\oplus _{p\in S}{\rm K}_{0} ( G,\X_p) \stackrel{\oplus
f_{p\ast }}{\rightarrow }\oplus _{p\in S}{\rm K}_{0} ( {\bf F}
_{p}[G] ) \rightarrow {\rm K}_{0}{\rm T}(\Z[G]) . 
\]
 Here the first map is induced by the structure maps $f_{p}:%
\X_p\rightarrow \Spec({\bf F}_p)$ for
$p\in S$ (see Theorem 1.1 in [C]), and the second map is as
described in 4.B.}
\medskip

\no{\sc Proof.} This follows at once from the definition of $\chi  (
\Rr\Gamma  (\X,\F) ,\left| -\right| _{\bullet
} ) $ and from 4.B. (Note that this, in part, justifies the choice of
convention in 4.B.) \ \ $\Box $
\medskip

\no{\bf 5.B \ Rings of integers.}
\medskip

The remainder of this article is devoted to the study of images of
arithmetic classes in various arithmetic situations. In this sub-section we
shall consider the case where $\mathcal{X}$ is the spectrum of a ring of
integers; thus in\ this sub-section we consider the case $\X=\Spec(\Oo_N)$ 
for a ring of integers $\Oo_{N}$ of a number
field $N$ which is at most tamely ramified over a number field $K$, with $%
N/K $ Galois and $G=\Gal(N/K)$.

Our main result here is Theorem 5.9, which is closely related to the work of
Fr\"{o}hlich in Chap. VI of [F2] and to the proof of the Second Fr\"{o}hlich
Conjecture in [CNT].

\medskip

Suppose that $\frak{a}$ is a $G$-stable $\Oo_{N}$-ideal and let $\F=
\widetilde{\frak{a}}$ \ be the associated $G$-$\X$ sheaf viewed as a
complex concentrated in degree zero. As $\mathcal{X}$ is affine 
\[
\Hh^{i}(\X,\F) =\left\{ 
\begin{array}{l}
{\ \frak{a\,\;\;}}\hbox{\rm if\ }i=0 \\ 
\{0\}\hbox{\rm \ \ if}\;i>0.
\end{array}
\right. 
\]
We endow ${\frak{a}}_{\C}=\C\otimes_\Z{\frak{a}}=\C\otimes_{\Q}N$ with the $G$-invariant positive definite \textit{Hecke form}
$h: \C\otimes_{\Q}N\times \C\otimes_{\Q}N   \rightarrow \C$ which is defined by the
rule 
\[
h\left( \lambda \otimes m,\nu \otimes n\right) =\lambda \overline{\upsilon }%
\sum_{\sigma }\sigma \left( m\right) \overline{\sigma \left( n\right) } 
\]
where the sum extends over the embeddings $\sigma :N\rightarrow \C$.
Thus, as in (2.5), $h$ determines metrics on the $\det  (  ( \C\otimes_{\Q}N ) _{\phi } ) $ for $\phi \in \widehat{G};$
we denote this set of metrics by $\det h_{\bullet }.$

\smallskip

\no{\bf Remark:} We refer to the form $h$ as the Hecke form since
this form was introduced by Hecke in his proof of the functional equation
for L-functions; see for instance 9.3 in [He].

\medskip

We write $\mu _{K}$ for the $G$-invariant positive hermitian form on $\C\otimes_\Q K[G]$ given by the rule 
\[
\mu _{K} ( \tsum_{g}x_{g}g,\tsum_{h}y_{h}h ) =\tsum_{\rho
}\tsum_{g,h}\delta _{g,h}\rho \left( x_{g}\right) \overline{\rho \left(
y_{h}\right) } 
\]
where the first right-hand sum extends over all embeddings $\rho $ of $K$
into $\Bbb{C}$. Again as per (2.5) $\mu _{K}$ induces metrics $\det \mu
_{K,\phi }\;$on the $\det  (  (
\C\otimes_\Q K[G]  )
_{\phi } ) $ for each $\phi \in \widehat{G}$; we denote this set of
metrics by $\det \mu _{K\bullet }$, or $\det \mu _{\bullet }$\ when $K$ is
clear from the context.

In the sequel, since $\X=\Spec(\Oo_N)$ is
affine, for brevity we shall write $\chi ( {\frak{a}}, \det h_{\bullet } ) $ in place of $\chi  ( R\Gamma  ( 
\X, {\frak{a}}), \det h_{\bullet
} ) $ etc. The following result is an equivariant version of the usual
discriminant-index theorem:
\medskip

\no{\bf Proposition 5.7} \ {\sl The following equality holds
in $A(\Z[G])$
\[
\chi  ( \Oo_{N} , \det h_{\bullet } )\cdot  \chi  ( {\frak{a}} 
, \det h_{\bullet } ) ^{-1}=\nu  ( \Oo_{N}/{\frak{a}} ) 
\]
where $\nu $ is the map on torsion classes of 4.B.}
\medskip

\no {\sc Proof.} This follows from Propositions 5.5 and 5.6 applied to the exact
sequence 
\[
\ \;0\rightarrow {\frak{a}}\rightarrow \Oo_{N}\rightarrow \Oo_{N}/{\frak{a}}%
\rightarrow 0.\;\;\;\;\Box 
\]

\medskip

\no{\bf Definition 5.8} \ For a given prime ideal $\frak{p}$ of 
$\Oo_{K},$ let $f_{\frak{p}}$ denote the residue class extension degree of $%
\frak{p}$ in $K/\Q$, denote by $I_{\frak{p}}$ the inertia group of a
chosen prime ideal of $\Oo_{N }$ above $\frak{p}$; and let $u_{\frak{p}%
} $ denote the augmentation character of $I_{\frak{p}}$ (that is to say the
regular character minus the trivial character). Define 
\[
\Pf_{p}\left( \Oo_{N}\right) :R_{G}^{s}\rightarrow \left( -p\right) ^{\Bbb{Z}}%
\hbox{\rm \ by the rule \ }\Pf_{p} ( \Oo_{N} ) ( \psi ) =\prod_{%
{\frak{p}}\mid p}( -p)^{\frac{1}{2}f_{\frak{p}} ( \psi ,\text{%
Ind}_{I_{\frak{p}}}^{G}u_{\frak{p}} ) } 
\]
for $\psi \in R_{G}^{s}$, where $ ( , ) $ denotes the standard
inner product on $R_{G}$. Note that $ ( \psi ,{\rm Ind}_{I_{\frak{p}%
}}^{G}u_{\frak{p}} ) = ( \psi \mid _{I_{\frak{p}}},u_{\frak{p}%
} ) $ is an even integer, since $\psi \mid _{I_{\frak{p}}}$ is a
symplectic character of the cyclic group $I_{\frak{p}}\;$and is therefore a
sum of characters of the form $\theta +\overline{\theta }$. We then define $%
\Pf\left( \Oo_{N}\right) $ to be the idele valued function, defined on
symplectic characters, which is $\Pf_{p}\left( \Oo_{N}\right) $ at primes over $%
p$ and which is 1 at the archimedean primes. 
\medskip

Let $\delta _{K}\in \Hom ( R_{G}^{s},\R_{>0} ) $ be the
homomorphism 
\[
\delta _{K}(\psi )= ( \left| G\right| ^{\left[ K:\Bbb{Q}\right] }\left|
d_{K}\right|  ) ^{\psi (1)/2} 
\]
where $d_{K}$ is the discriminant of $K/\Q$.

The main result of this sub-section is the following description of the tame
arithmetic classes $\chi ^{s} ( \Oo_{N},\det h_{\bullet } )$ and 
$\chi ^{s} ( \Oo_{K}G, \det \mu _{K\bullet })$.
\medskip

\no{\bf Theorem 5.9} \ {\sl (a) The class $\chi ^{s} ( \Oo_{N},\det h_{\bullet } ) $
 in $A^s_T(\Z[G])$ is represented by
the homomorphism $\widetilde{\varepsilon }_{\infty }^{s} ( K )
^{-1}Pf ( \Oo_{N}) ^{-1}\times \delta _{K},$ where, for a
symplectic character $\psi ,$ $\widetilde{\varepsilon }_{\infty }^{s} (
K )  ( \psi  ) $  is the archimedean epsilon factor $
\varepsilon _{\infty } ( K,\psi -\psi  ( 1 )\cdot 1 ) $ 
\ (see Section 7);

 (b) the class $\chi ^{s} ( \Oo_{K}G,\det \mu _{K\bullet } ) $ 
in  $A^s_T(\Z[G])$  is represented by
the homomorphism $1\times \delta _K$.}
\medskip

Before proceeding with the proof of the theorem, we first introduce some
notation and establish some preparatory results.

For a prime number $p,$ let $\beta _{p}$ be an $\Oo_{K,p}[G]$-basis of $\Oo_{N,p}$
and let $b$ be a $K[G]$-basis of $N$ (so that $b$ is a so-called normal basis
of $N/K$). Recall (see I.4 of [F2]) that for a character $\psi $ of $G$ the
resolvent $\left( b\mid \psi \right) $ is defined to be the value $\Det(
\sum_{g\in G}g(b)g^{-1} ) \left( \psi \right) $; note that, with the
notation of (3) in 3.B, $\left( b\mid \psi \right) =\Det\left( r(b)\right) (%
\overline{\psi })$ and so by (4) we have proved the following particular
instance of the Galois action formula for resolvents (cf. Theorem 20A in
[F1]) 
\begin{equation}
\left( g\left( b\right) \mid \psi \right) =\left( b\mid \psi \right) \cdot \det(\psi)(g).
\end{equation}
The local resolvents $\left( \alpha _{p}\mid \psi \right) \;$are defined
similarly (see \textit{loc. cit.}).

Set $\Omega _{K}=\Gal(\overline{\Bbb{Q}}/K)$ and recall that we write $%
\Omega $ for $\Omega _{\Bbb{Q}}$. For an $\Omega$-module $A,$ let 
\[
{\cal N}_{K/\Bbb{Q}}: \Hom_{\Omega _{K}} ( R_{G},A )
\rightarrow \Hom_{\Omega } ( R_{G},A ) 
\]
denote the co-restriction map of (3.3) in II.3 of [F1]; we extend the domain
of this map to include resolvents, which are not in general $\Omega_{K}$%
-equivariant, as per (3.1) in III.3 in [F1].

Next we recall the $p$-adic absolute value function from [CEPT2]. Let $L=%
\Bbb{Q}\left( \zeta _{p}\right) .\;$By Lemma 3.1 \textit{loc. cit} we know
that we can find $\lambda \in L_{p}$ such that $\lambda ^{p-1}=-p.$ Once and
for all we fix such a choice of $\lambda $ and define $\left\| -\right\| :%
\overline{\Bbb{Q}}_{p}^{\times }\rightarrow \lambda ^{\Bbb{Q}}$ by
stipulating that $x\cdot\left\| x\right\| $ is to be a $p$-adic unit.

We now recall some related constructions from [CEPT2]; for full details see
(3.1) and (3.2) in \textit{loc. cit. }Let $R_{G}\left( \overline{\Bbb{Q}}%
_{p}\right) $ denote the ring of $\overline{\Bbb{Q}}_{p}$-characters of $G$
and set $\Omega _{p}=\Gal ( \overline{\Q}_{p}/\Q_{p} ) .$
For $g\in \Hom ( R_{G} ( \overline{\Q}_{p} ) ,\overline{\Q}_{p}^{\times } ) $ 
define $\left\| g\right\| ( \varphi
 ) =\left\| g ( \varphi  ) \right\| $. We shall say that $g$ is%
\textit{\ well-defined} if $\left\| g\right\| $ takes values in $\lambda ^{%
\Bbb{Z}}.$

Once and for all we fix a field embedding $h:\overline{\Bbb{Q}}\rightarrow 
\overline{\Bbb{Q}}_{p}$. From II.2.1 in [F1], $h$ induces an isomorphism 
\[
h^{\ast }:{\Hom}_{\Omega_{L}} ( R_{G}, ( \overline{\Bbb{Q}}%
\otimes \Q_{p} ) ^{\times } ) \cong {\Hom}_{\Omega
_{L_{p}}} ( R_{G} ( \overline{\Bbb{Q}}_{p} ) ,\overline{\Bbb{Q}}%
_{p}^{\times } ) . 
\]
For $f\in {\Hom}_{\Omega_{L}} ( R_{G}, ( \overline{\Bbb{Q}}%
\otimes \Q_{p} ) ^{\times } ) $ define $\left\| f\right\|
=h^{\ast -1}\left( \left\| h^{\ast }f\right\| \right) $; we shall say that $%
\left\| f\right\| $ is well-defined when $\left\| h^{\ast }f\right\| $ is
well-defined.

In the sequel we employ a standard abuse of notation and write $\Det (
\Oo_{L_{p}}[G]^{\times }) $ for $h^{\ast } ( \Det (
\Oo_{L_{p}}[G]^{\times }))$.

\medskip

\no{\bf Theorem 5.10} \  {\sl For $\psi \in R_{G}^{s}$, \ 
$
{\rm sign}\left({\cal N}_{K/\Q} ( b\mid \psi  )\right) =\varepsilon
_{\infty } (K, \psi -\psi  ( 1 ) 1) . 
$}
\medskip

\no{\sc Proof.} This is III.4.9 of [F1].\ \ \ \ $\Box $

\medskip

\no{\bf Proposition 5.11} \ {\sl We have\ 
${\cal N}_{K/\Q} ( \beta _{p}\mid - ) ^{s}\cdot \Pf_{p} (
\Oo_{N} ) ^{-1}\in {\Det}^{s} ( \Oo_{T,p}G^{\times })$  
where we recall from 4.B that $T$  denotes the maximal
abelian tame extension of $\Bbb{Q}$ in }$\overline{\Q}$.
\medskip

\no{\sc Proof.} Let $\tau ^{\ast }$ denote the adjusted Galois Gauss sum of (3.9) in
[T1] (or see IV.1.7. in [F1]). From the discussion following Theorem 2 in
[T1] we know that we can find $z_{p}\in \Z[G]^{\times }$ such that 
for all $\phi \in R_{G}$%
\[
{\cal N}_{K/\Q} ( \beta _{p}\mid \phi  ) =\Det (
z_{p} )  ( \phi  ) \tau ^{\ast } ( \phi  ) . 
\]

Recall that we have fixed a choice of field embedding $h:\overline{\Q}%
\rightarrow \overline{\Q}_{p}$. By Theorem 4 in [CEPT2] we know that $%
\left\| \tau _{p}^{\ast }\right\| =\left\| \widetilde{\varepsilon }%
_{0,p}\right\| $ is well-defined. Writing 
\[
\tau ^{\ast }=\tau ^{\prime }\tau _{p}^{\ast }\hbox{\rm \ where \ }\tau
^{\prime }=\tprod_{q\neq p}\tau _{q}^{\ast } 
\]
we get 
\[
{\cal N}_{K/\Q} ( \beta _{p}\mid - )
 =\Det(
z_{p} ) \tau ^{\prime }\tau _{p}^{\ast }={\Det}(z_{p})
\tau ^{\prime }\tau _{p}^{\ast }\left\| \tau _{p}^{\ast }\right\| \left\|
\tau _{p}^{\ast }\right\| ^{-1} 
\]
and by Theorem 4 in \textit{loc. cit} 
\[
\tau _{p}^{\ast }\left\| \tau _{p}^{\ast }\right\| \in {\Det} (
\Oo_{L_{p}}[G]^{\times } ) \;\hbox{\rm \ and \ }\tau ^{\prime }\in \Det
( \Oo_{T,p}[G]^{\times } ) . 
\]

>From Theorem 7.4 in [M] we know that each value of $\tau _{p}^{s}$ is plus
or minus an integral power of $p$. Thus for a symplectic character $\psi $
of $G,$%
\[
\left| \tau _{p}^{\ast }\left( \psi \right) \right| =N{\frak{f}}_{p}\left(
\psi \right) ^{\frac{1}{2}}=\tprod_{{\frak{p}}\mid p}N{\frak{p}}^{\frac{1}{2}%
 ( \psi ,\text{Ind}_{I_{\frak{p}}}^{G}\;u_{\frak{p}} ) }=\pm
\Pf_{p}\left( \Oo_{N}\right)(\psi) . 
\]
where $N\frak{f}_{p} ( \psi ) $ denotes the $p$-part of the
absolute norm of the Artin conductor of $\psi $. As$\;\left\| \tau
_{p}^{\ast s}\right\| $ and $Pf_{p} ( \Oo_{N} ) $ are both integral
powers of $-p$, we deduce that $\left\| \tau _{p}^{\ast s}\right\| ^{-1}$= $\Pf_{p}\left( \Oo_{N}\right) $ as required.\ \ \ \ $\Box $

\medskip

\no{\sc Proof of Theorem 5.9. } \ We begin by proving (a). Let $ \{
x_{i} \} $ denote a $\Z$-basis of $\Oo_{K}$. Then $\left\{
x_{i}b\right\} $ resp. $\left\{ x_{i}\beta _{p}\right\} $ is a $\Q[G]$%
-basis resp. a $\Z_p[G]$-basis for $N$ resp. $\Oo_{N,p}$. With the
previous notation choose $\lambda _{p}\in $ $K_{p}[G]$ such that $b=\lambda
_{p}\beta _{p}$ and write $x_{i}\lambda _{p}=\tsum_{j}\lambda _{p}^{ij}x_{j}$
with $\lambda _{p}^{ij}\in \Q_{p}[G].$ Then 
\[
x_{i}b=x_{i}\lambda _{p}\beta _{p}=\tsum_{j}\lambda _{p}^{ij}x_{j}\beta _{p} 
\]
and so the matrix $ ( \lambda _{p}^{ij} ) _{ij}$ transforms the $%
\Q_{p}[G]$-basis $\left\{ x_{i}\beta _{p}\right\} $ into the basis
$\left\{ x_{i}b\right\} ;$ therefore the finite coordinate of the
representing homomorphism of the arithmetic class $\chi ( \Oo_{N},
\det h_{\bullet } )$ is
\[
\prod_{p}{\Det}(\lambda _{p}^{ij})=\prod_{p}{\cal N}_{K/\Q
}(\lambda _{p})=\prod_{p}{\cal N}_{K/\Q}(b\mid -)\cdot {\cal N}_{K/\Bbb{%
Q}}(\beta _{p}\mid -)^{-1}. 
\]
To obtain the archimedean coordinate for a chosen irreducible character $%
\phi $ we have to extend our notation and choose a positive integer $n_{\phi}$ 
such that $\det(\phi)^{n_{\phi }\phi \left(
1\right) }$ is trivial. We then write $\psi =n_{\phi }\phi (1)\overline{\phi 
}$  and set 
\[
W_{\psi }=W_{\phi }^{n_{\phi }} 
\]
where $W_{\phi }^{n_{\phi }}\;$denotes the direct sum of $n_{\phi }$ copies
of $W_{\phi }.$ We endow $W_{\psi }$ with the hermitian form, $\nu _{\psi }$
say, given by the orthogonal sum of the hermitian forms on the $W_{\phi }$,
and we let $\left\{ w_{\psi ,k}\right\} $ denote the basis of $W_{\psi }\;$%
derived from the bases $\left\{ w_{\phi ,l}\right\} $ of $W_{\phi }.$ $\;$We
must now consider the wedge product 
\begin{equation}
\bigwedge _{i,k}\left( x_{i}\cdot r(b)(1\otimes w_{\psi ,k}\right) )=\bigwedge
_{i}\bigwedge _{k} ( \tsum_{g}x_{i}g\left( b\right) \otimes gw_{\psi
k} ) =\bigwedge _{i,k}y_{i}(1\otimes w_{\psi ,k})
\end{equation}
where $y_{i}=x_{i}(b\mid \psi )$; here we obtain the second equality from
the fact that 
\[
\bigwedge_{k} ( \tsum_{g}g\left( b\right) \otimes gw_{\psi k} ) =%
{\Det} ( r ( b )  )  ( \overline{\psi } )
\tbigwedge_{k}w_{\psi k}=(b\mid \psi )\tbigwedge_{k}w_{\psi k}. 
\]
A priori ($b\mid \psi )\in \C\otimes _{\Q}N$; however, because $%
\det \left( \psi \right) =1$, by the Galois action formula (10), ($b\mid
\psi )\in \C\otimes _{\Bbb{Q}}K$. Therefore 
\[
\bigwedge_{i}x_{i}(b\mid \psi )=N_{K/\Bbb{Q}}(b\mid \psi )\bigwedge _{i}x_{i}. 
\]
To complete the proof of (a), note first that as the $\left\{ x_{i}\right\} $
are fixed by $G$ and as $\left\{ w_{\psi ,k}\right\} \;$is an orthonormal
basis for the form $\nu _{\psi },$%
\begin{equation}
h\otimes \nu _{\psi } ( x_{i}\otimes w_{\psi ,k},x_{j}\otimes w_{\psi
,l} ) =\tsum_{\sigma }\sigma(x_{i})\overline{\sigma \left( x_{j}\right) }%
\delta _{k,l}=\left| G\right| \tsum_{\rho }\rho \left( x_{i}\right) 
\overline{\rho \left( x_{j}\right) }\delta _{k,l}.
\end{equation}
In the sequel we shall write $ ( h\otimes \nu _{\psi } ) ^{G}$ for
the restriction of $h\otimes \nu_{\psi }$ from $(\C\otimes _\Q N)\otimes W$
to $((\C\otimes_\Q N)\otimes W) ^{G}$. Hence the archimedean
coordinate of the representing homomorphism of $\chi  ( \Oo_{N},\det
h_{\bullet } )$ at $\phi $ is the $n_{\phi }\phi ( 1 ) $-st
root of 
\[
\det (h_{\psi } ( \tbigwedge_{i,k}x_{i}r(b)\left( 1\otimes w_{\psi
,k}\right) ))=
\]
\[
=\det \left( \left( h\otimes \nu _{\psi }\right) ^{G}\left( x_{i}r\left(
b\right) \otimes w_{\psi ,k},x_{j}r\left( b\right) \otimes w_{\psi
,l}\right) \right) ^{1/2}
\]
\[
=\left| {\cal N}_{K/\Bbb{Q}}(b\mid \psi )\right| \det \left( \left(
h\otimes \nu _{\psi }\right) ^{G}\left( x_{i}\otimes w_{\psi
,k},x_{j}\otimes w_{\psi ,l}\right) \right) ^{1/2}
\]
\[
=\left( \left| G\right| ^{\left[ K:\Bbb{Q}\right] }\left| d_{K}\right|
\right) ^{\psi \left( 1\right) /2}\left| {\cal N}_{K/\Bbb{Q}}(b\mid \psi
)\right| .
\]
 
Note that the square roots in the above right-hand terms (which are of
course taken to be positive) arise since we are dealing with the metrics
which are, of course, given by the square root of the corresponding positive
definite hermitian forms. This then shows that the class $\chi (\Oo_N,\det h_{\bullet } )$ 
is represented by the homomorphism
which maps an irreducible character $\phi $ to the value 
\[
{\cal N}_{K/\Bbb{Q}}(b\mid \phi )\cdot{\cal N}_{K/\Bbb{Q}}(\beta _{p}\mid
\phi )^{-1}\times \left| {\cal N}_{K/\Bbb{Q}}(b\mid \phi )\right|\cdot \delta
_{K}\left( \phi \right) . 
\]

\smallskip

We now consider $\chi ^{s}( \Oo_{N},\det h_{\bullet } ) .$ Then by Proposition 5.10
this class is represented by the homomorphism which
maps a symplectic character $\psi $ to the value 
\[
{\cal N}_{K/\Bbb{Q}}(b\mid \psi ){\cal N}_{K/\Bbb{Q}}(\beta _{p}\mid
\psi )^{-1}\times \widetilde{\varepsilon }_{\infty }( K,\psi  ) 
{\cal N}_{K/\Bbb{Q}}(b\mid \psi )\delta _{K} ( \psi  ) . 
\]
Since 
\[
 ( \psi \mapsto \widetilde{\varepsilon }_{\infty } ( K,\psi  )\cdot 
{\cal N}_{K/\Bbb{Q}}(b\mid \psi ) ) \ \in\  {\Hom}_{\Omega
}^{+} ( R_{G}^{s},\overline{\Bbb{Q}}^{\times } ) ={\Det}^{s}(%
\Q[G]^{\times }), 
\]
we conclude that the class is also represented by 
\[
\psi \mapsto \widetilde{\varepsilon }_{\infty } ( K,\psi  )
^{-1}\prod_{p}{\cal N}_{K/\Bbb{Q}}(\beta _{p}\mid \psi )^{-1}\times
\delta _{K} ( \psi  ) 
\]
and the result then follows from 5.11.

The proof for (b) is similar, but considerably easier, because we may
replace $b$ and all the $\beta _{p}$ by 1 throughout in the above.\ Indeed,
we see immediately that, with these choices, the finite coordinate is 1.
Since 
\[
 \left( \mu _{K}\otimes \nu _{\phi }\right) ^{G} (
\tsum_{g}x_{i}g\otimes gw_{\phi ,k},\tsum_{h}x_{j}h\otimes hw_{\phi
,l} )
\]
\[
=\tsum_{\rho }\tsum_{g,h}\rho  ( x_{i} ) \overline{\rho \left(
x_{j}\right) }\mu \left( g,h\right) \nu  ( gw_{\phi ,k},hw_{\phi
,l} )
\]
\[
=\delta _{k,l}\left| G\right| \tsum_{\rho }\rho \left( x_{i}\right) 
\overline{\rho \left( x_{j}\right) } 
\]

we have 
\[
\det \mu _{K,\phi } ( \tbigwedge_{i,k}\tsum_{g}x_{i}g\otimes gw_{\phi
,k} ) =\delta _{K} ( \phi  ). \;\;\;\;\Box 
\]
\bigskip

\ \ \ 

\section{Equivariant Quillen Metrics}

\no{\bf 6.A \ Definition of arithmetic classes.}
\medskip

In this section we consider an\ arithmetic variety $\mathcal{X}$ with
fibral dimension $d$. Since $G$ acts tamely on $\mathcal{X}$, $G$ must act
freely on the complexified generic fibre $X:={\X}\times _{\Bbb{Z}}\Bbb{%
C}$; in the sequel we shall abuse terminology and identify $X$ with the complex manifold 
$\X(\C)$ of its complex points. We fix a $G$-invariant
K\"{a}hler metric $h=h^{TX}$ on $X$
which is invariant under complex conjugation.

A hermitian $G$-bundle on $\mathcal{X}$ is a pair $ ( \F
,f ) $, where $\mathcal{F}$ is a locally free $G$-$\X$ sheaf
with the property that the induced  holomorphic vector bundle  
$\F_\C$ over $X$ supports a $G$-invariant
hermitian metric $f$, which is invariant with complex conjugation.

The complex lines $\det ( \Hh^{\bullet } ( \Rr\Gamma  ( 
\X,\F )  )  ) _{\phi },$ for $\phi \in 
\widehat{G},$ carry metrics $f_{L^{2},\phi }$ coming from the $L^{2}$-metric
of Hodge theory for the Dolbeault resolution. As per Section II in [B], the $%
f_{L^{2},\phi }$ can be transformed to equivariant Quillen metrics $%
f_{Q,\phi }$ for $\phi \in \widehat{G}.$ One of the main objectives of this
article is the study of the arithmetic classes 
\[
\chi \left( \Rr\Gamma \left( \X, \F\right) , f_{Q\bullet
}\right) \;\hbox{\rm in\ \ }A(\Z[G]) . 
\]
More generally, we shall also consider a bounded complex $\mathcal{G}%
^{\bullet }$ of hermitian $G$-bundles on $\mathcal{X}$, with $g^{i}$
denoting the hermitian form on $\G^{i}$. Then the $g^{\bullet}$
induce metrics $g_{Q,\phi }^{\bullet }$ on the equivariant determinant of
the hypercohomology of $\G^{\bullet }$, and so the arithmetic class 
\[
\chi  (\Rr\Gamma  ( \X, \G^{\bullet } )
,g_{Q\bullet }^{\bullet }  )\ \hbox{\rm \ in\ \ }A(\Z[G]) 
\]
is defined; explicitly, we may identify the equivariant determinant of
 $\det  ( \Hh^{\bullet }( \Rr\Gamma  ( \X, \G^{\bullet }))) _{\phi }$ with the product 
\[
\det  ( \Hh^{\bullet } ( \Rr\Gamma (\X,\G^{\bullet }))) _{\phi }=\bigotimes _{i}\det ( \Hh^{\bullet } ( 
\Rr\Gamma ( \X, \G^{i}))) _{\phi }^{\left( -1\right) ^{i}} 
\]
and so 
\[
\chi  ( \Rr\Gamma \left( \mathcal{X},\mathcal{G}^{\bullet }\right)
,g_{Q\bullet }^{\bullet } ) =\tprod_{i}\chi  ( \Rr\Gamma \left( 
\X, \G^{i}\right) ,g_{Q\bullet }^{i} ) ^{\left(
-1\right) ^{i}}. 
\]
In the sequel we shall write $g_{Q\bullet }$ for the metrics on the
equivariant determinant of hypercohomology induced by the $ \{
g_{Q\bullet }^{i} \} .$

\medskip

\no{\bf 6.B \ $1$-dimensional subschemes}
\medskip

In this sub-section we place ourselves in the situation described in 5.A. We
now suppose further that $\mathcal{X}$ and $\mathcal{Y}$ are regular
arithmetic varieties, with $\mathcal{Y}$ being connected and having the
property that at, any prime $p$ of non-smooth reduction, $\Y_{p}^{%
\text{red}}$ is a union of smooth irreducible components with strictly
normal crossings and with each such component having multiplicity in $\Y_{p}$ coprime to $p$. As noted previously, since $G$ acts tamely
on $\mathcal{X}$, the branch locus $b$ of the cover $\mathcal{X}/\mathcal{Y}$
is a vertical divisor. By hypothesis $b$ is a divisor with strictly normal
crossings. Here we consider an irreducible regular connected closed
horizontal sub-scheme $\cal Z$ of $\mathcal{X}$ of dimension one; we
may therefore write ${\cal Z}=\Spec(  \Oo_{N} ) $ for some ring
of integers $N$, where $G$ acts tamely on $N$. As previously we put $K=N^{G}$
consider $\F=\Oo_{\mathcal{Z}}$ and endow $\mathcal{F}$
with the Hecke form $h$ of 5.B.

Next we recall the Pfaffian divisor from Section 2 of [CEPT2]: for each
symplectic character $\psi$ of $G$, the Pfaffian divisor $\Pf(\X,
\psi )$ is a divisor on $\mathcal{Y}$ which is supported on the branch
locus $b$. Let $\mathcal{W}$ $=\pi( \mathcal{Z}) $ so that $%
\mathcal{W}$ is a closed sub-scheme of $\mathcal{Y}$. Throughout this
sub-section we shall suppose that $\mathcal{W}$ meets $b$ transversely and
at smooth points of $b.$ As we shall see in the next section, in practice we
can often reduce to this situation by means of a moving lemma - subject to
certain base extensions.

>From Theorem 5.9 we know that $\chi ^{s}( O_{\mathcal{Z}},\det
h_{\bullet }) $ is represented by the homomorphism $\widetilde{\varepsilon }_{\infty }
( K) ^{-1}\Pf( \Oo_{N})
^{-1}\times \delta_{K}$. Let $\{ b_{i}\} $ denote the
irreducible components of $b$; let $\eta _{i}$ denote the generic point of
an irreducible component, $B_{i}$ say, of $\pi ^{-1}(b_{i})$; let $I_{i}$
denote the inertia group of $\eta _{i}$ and recall that $u_{i}$ denotes the
augmentation character of $I_{i}$. From (2.1) in [CEPT2] we know that for 
$\psi \in R_{G}^{s}$
\begin{equation}
\Pf\left( {\X}, \psi \right) =\frac{1}{2}\sum_{i}( \psi
, {\rm Ind}_{I_{i}}^{G}\;u_{i}) b_{i}.
\end{equation}
A closed point $\frak{p}$ of $\mathcal{W}$ (above $p$, say) is ramified in $
\mathcal{Z}$/$\mathcal{W}$ if and only if it is a point of intersection of 
$\mathcal{W}$ and some $b_{i}$. Since we have assumed that $\mathcal{W}$
intersects $b$ transversely at smooth points of $b$,  $I_{\frak{p}}$ is a
conjugate of $I_{i}$ and recall that we denote the residue class degree of
the point $\frak{p}$ by $f_{\frak{p}}$. In the sequel for such a point $%
\frak{p}$ we write $n\left( \frak{p}\right) =i.$ By (5.8) for $\psi \in
R_{G}^{s}$ we have 
\[
\Pf_{p}({\cal Z}, \psi) =\tprod_{\frak{p}}(-p)^{\frac{1}{2}f_{\frak{p}}
( \psi , {\rm Ind}_{I_{n( {\frak{p}})
}}^{G}\;u_{n( {\frak{p}}) }) } 
\]
where the product extends over all points of intersection of $\mathcal{W}$
with the fibre of $b$ above $p$. We therefore denote the right-hand
expression by $\deg \left( {\cal W}\cdot \Pf_{p}\left( \X,
\psi \right) \right)$, and we let $\deg \left( {\cal W}\cdot \Pf\left( 
\X , \psi \right) \right)$ denote the finite idele whose 
$p$-th component is $\deg \left({\cal W}\cdot \Pf_{p}\left( \X , \psi \right) \right)$.
(Note that almost all $p$-components are 1
and that the use of $-p$ in place of\ $p$ means that of course we are using
degree in a non-standard way.) Writing $\widetilde{\varepsilon }_{\infty
}\left( \mathcal{W}\right) $ for $\widetilde{\varepsilon }_{\infty }\left(
K\right)$, we have now shown that the class $\chi ^{s}( \Oo_{\mathcal{Z}
}, \det h_{\bullet }) $ is represented by the homomorphism 
\[
\widetilde{\varepsilon }_{\infty }\left( \mathcal{W}\right) ^{-1}\cdot \deg
\left({\cal W}\cdot \Pf\left( \X\right) \right) ^{-1}\times
\delta _{K}. 
\]
Since the Dolbeault complex of a point is trivial, the equivariant Quillen
metrics associated to the metrics $h_{\bullet }$ are precisely the $\det
\left( h_{\bullet }\right) \;$(cf. Definitions 2.1 and 2.2 in [B]). So
finally we have now \medskip established the main result of this sub-section
\medskip

\no{\bf Theorem 6.1} \ {\sl The symplectic arithmetic class $\chi
^{s}( \Oo_{\mathcal{Z}},\det h_{Q\bullet }) $ is
represented by the homomorphism } 
\[
\widetilde{\varepsilon }_{\infty }( {\cal W})^{-1}\deg ( 
{\cal W}\cdot \Pf(\X)) ^{-1}\times \delta_{K}. 
\]
\medskip

\no{\bf 6.C \ Invariance under passage to degree zero.}
\medskip

In this sub-section we establish a number of results concerning the
independence, with respect to the choice of hermitian metric, of arithmetic
classes after passage to degree zero by the method described in
4.E.  Recall that we denote the complexified generic fibre of 
$\mathcal{X}$ by $X.$
\medskip

\no{\bf Theorem 6.2} \ {\sl Suppose that $\mathcal{F}$ is a
hermitian $G$-bundle on $X$ and let $f,f'$
 be two $G$-invariant hermitian metrics on $\mathcal{F}$. Then
there exists a positive real number $c$ such that
for each $\phi \in \widehat{G}$
\[
f_{Q,\phi }=c^{\phi \left( 1\right) ^{2}}f_{Q,\phi }' 
\]
and so }
\[
\widetilde{\chi }( \Rr\Gamma \F, f_{Q\bullet }) =%
\widetilde{\chi }( \Rr\Gamma \F, f_{Q\bullet }'). 
\]
\medskip

\no{\sc Proof.} For each $\phi \in \widehat{G}$, let $\beta _{\phi }$ be the positive
real number such that $\beta _{\phi }f_{Q,\phi }=f_{Q,\phi }'.$ As
previously we extend $\beta $ to $R_{G}$ by setting $\beta(
\phi) =\beta _{\phi }^{1/\phi \left( 1\right) }$, $\beta \left( \phi
+\psi \right) =\beta \left( \phi \right) \beta \left( \psi \right) $ etc.
In [B], Bismut considers the central function $\sigma $ on $G$
\[
\sigma =\tsum_{\phi \in \widehat{G}}2\log \left( \beta _{\phi }\right) \phi
\left( 1\right) ^{-1}\phi ; 
\]
the Anomaly Formula in Theorem 2.5 of [B] shows that $\sigma(g)$
may be evaluated in terms of integrals over the fixed points of $g$.
However, since $G$ acts freely on $X$, for each $g\in G,\;g\neq 1_{G},$ the
sub-variety of fixed points $X^{g}=\{ x\in X(\C)\mid 
x^{g}=x\}$ is empty. Thus we immediately deduce that $\sigma
\left( g\right) =0$ whenever $g\neq 1_{G}.\;$This then shows that $\sigma $
is a scalar multiple of the regular character and the result follows.\ \ \ $%
\Box $

\medskip

Next we consider the direct image of a hermitian bundle on a closed
sub-scheme of a regular arithmetic variety $\mathcal{X}$. The formation of
standard (i.e. non-hermitian) Euler characteristics respects closed
immersions; however, this need not be the case for arithmetic classes, as
the associated Quillen metrics may change. The precise variation in the
arithmetic classes, that we wish to consider, was determined in Theorem 0.1
in [B].

We begin by considering a $G$-equivariant closed immersion $i:{\cal 
Z}\rightarrow {\X}$ of an arithmetic variety $\mathcal{Z}$ which also supports
a tame action by $G$. Let $\mathcal{F}$ denote a locally free $G$--$\mathcal{Z}$
sheaf. Since $\mathcal{X}$ is regular, we may resolve $i_{*}\mathcal{F}$
by a bounded complex ${\G}^{\bullet }$ of locally free coherent $G$--$\mathcal{X}$ 
modules. We then have natural isomorphisms in the derived
category of $\Z[G]$-modules 
\begin{equation}
\Rr\Gamma ( {\cal Z}, \F) \cong 
\Rr\Gamma ( {\X}, i_{* }\F) \cong 
\Rr\Gamma ( \X , {\G}^{\bullet })
\end{equation}
and hence, for each $\phi \in \widehat{G}$, we obtain isomorphisms 
\begin{equation}
\sigma _{\phi }:\det \Hh^{\bullet }( \Rr\Gamma ( 
\X, \G^{\bullet }))_{\phi }\cong \det 
\Hh^{\bullet }( \Rr\Gamma ( {\cal Z},
\F))_{\phi }.
\end{equation}

In order to describe the relevant metrics that we wish to place on these
determinants of cohomology, we need some further notation. 
Let $Z={\cal Z}_\C$ and let $TZ$ denote the tangent bundle of $Z.$ We let $h^{TZ}$ denote the restriction of $h$ to $TZ$. Let $N_{Z\mid X}$ denote
the normal bundle to $Z$ in $X$ and let $h^{N_{Z\mid X}}$ be the metric
on $N_{Z\mid X}$ induced by $h$. Let $f$ denote a given $G$-invariant metric
on $\mathcal{F}$; we then endow each term ${\cal G}^{i}$ of ${\G}^{\bullet }$ with a $G$-invariant hermitian metric $g^{i}$ in such a way
that the metrics $\{ g^{i}\}$ satisfy Bismut's Condition A with
respect to $h^{N_{Z\mid X}}$ and $f$.

We now wish to compare the arithmetic classes $\chi ( \Rr\Gamma( {\cal Z}, F) ,f_{Q\bullet })$ and $%
\chi (\Rr\Gamma ( \X, \G^{\bullet
}) ,g_{Q\bullet }).$

Let $\alpha _{\phi }$ be the unique positive real number such that under the
isomorphism $\sigma _{\phi }$ of (15) 
\[
\sigma _{\phi }^{* }\left( f_{Q,\phi }\right) =\alpha _{\phi }g_{Q,\phi
}. 
\]
Then by Proposition 5.4 we see that the arithmetic class 
\[
\chi (\Rr\Gamma ( \X, \G^{\bullet
}) , g_{Q\bullet })\cdot \chi ( \Rr\Gamma ( {\cal Z}, \F) , f_{Q\bullet })^{-1} 
\]
is represented by the homomorphism $1\times \alpha ^{-1}\in \Hom_{\Omega _{%
\Bbb{Q}}}( R_{G}^{s}, J_{f}) \times \Hom\left( R_{G}^{s},\Bbb{R}%
_{>0}\right) $ which maps the character $\phi $ to $1\times \alpha _{\phi
}^{-1/\phi \left( 1\right) }$ (so that of course $\alpha \left( \phi \right)
=\alpha _{\phi }^{1/\phi \left( 1\right) }$).
\medskip

\no{\bf Theorem 6.3} \ {\sl With the above notation and hypotheses there
is a positive real number $b$ such that for each $\phi \in 
\widehat{G}$, $\alpha _{\phi }=b^{\phi \left( 1\right) ^{2}}$
and so 
\[
\widetilde{\chi }(\Rr\Gamma ( \X, \G^{\bullet }) , g_{Q\bullet }) =\widetilde{\chi }( \Rr
\Gamma ( {\cal Z}, \F) ,f_{Q\bullet }) . 
\]}

\no{\sc Proof.} In Theorem 0.1 in [B] Bismut considers the central function $\tau$
\[
\tau =\tsum_{\phi \in \widehat{G}}2\log \left( \alpha _{\phi }\right) \phi
\left( 1\right) ^{-1}\phi . 
\]
and shows that $\tau \left( g\right) $ may be evaluated in terms of
integrals over the fixed points of $g.$ As in the proof of 6.2 we deduce
that $\tau \left( g\right) =0$ whenever $g\neq 1_{G}.\;$This then shows that 
$\tau $ is again a scalar multiple of the regular character and the result
follows.\ \ \ $\Box \smallskip $

We now interpret the above results in terms of arithmetic classes.
\medskip

\no{\bf Proposition 6.4} \ {\sl Let $\left( \F_{j}, f_{j}\right) 
$ for $j=1,\ldots, n$ and $\left( \G_{k}, g_{k}\right) $ for $k=1,\ldots, m$ be hermitian bundles
on closed $G$-subschemes $i_{j} : {\cal Z}_{j}\rightarrow 
\X$, $i_{k}:{\cal W}_{k}\rightarrow \X$
such that  
\[
\tsum_{j}\left[ i_{j*}\F_{j}\right] =\tsum_{k}\left[ i_{k*}
\G_{k}\right] \ \ \hbox{\rm in \ K}_{0}( G,\X) . 
\]
Then there is an equality of classes in $A(\Z[G])$
\[
\tprod_{j}\widetilde{\chi }( \Rr\Gamma ({\cal Z}_{j}, \F_{j}) , f_{j, Q\bullet }) =\tprod_{k}\widetilde{%
\chi }( \Rr\Gamma ( {\cal W}_{k}, \G_{k}) ,g_{k,Q\bullet })\ . 
\]}

\no{\sc Proof.} We first choose resolutions by locally free $G$--$\X$ sheaves 
\[
{\cal A}_{j}^{\bullet }\rightarrow i_{j\ast }\F_{j}\; ,\;\;\;
{\cal B}_{k}^{\bullet }\rightarrow i_{k\ast }\G_{k}. 
\]
>From the definition of ${\rm K}_{0}(G, \X)$, we can find locally free $G$--$\X$ sheaves 
$D_{a,b},\;E_{c,d}$ and an isomorphism, which we
henceforth treat as an equality, 
\[
\oplus _{b}D_{2,b}\oplus _{d}E_{1,d}\oplus _{i}E_{3,d}\oplus _{j,a\;{\rm
even}}{\cal A}_{j}^{a}\oplus _{k,b\;{\rm odd}}{\cal B}_{k}^{b} 
\]
\begin{equation}
=\oplus_{d}E_{2,d}\oplus_{b}D_{1,b}\oplus_{j}D_{3,b}\oplus
_{j,a\;{\rm odd}}{\cal A}_{j}^{a}\oplus _{k,b\;{\rm even}}{\cal B}_{k}^{b}
\end{equation}
where the $G$--$\X$ sheaves $D_{a,b}$, $E_{c,d}$ fit into exact
sequences 
\[
0\rightarrow E_{1,d}\rightarrow E_{2,d}\rightarrow E_{3,d}\rightarrow 0 
\]
\[
0\rightarrow D_{1,b}\rightarrow D_{2,b}\rightarrow D_{3,b}\rightarrow 0. 
\]
We then endow the sheaves $E_{3,d}$ and $D_{3,b}$ with arbitrary 
$G$-invariant metrics $\xi_{3,d}$ and $\eta_{3,b}$; we then choose 
$G$-invariant metrics $\xi _{1,d}$, $\xi _{2,d}$, $\eta _{1,b}$, $\eta _{2,b}$ on
$E_{1,d}$, $E_{2,d}$, $D_{1,b}$, $D_{2,b}$ satisfying Condition A as above, so that
by Theorem 6.3: 
\[
\widetilde{\chi }( \Rr\Gamma D_{1,b}, \xi_{1,d,Q})\cdot \widetilde{\chi }%
( \Rr\Gamma D_{3,b}, \xi_{3,d,Q}) =\widetilde{\chi }( \Rr\Gamma
D_{2,b}, \xi_{2,d,Q}) \;\;{\rm etc.}
\]
We then endow the sheaves ${\cal A}_{j}^{a}, {\cal B}_{k}^{b}$ with 
$G$-invariant metrics $\alpha _{j}^{a}$, $\beta _{k}^{b}$ satisfying Condition A,
so that by Theorem 6.3 
\[
\widetilde{\chi }( \Rr\Gamma {\cal A}_{j}^{\bullet }, 
\alpha_{j,Q\bullet }) =\widetilde{\chi }( \Rr\Gamma ( i_{j\ast }
\F_{j}), f_{j,Q\bullet }) \;\;{\rm etc.} 
\]
The desired equality then follows from (16) and Theorem 6.2.\ \ \ $\Box$

\bigskip

\section{Logarithmic Differentials}

In this section we consider the Arakelov-Euler characteristic associated to
the logarithmic de Rham complex of an arithmetic variety $\mathcal{X}$ with
fibral dimension $d$. We begin by relating this class to an arithmetic class
associated to the top Chern class of the logarithmic differentials of 
$\mathcal{X}$. After allowing for various innocuous base field extensions,
we shall use the moving techniques of [CPT1] to express this top Chern class
as a difference of two horizontal 1-cycles together with a relatively
innocuous fibral term. We shall then be able to use the results of 5.B to
show that the arithmetic class associated to the logarithmic de Rham
complex of $\mathcal{X}$ has the remarkable property of characterising
symplectic $\varepsilon_{0}$-constants of $\mathcal{X}$. 

Recall that we have fixed a K\"{a}hler metric $h$\ on the tangent bundle of 
$X=\X(\C)$ and we let $h^{D}$ denote the dual
metric induced on the cotangent bundle of $X.$

In this section we again suppose, that $\mathcal{X}$ and $\mathcal{Y}$ are
as described in 6.B. Let $S$ denote a finite set of prime numbers which
contains all the primes which support the branch locus, together with all
primes $p$ where the fibre $\Y_{p}$ fails to be smooth. We put 
$S'=S\cup \left\{ \infty \right\} .$

Let $\chi (\Y_{\Q}) =\chi ( \Y(\C)))$ denote the Euler characteristic of the generic
fibre of $\Y$. Note that in all cases $d\cdot\chi ( \Y_{\Q})$ is an even integer, 
so that we may define $\xi_{S}: R_{G}\rightarrow \Q^{\times }$ by the rule 
\[
\xi _{S}( \theta) =\tprod_{p\in S}p^{-\theta (1)
\cdot d\cdot\chi ( \Y_{\Q}) /2}\ . 
\]

Let $\Omega_{\Y/\Z}^{1}( \log \Y_{S}^{\rm red}/\log S) $ denote the sheaf of degree one
relative logarithmic differentials with respect to the morphism $( 
\Y, \Y_{S}^{\rm red}) \rightarrow ( 
\Spec(\Z), S)$ of schemes with
log-structures (see [K]). Under our hypotheses $\Omega_{\Y/\Z}^{1}( \log \Y_{S}^{\rm red}/\log S)$
is a
locally free $\Y$-sheaf of rank $d$, and furthermore the
cover $\X$/$\Y$ is log-\'{e}tale, so that 
\begin{equation}
\Omega_{\X/\Z}^{1}( \log \X_{S}^{\rm red}/\log S)=\pi^*
\Omega_{\Y/\Z}^{1}( \log \Y_{S}^{\rm red}/\log S)\ .
\end{equation}

The main goal of this section is the study of the arithmetic class (see 6.A) 
\[
{\frak{c}}=\chi ( \Rr\Gamma ( \wedge^{\bullet }\omlogX  , \wedge ^{\bullet }h_{Q}^{D}) 
\]
\[
=\tprod_{i=0}^{d}\chi ( \Rr\Gamma ( \wedge^i\omlogX ,\wedge ^{i}h_{Q}^{D}) ^{\left( -1\right) ^{i}}. 
\]

To explain our main result we need to introduce some notation on $%
\varepsilon_{0}$-constants. For a more detailed account see Sect. 4 in
[CEPT2] and Sect. 2 and Sect. 5 in [CEPT1]. For a given prime number $p$,
we choose a prime number $l=l_{p}$ which is different from $p$ and we fix
a field embedding $\Q_{l}\rightarrow \C$; then, following the
procedure of Sect. 8 in [D], each of the \'{e}tale cohomology groups 
$\Hh_{\acute{e}t}^{i}(\X\times\bar\Q_p, \Q_{l})$ for 0$\leq i\leq 2d$, affords a continuous complex
representation of the local Weil-Deligne group. Thus, after choosing both an
additive character $\psi_{p}\;$of $\Q_p$ and a Haar measure $%
dx_{p}$ of $\Q_{p}$, for each complex character $\theta$ of $G$ the
complex number $\varepsilon_{0,p}( \Y, \theta , \psi_{p}, dx_{p}, l_{p})$
is defined. (For a representation $\ V$ of $G$
with character $\theta$ this term was denoted $\varepsilon_{p,0}(X\otimes_{G}V, \psi_{p}, dx_{p}, l)$
in 2.4 of [CEPT1].) Setting
$
\widetilde{\varepsilon }_{0,p}( \Y, \theta , \psi
_{p}, dx_{p}, l_{p}) =\varepsilon_{0,p}( \Y, 
\theta -\theta (1)\cdot 1, \psi _{p}, dx_{p}, l_{p})$, 
by Corollary 1 to
Theorem 1 in [CEPT2] we know that when $\theta$ is symplectic, $%
\widetilde{\varepsilon}_{0,p}( \Y, \theta , \psi
_{p}, dx_{p}, l_{p})$ is a non-zero rational number, which is
independent of choices, and 
$\theta \mapsto \widetilde{\varepsilon }_{0,p}( \Y, \theta)$ defines an element 
\[
\widetilde{\varepsilon}_{0,p}^{s}( \Y) \in \Hom_{\Omega }( R_{G}^{s},\Q^{\times }) .
\]
In the case where $\mathcal{X}$ is the spectrum of a ring of integers $\Oo_N$
of a number field $N$ and $K=N^{G}$, we shall write $\varepsilon_{0,p}(K) $ for 
$\varepsilon_{0,p}(\Y)$.

Analogously, for the Archimedean prime $\infty $ of $\Bbb{Q}$, Deligne
provides a definition for $\varepsilon_{\infty }( \Y) $
and from 5.5.2 and 5.4.1 in [CEPT1] we recall that 
\[
\widetilde{\varepsilon }_{\infty }^{s}\left( \mathcal{Y}\right) \in \Hom_{\Omega }(R_{G}^{s}, {\pm }1) . 
\]

For $\phi \in R_{G}^{s}$ almost all $\widetilde{\varepsilon }%
_{0,v}^{s}\left( \mathcal{Y},\phi \right) $ are equal to 1; the global $%
\widetilde{\varepsilon }_{0}$-constant of $\phi $ is 
\[
\widetilde{\varepsilon }_{0}^{s}\left( \mathcal{Y},\phi \right) =\tprod_{v}%
\widetilde{\varepsilon }_{0,v}^{s}\left( \mathcal{Y},\phi \right) 
\]
and we define 
\[
\varepsilon _{0,S}^{s}\left( \mathcal{Y},\phi \right) =\widetilde{%
\varepsilon }_{0}^{s}\left( \mathcal{Y},\phi \right) \tprod_{v\in S^{\prime
}}\varepsilon _{0,v}\left( \Y,\phi \left( 1\right) \right) . 
\]
The main result of this section is

\medskip

\no{\bf Theorem 7.1} \ {\sl The arithmetic class ${\frak{c}}^{s}$
lies in the group of rational classes $R_{T}^{s}(\Z[G])$
and } 
\[
\theta ( {\frak{c}}^{s}) =\xi_{S}^{s}\cdot \varepsilon_{0,S}^{s}( 
\Y)^{-1}. 
\]

By way of preparation for the proof of Theorem 7.1, we shall initially work
with an arbitrary locally free $\mathcal{Y}$--sheaf $\mathcal{E}$; only
towards the end of the section shall we need to specialise to the case where 
${\cal E}=\Omega_{\Y/\Z}^{1}( \log \Y_{S}^{\rm red}/\log S)$. Throughout this section we adopt the
notation and hypotheses of [CPT1]. For $i\geq 0$ let $c^{i}\left( {\cal E}
\right) =\gamma ^{i}( {\cal E}-{\rm rk}( {\cal E})) $
which lies in $F^{i}K_{0}( \Y)$, the $i$-th component of
the Grothendieck $\gamma$-filtration. We define $\underline{c}^{i}({\cal E}
)$ to be the class 
\[
\underline{c}^{i}({\cal E})\equiv c^{i}({\cal E})\ {\rm mod}\
F^{i+1}K_{0}(\Y) . 
\]
\smallskip

\no{\bf Lemma 7.2} \ {\sl Let $\mathcal{E}$ be as
above, let $\mathcal{L}$ denote an arbitrary line bundle on $
\mathcal{Y}$ and suppose that $n_{0}$ is a given negative
integer. Then there exists an integer $n_{1}\leq n_{0}$ and
integers $l_{n}$ for $n_{1}\leq n\leq n_{0}$, which
depend only on ${\rm rk}({\cal E})$, such that for all $i\geq 0$ }
\begin{equation}
\underline{c}^{i}({\cal E})\equiv \sum_{n=n_{1}}^{n_{0}}l_{n}\underline{c}
^{i}({\cal E}\otimes {\cal L}^{n})\;{\rm mod}\;F^{i+1}K_{0}( 
\Y) .
\end{equation}
\medskip
\no{\sc Proof.} This is Lemma 5.3 in [CPT1] with $\mathcal{L}$ replaced by ${\cal L
}^{-1}$.\ \ $\Box$
\medskip

\no{\bf Definition 7.3} \ Let $m$ be a given positive integer. We shall call a
finite Galois extension $M$ of $\Bbb{Q}$ {\sl harmless for }$m$, if $M/%
\Bbb{\,Q\;}$is non-ramified at $S$ and if the extension degree $[M:
\Q]$ is congruent to $1\;{\rm mod}\;m$.
\medskip

\no{\bf Remark} \ From 9.1.2 in [CEPT1] we know that we can construct harmless
for $m$ extensions whose residue class fields over $S$ are arbitrarily
large.
\medskip

If $m$ is a positive integer and if $M$ is harmless for $m,$ let 
$e: \Spec(\Oo_M)\rightarrow \Spec(\Z)$ be the structure morphism,
write $\Y'$ for the base extension $\Y\times_\Z\Oo_N$, and
$\mathcal{E}^{\prime }$ for the pullback of $\mathcal{E}$ to $\mathcal{Y}%
^{\prime }$.\smallskip

Suppose now that an integer $m$ is given and that $\mathcal{E}$ has rank $d$;
for an integer $n$ we put ${\cal E}(n)={\cal E}\otimes O_{\Y
}(n)$. From 5.1 in [CPT1] we know that we can find a negative integer $%
n_{0} $ with the following property: let $n_{1}$ be an integer chosen as in
Lemma 7.2 with respect to the negative integer $n_{0};$ then for each $%
n,\;n_{1}\leq n\leq n_{0},$ there is an open subset $U_{n}$ of Spec$\left( 
\Bbb{Z}\right) $, which contains $S$, an extension$\;M$ which is harmless
for $m,$\ and a (possibly non-effective) 1-cycle $D_{n}^{\prime }$ on $%
\mathcal{Y}^{\prime }$ whose irreducible components are horizontal and meet $%
b$ transversely and at points which are smooth points of both $D_{n}^{\prime
}$ and $b,$ such that 
\[
c^{d}({\cal E}^{\prime }(n))|_{U_{n}}=[ O_{D_{n}^{\prime }\times
U_{n}}]+T_{n}\ \ \hbox{\rm in}\ \ {\rm K}_{0}(\Y^{\prime
}\times U_{n}) 
\]
where $T_{n}$ is supported on closed points.\ With the notation of Lemma
7.2, we set $U=\tbigcap_{n=n_{1}}^{n_{0}}U_{n}$, so that for all $%
n,\;n_{1}\leq n\leq n_{0}$ we have 
\[
c^{d}({\cal E}^{\prime }(n))|_{U}=[ O_{D_{n}^{\prime }\times U}%
] +T_{n}|_{U}\ \ \hbox{\rm in}\ \ {\rm K}_{0}(\Y^{\prime }\times
U). 
\]
Pushing forward by $e$ and using the fact that $\Oo_{M}$ is free over $\Z$, we get 
\begin{equation}
[ M:\Q]\cdot c^{d}({\cal E}(n))| _{U}=[ \Oo_{e_{\ast
}D_{n}^{\prime }\times U}]+e_{\ast }( T_{n}|
_{U})\ \ \hbox{\rm in}\ \ {\rm K}_{0}(\Y\times U).
\end{equation}
In the sequel we work with a chosen such extension $M$. We write $\widehat{
\pi ^{* }D_{n}'}$ for the normalisation of $\pi ^{\ast
}D_{n}^{\prime }\;$and we endow both $e_{\ast }\pi ^{\ast }D_{n}^{\prime }$
and $e_{\ast }\widehat{\pi ^{\ast }D_{n}^{\prime }}$ with the Hecke form $%
h_{n}$ of 6.B; we denote their arithmetic classes by $\chi \left( e_{\ast
}\pi ^{\ast }D_{n}^{\prime },\det h_{n\bullet }\right) \ $\ and $\widetilde{%
\chi }\left( e_{\ast }\widehat{\pi ^{\ast }D_{n}^{\prime }},\det h_{n\bullet
}\right) $. Proposition 5.7, together with Lemma 7.4 below, shows that the
two resulting classes coincide after passage to degree zero 
\begin{equation}
\widetilde{\chi }\left( e_{\ast }\pi ^{\ast }D_{n}^{\prime },\det
h_{n\bullet }\right) \ =\widetilde{\chi }\left( e_{\ast }\widehat{\pi ^{\ast
}D_{n}^{\prime }},\det h_{n\bullet }\right) .
\end{equation}

\medskip

\no{\bf Lemma 7.4} \ {\sl Suppose that $\F$ is a
coherent $\Y$-sheaf  which is supported on a single prime 
$p$ and suppose: either that $p\notin S$; or that, if 
$p\in S,$ then $\mathcal{F}$ is supported over a finite
number of points of $\mathcal{Y}$. Then $f_{p\ast }( \pi
^{\ast }\mathcal{F}) $ is a free class in K$_{0}( {\bf F}_{p}[G])$.}
\medskip

\no{\sc Proof.} Let $h:\Y\rightarrow\Spec(\Z)$ denote
the structure morphism of $\Y$ and suppose first that $\mathcal{F}$
is the coherent $\mathcal{Y}$-sheaf given by the structure sheaf of a closed
point of $\mathcal{Y}$. As $f_{\ast }=h_{\ast }\pi_{\ast }$ and $\pi_{\ast
}\pi ^{\ast }\F=\F\otimes_{\Oo_{\mathcal{Y}}}\pi _{\ast }\Oo_{\mathcal{X%
}}$ in G$_{0}(G,\Y) $, the result follows readily from
the normal basis theorem.

Suppose now that $p\notin S.$ For a $p$-regular element $g\in G,g\neq 1,%
{\X}_{p}^{g}=\emptyset$, since $G$ acts freely away from $S$.
Thus by the Lefschetz-Riemann-Roch theorem, we know that the Brauer trace of 
$g$ on $f_{p\ast }(\F) $ is zero; hence we may conclude
that $f_{p\ast }(\F)$ is a free class. \ \ \ $\Box
\bigskip $

Recall that $\widetilde{\frak{c}}$ denotes the arithmetic class obtained
from ${\frak{c}}$ by passage to degree zero, as per 4.E. As an intermediate
step towards proving Theorem 7.1, we first show that the result holds in
degree zero:
\medskip

\no{\bf Theorem 7.5} \ {\sl The arithmetic class $\widetilde{\frak{c}}^{s}$
lies in the group of rational classes $R_{T}^{s}(\Z[G])$ and 
\[
\theta ( \widetilde{\frak{c}}^{s}) =\widetilde{\varepsilon_{0}^{s}}( \Y) ^{-1}. 
\]}

\no{\sc Proof.} We apply the above work where we now take $\E=\omlogY$
 and where we take $n_{0}$ sufficiently small and
negative to guarantee that $\E^{D}\left( -n\right) $ has a regular
section for all $n\leq n_{0}$. Recall that $\pi ^{\ast }\E$ is
endowed with the metric $h^{D}$, the dual of the K\"{a}hler metric; we endow 
$\pi ^{\ast }\Oo_{\mathcal{Y}}( n)$ with a chosen $G$-invariant
metric $\nu _{n}$.

By 7.2 together with 7.4 and Proposition 6.4, we know that 
\begin{equation}
\widetilde{\chi }\left( \Rr\Gamma ( \wedge ^{\bullet }\pi ^{\ast }%
\E) ,( \wedge ^{\bullet }h^{D}) _{Q}\right)
=\tprod_{n=n_{0}}^{n_{1}}\widetilde{\chi }\left( \Rr\Gamma ( \wedge
^{\bullet }\pi ^{\ast }{\E}( n) ) , ( \wedge
^{\bullet }h^{D}\otimes \nu _{n}) _{Q}\right)^{l_{n}}.
\end{equation}
Let  ${\cal W}_{n}$ denote the closed one dimensional sub-scheme of $%
\mathcal{Y}$ cut out by the regular section of $\E^{D}(-n) $ and put ${\cal Z}_{n}=
\pi^*{\cal W}_{n};$  so that
we have the Koszul quasi-isomorphism 
\[
\wedge^{\bullet }\E( n)\rightarrow \Oo_{{\cal W}_{n}}.
\]
By 6.4 we know that 
\begin{equation}
\widetilde{\chi }\left( \Rr\Gamma ( \wedge ^{\bullet }\pi ^{\ast }%
\E( n)) , ( \wedge ^{\bullet }h^{D}\otimes
\nu _{n}) _{Q}\right) =\widetilde{\chi }\left( \Rr\Gamma \Oo_{{\cal Z}_{n}},j_{n}\right) 
\end{equation}
where $j_{n}$ denotes the Hecke metric on $O_{{\cal Z}_{n}}.$ Next
observe that $\left[ \wedge^{\bullet }\pi^{\ast }\E(
n) \right] =\left( -1\right)^{d}c^{d}\left( \pi ^{\ast }\E
( n) \right) ,$ and also by 7.4 we know that $\widetilde{\chi }%
\left( e_{\ast }\pi ^{\ast }T_{n},\pi ^{\ast }\left|-\right| \right)=0$;
hence by Proposition 6.4, together with (19), (20), (21) and (22), we may
conclude that 
\[
\tprod_{n=n_{1}}^{n_{0}}\widetilde{\chi }( \Rr\Gamma O_{{\cal Z}_{n}},j_{n}) ^{l_{n}[M:\Q]}=\tprod_{n=n_{1}}^{n_{0}}\widetilde{%
\chi }( e_{\ast }\pi ^{\ast }D_{n}^{\prime },\det h_{n\bullet })
^{\left( -1\right) ^{d}l_{n}}
\]
\begin{equation}
=\tprod_{n=n_{1}}^{n_{0}}\widetilde{\chi }( e_{\ast }\widehat{\pi ^{\ast }
D_{n}^{\prime }},\det h_{n\bullet }) ^{\left( -1\right)
^{d}l_{n}}.
\end{equation}
Let ${\cal C}=\left( -1\right) ^{d}c^{d}\left( \E\right) $ and
consider the restriction of $\mathcal{C}$ to an irreducible component 
$b_{i}$ of $b$ over $p$; in this way we obtain a punctual virtual sheaf whose
length we denote by $n_{i}$. If, for $\psi \in R_{G}^{s},$ we have 
$\Pf_{p}\left( \mathcal{X},\psi \right) =\tsum_{i}q_{i}b_{i}$ (see (13)), then
we may define $\deg \left( {\cal C}\cdot \Pf\left( \X\right)
\right) \left( \psi \right) \in J_{f}$ to be the idele whose component at
primes over $p$ is $\left( -p\right) ^{\Sigma _{i}q_{i}n_{i}}.$ We then use
this construction to define the symplectic arithmetic class ${\frak{h}}\in
A_{T}^{s}(\Z[G])$ to be that class which is represented by
the homomorphism 
\[
\left( \widetilde{\varepsilon }_{\infty }^{s}(\Y)\cdot\deg \left( 
{\cal C}\cdot \Pf\left( \X\right) \right) \right) \times 1.
\]
In fact, from Theorem 1 in [CEPT2], we know that $\frak{h}$ is a rational
class and that moreover 
\begin{equation}
\theta \left( {\frak{h}}\right) =\widetilde{\varepsilon }_{\infty }^{s}\left( 
\Y\right) \tprod_{p<\infty }\widetilde{\varepsilon _{0,p}^{s}}%
\left( \Y\right) =\widetilde{\varepsilon_{0}^{s}}\left( \mathcal{Y}%
\right) .
\end{equation}
By Theorems 5.9 and 6.1, the left-hand arithmetic class in (23) above is
represented by the character function given on characters of degree zero by 
\[
\tprod_{n=n_{1}}^{n_{0}}\left( \widetilde{\varepsilon }_{\infty }^{s}(\widehat{\pi
^{\ast }D_{n}^{\prime }})\deg ( \widehat{\pi ^{\ast }
D_{n}^{\prime }}\cdot \Pf( \X^{\prime }) )
\right) ^{-\left( -1\right) ^{d}l_{n}}\times 1
\]
\begin{equation}
=\tprod_{n=n_{1}}^{n_{0}}\left( \widetilde{\varepsilon }_{\infty
}^{s}(e_{\ast }\pi ^{\ast }D_{n}^{\prime })\deg \left( e_{\ast }\pi ^{\ast
}D_{n}^{\prime }\cdot\Pf\left( \X\right) \right) \right)
^{-\left( -1\right) ^{d}l_{n}}\times 1.
\end{equation}
By (18) and (19) together with Theorem 5.5.2 in [CEPT1] we know that 
\[
\tprod_{n=n_{1}}^{n_{0}}\widetilde{\varepsilon }_{\infty }^{s}(e_{\ast }\pi
^{\ast }D_{n}^{\prime })^{l_{n}}=\widetilde{\varepsilon }_{\infty }^{s}(%
\Y)^{[M:\Bbb{Q}]}.
\]
Again by (18) and (19) 
\begin{equation}
\tprod_{n=n_{1}}^{n_{0}}\left( \deg \left( e_{\ast }\pi ^{\ast
}D_{n}^{\prime }\cdot \Pf\left( \X\right) \right) \right)
^{l_{n}}=\deg ( ( -1) ^{d}{\cal C}\cdot\Pf( 
\X))^{\left[ M:\Q\right] }.
\end{equation}
By (21), (22) and (23) we know that 
\[
\widetilde{\chi }\left( \Rr\Gamma ( \wedge ^{\bullet }\pi ^{\ast }%
\E) ,( \wedge ^{\bullet }h^{D}) _{Q}\right) ^{[M:\Q]}
=\tprod_{n=n_{1}}^{n_{0}}\widetilde{\chi }\left( e_{\ast }\pi
^{\ast }\widehat{D_{n}^{\prime }},\det h_{n\bullet }\right) ^{\left(
-1\right) ^{d}l_{n}}
\]
and by the above work the right-hand class is represented by the same
homomorphism as ${\frak{h}}^{-[M:\Q]}.$ Thus, by varying $M$, we see that 
\[
\widetilde{\chi }( \Rr\Gamma ( \wedge ^{\bullet }\pi ^{\ast }%
\E) , ( \wedge ^{\bullet }h^{D}) _{Q}) ={\frak{h}}^{-1}
\]
and so by (24) 
\[
\theta \left( \widetilde{\chi }( \Rr\Gamma ( \wedge ^{\bullet }\pi
^{\ast }\E) , ( \wedge^{\bullet }h^{D})
_{Q}) \right) =\widetilde{\varepsilon _{0}^{s}}( \Y) ^{-1}.\;\;\;\;\Box 
\]
\medskip 

$\;\;$Before embarking on the proof of Theorem 7.1, we first need a number
of preliminary results.\smallskip

\medskip

\no{\bf Lemma 7.6} \ {\sl (a)\ For a coherent $G$--$\X$ sheaf $\F$
there is a quasi-isomorphism of $\Z[G]$-complexes 
\[
\left( \Rr\Gamma \mathcal{F}\right) ^{G}\cong R\Gamma (\F^G) . 
\]

(b) If $\left(\F ,f\right)$  is a hermitian $G$-bundle 
on $\X$, then there is an equality in $A(\Z[G])$
\[
\chi \left( ( \Rr\Gamma \F) ^{G}, f_{Q,1}\right) =\chi
\left( R\Gamma ( \F^{G}) , ( f^{G})
_{Q}\right) . 
\]}
\smallskip

\no{\sc Proof.} Part (a) follows at once on expressing $\Rr\Gamma \mathcal{F}$
and $\Rr\Gamma ( \F^{G}) $ in terms of Cech complexes for
a given affine cover of $\mathcal{Y}$  (which pulls back to an affine
cover of $\mathcal{X}$, since $\mathcal{X}$/$\mathcal{Y}$ is finite) and
then taking $G$-invariants of the first complex. Part (b) is then immediate since 
$f_{Q,1}$ (the Quillen metric for the trivial character) is constructed by
forming the Quillen metric associated to the restriction of $f$ to the
trivial isotypical component of $\F_{\C},$ namely $(\F_\C)^{G}$.
See II.a in [B] for further details.\ \ \ \ $\Box $

\medskip

Next we note the following elementary result from 3.2:
\medskip

\no{\bf Lemma 7.7} \ {\sl When $G$ is the trivial group, then
there is an isomorphism $\gamma :A\left(\Z\right) \rightarrow \R_{>0},$ 
(which coincides with the degree map on p. 162 of [S]). Furthermore, if a class ${\frak{e}}\in A\left( \Z
\right)$ has $\gamma \left( \mathit{\frak{e}}\right) ^{2}\in \Q_{>0},\;$ then the symplectic class ${\frak{e}}^{s}$ is a
rational class and $\theta \left(  {{\frak{e}}}^{s}\right)=\gamma 
\left(  {\frak{e}}\right) ^{2}$.}

\medskip

\no{\sc Proof.} For a rational finite idele $j\in J_{\Q,f}$ we write $c\left(
j\right) $ for the positive rational number which generates the fractional $%
\Bbb{Z}$-ideal given by the content of $j.\;$The first part of the lemma
then follows from 3.2 on noting that the map from $J_{\Q,f}\times \R_{>0}\rightarrow \R_{>0}$ given by mapping $\left( j,r\right)
\longmapsto c\left( j\right) r^{-1}\;$has kernel $( \widehat{\Z}%
^{\times }\times 1)\cdot \Delta ( \Q^{\times }) .\;$\ To
show the second part of the lemma we first note that the class $\frak{e\;}$%
is represented by $1\times \gamma \left({\frak{e}}\right)
^{-1}\;$ and that the symplectic characters of the trivial group are the
even multiples of the trivial character. It therefore follows that \textit{$%
{\frak{e}}^{s}\;$}is represented by $1\times \gamma \left( 
{\frak{e}}\right) ^{-2},$ which has the same class in $A_{T}^{s}\left( \Z
\right) $ as $\gamma \left( {\frak{e}}\right)^{2}\times 1.$ \ \ $%
\Box $

\medskip

\bigskip We denote by $\Y_{S}^{\text{red}}$ the disjoint union of
the reduced fibres of $\mathcal{Y}$ over $p\in S.$ Let $\Y_{i},$
for $i\in {\cal I}$, denote the irreducible components of $\Y
_{S}^{\text{red}},$ so that 
\[
\Y_{S}^{\text{red}}=\cup _{i\in \mathcal{I}}\Y_{i}. 
\]
Let $p_{i}$ denote the prime which supports $\Y_{i}$ and let $\chi_{c}( \Y_{i}^{\ast })$
denote the $\ell $-adic Euler
characteristic with compact supports of $\Y_{i}^{\ast }=
\Y_{i}-\cup _{j\neq i}\Y_{j}$, the non-singular part of $\Y_{i}$.

Thanks to Theorem 7.5, in order to prove Theorem 7.1, we need only show
that, with the notation of 4.2, ${\frak{c}}_0^{s}$ is a rational class and
that 
\[
\theta ({\frak{c}}_{0}^{s})=\xi _{S}\left( 2\right) \tprod_{v\in S^{\prime
}}\varepsilon_{0,v}\left( \Y,2\right) ^{-1}. 
\]
Therefore, by 7.7, it will suffice to show that 
\[
\gamma \left( {\frak{c}}_{0}\right) ^{2}=\xi _{S}\left( 2\right) \tprod_{v\in
S^{\prime }}\varepsilon _{0,v}\left( \Y,2\right) ^{-1}\in \Q_{>0}. 
\]
>From (17), we know that for all non-negative $j$ 
\[
\left( \bigwedge ^{j}\omlogX \right) ^{G}=\bigwedge ^{j}\omlogY
\]
hence 
\[
{\frak{c}}_0=\tprod_{j=0}^{d}\chi \left( \Rr\Gamma \wedge ^{j}\omlogY 
 , \wedge ^{j}h_{Q,1}^{D}\right) ^{\left( -1\right) ^{j}} 
\]
which we write more succinctly as 
\[
{\frak{c}}_0=\chi \left( \Rr\Gamma \wedge ^{\bullet } \omlogY , \wedge
^{\bullet }h_{Q,1}^{D}\right) . 
\]
We therefore see that it is enough to show the following two
results:
\smallskip

\no{\bf Theorem 7.8} \ {\sl
\[
\gamma \circ \chi \left( \Rr\Gamma \wedge ^{\bullet } \omlogY , \wedge
^{\bullet }h_{Q,1}^{D}\right) =\tprod_{i\in {\cal I}
}p_{i}^{-\left( m_{i}-1\right) \chi_c \left( \Y_{i}^{\ast }\right) } 
\]
for any K\"{a}hler metric $j$ on the complex tangent
bundle $TY$.}

\medskip

\no{\bf Remark} In fact this result can also easily be proved using the
arithmetic Riemann-Roch Theorem of Gillet and Soul\'{e}; this alternative
approach to the calculation of Arakelov-Euler characteristics is explored in
[CPT2]; here, however, we shall provide a direct proof, which is due to
Bismut and which was shown to us by C. Soul\'{e}.\medskip
\medskip

\no{\bf Theorem 7.9 } 
\ \ \ \ \ $
\varepsilon_{0,S}\left( \Y,2\right) =\xi _{S}\left( 2\right)
\tprod_{i\in \mathcal{I}}p_{i}^{\left( m_{i}-1\right) \chi_c \left( \Y_{i}^{\ast }\right) }\in \Q_{>0}. 
$

\medskip

We begin by proving Theorem 7.9. For a place $v$ of $\Q$ we calculate
the $\varepsilon _{0,v}$-constants with respect to the standard Haar
measures $dx_{v}$ of $\Z_{v}$ and with respect to the Tate-Iwasawa
additive character $\psi _{v}$ of $\Q_v$ (see [Ta] p. 316-319).

We first consider the case of a finite prime $p$. From Theorem 2 in [Sa] we
know that $$\varepsilon _{0,p}\left( 2,\Y,\psi_{p}\circ
p^{-1},\;pdx_{p}\right) =\pm \tprod_{p_{i}=p}p^{\left( m_{i}-1\right) \chi
_{c}^{\ast }\left( \Y_{i}\right) }.$$ Thus by the standard
transformation formulae for $\varepsilon $-constants (see 5.3 and 5.4 in
[D]) 
\[
\varepsilon_{0,p}\left( 2,\Y,\psi _{p},dx_{p}\right) =\pm \sigma
^{2}\left( p\right) \tprod_{p_{i}=p}p_{i}^{\left( m_{i}-1\right) \chi
_{c}^{\ast }\left( \Y_{i}\right) }
\]
where $\sigma $ denotes the determinant of the motive of $\X\otimes_{G}V$ and where $V$ denotes the trivial representation of $G$. From
Proposition 2.2.1.a,c in [CEPT1] we know that $\sigma ^{2}\left( p\right)
=p^{-d\chi ( Y_{\Q}) }.$

Next, we consider the archimedean prime $v=\infty .$ From Lemma 5.1.1 in
[CEPT1], we know that $\varepsilon _{0,\infty }\left( 2,\Y,-\psi
_{\infty },dx_{\infty }\right) =\pm 1.$

To show that $\varepsilon_{0,S}\left( \Y,2\right) $ is positive,
note that from (2.2) in 2.4 of [CEPT1] we know that in all cases 
\[
{\rm sign}\left( \varepsilon _{0,v}\left( \Y,2\right) \right)
=\det \left( \sigma \right) \left( -1_{v}\right) .
\]
Thus by global reciprocity $1=$ $\tprod_{v\in S^{\prime }}\det \left( \sigma
\right) \left( -1_{v}\right) $ and so we have indeed now shown that $%
\varepsilon_{0,S}\left( \Y,2\right) $ is a positive rational
number.\ $\;\Box$
\medskip

Prior to proving Theorem 7.8, we note that we have:\smallskip

\no{\bf Lemma 7.10} \ {\sl Writing $\omega _{\Y/\Z}$ for the canonical sheaf of 
$\Y/\Z$,
there is a natural isomorphism between $\wedge^{d}\omlogY$ and 
$\omega _{\Y/\Z}( \Y_{S}^{\text{red}}-\Y_{S})$.}
\medskip

\no{\sc Proof.} Recall that $\{ \Y_{i}\} _{i\in \mathcal{I}}$
denote the irreducible components of  the disjoint union of the special
fibres $\Y_{S}^{\text{red}}.\;$From Proposition 3.1 in [CPT2] we
know that the natural morphism 
\[
\omega :\Omega _{\Y/\Z}\rightarrow 
\omlogY
\] 
has the same kernel and cokernel as the natural map 
\[
a:\oplus _{p\in S}\Oo_{\mathcal{Y}}/p\Oo_{\mathcal{Y}}\rightarrow \oplus _{i\in 
\mathcal{I}}\Oo_{\Y_{i}}. 
\]
The result then follows on taking determinants, since $\omega _{\Y/
\Z}\cong \det \Omega_{\Y/\Z}.$ $\ \;\;\Box $

\medskip

\no{\sc Proof of Theorem 7.8}. For brevity we regard the isomorphic degree
map $\gamma $ of 7.7 as an identification and we again put $\E=\omlogY$.

For $0\leq n\leq d$, the Duality Theorem in Section 11 of [H] gives a
quasi-isomorphism of complexes 
\[
\Rr\Gamma \left(\Hom_{\Oo_{\mathcal{Y}}}\left( \wedge ^{n}\E
,\omega _{\Y/\Z}[ d] \right)\right) \cong {\Hom}_{\Z}
\left( \Rr\Gamma(\wedge^{n}\E), \Z\right) 
\]
and by standard Hodge theory we know that the induced isomorphisms on
complex cohomology are isometries when the complex cohomology groups are
endowed with their $L^{2}$-metrics. Thus we see that 
\[
\chi _{L^{2}}\left( \Rr\Gamma \left(\Hom_{\Oo_{\mathcal{Y}}}( \wedge
^{n}\E,\omega _{\mathcal{Y}/\Bbb{Z}}[ d] )\right) \right)
=\chi _{L^{2}}\left( \Hom_{\Bbb{Z}}\left( \Rr\Gamma \wedge ^{n}%
\E, \Z\right) \right) 
\]
\begin{equation}
=\chi_{L^{2}}\left( \Rr\Gamma \wedge ^{n}\E\right)^{-1}.
\end{equation}
where for brevity we write $\chi _{L^{2}}\left( \Rr\Gamma \wedge ^{n}\E\right) $ in place of $\chi \left(\Rr
\Gamma \wedge ^{n}\E, \left\| {\ }\right\| _{L^{2}}\right) .$

Next we observe that by Lemma 7.10, we know that $\Hom_{\Oo_{\mathcal{Y}
}}\left( \wedge ^{n}\E, \omega _{\mathcal{Y}/\Bbb{Z}}( \Y_{p}^{\text{red}}-\Y_{p}) \right) $\linebreak $\cong
\wedge ^{d-n}\mathcal{E}.$ Thus we obtain a quasi-isomorphism 
\[
\Rr\Gamma (\wedge ^{d-n}\E)\cong \Rr\Gamma \left(\Hom_{\Oo_{%
\mathcal{Y}}}( \wedge ^{n}\mathcal{E}, \omega _{\mathcal{Y}/\Bbb{Z}}( \Y_{S}^{\text{red}}-\Y_{S}))
\right) 
\]
and, again by standard Hodge theory, we know that the induced isomorphisms
on complex cohomology are isometries with respect to their $L^{2}$-metrics.\
Thus we can write the number $\chi _{L^{2}}\left( \Rr\Gamma \left( \wedge
^{\bullet }\E\right) \right) ^{2}$ as: 
\[
\tprod_{n=0}^{d}\left[ \chi _{L^{2}}\left( \Rr\Gamma \left( \wedge ^{n}%
\E\right) \right) ^{\left( -1\right) ^{n}}\cdot \chi _{L^{2}}\left( 
\Rr\Gamma \left(\Hom_{\Oo_{\mathcal{Y}}} ( \wedge ^{n}\E,\omega _{%
\mathcal{Y}/\Bbb{Z}}( \Y_{S}^{\text{red}}-\Y_{S}) ) \right)\right) ^{\left( -1\right) ^{d-n}}\right] . 
\]
But this latter product can be rewritten as $\Pi _{1}\cdot \Pi _{2}$ where $\Pi
_{1}$ resp. $\Pi _{2}$ is the first resp. second of the following
expressions: 
\[
\tprod_{n=0}^{d}\left[ \chi _{L^{2}}\left( \Rr\Gamma \left( \wedge ^{n}%
\E\right) \right) ^{\left( -1\right) ^{n}}\cdot\chi _{L^{2}}\left( \Rr\Gamma \left(\Hom_{O_{\mathcal{Y}}}( \wedge ^{n}\E,\omega _{%
\mathcal{Y}/\Bbb{Z}}[ d]) \right) \right) ^{\left( -1\right) ^{n}}%
\right] 
\]
\[
\tprod_{n=0}^{d}\left[ \chi _{L^{2}}\left( \Rr\Gamma \left(\Hom_{\Oo_{%
\mathcal{Y}}}( \wedge ^{n}\mathcal{E},\omega _{\mathcal{Y}/\Bbb{Z}%
})\right) \right) ^{-1}\cdot\chi_{L^{2}}\left( \Rr\Gamma \left(\Hom_{\Oo_{%
\mathcal{Y}}}( \wedge ^{n}\E, \omega _{\Y/\Z
}( \Y_{S}^{\text{red}}-\Y_{S}))\right) \right)
 \right] ^{\left( -1\right) ^{d+n}} 
\]
and we note that (27) implies that $\Pi _{1}=1.$ Hence we may conclude that 
$\chi_{L^{2}}( R\Gamma\wedge^{\bullet}\E)
^{2} $ is equal to $\Pi _{2}.$ In order to evaluate $\Pi _{2}$ we consider
the exact sequences 
\[
0\rightarrow \omega _{\mathcal{Y}/\Bbb{Z}}( \Y_{S}^{\text{red}}-\Y_{S}) 
\rightarrow \omega _{\mathcal{Y}/\Bbb{Z}%
}( \Y_{S}^{\text{red}}) \rightarrow \omega _{\mathcal{Y}/%
\Bbb{Z}}( \Y_{S}^{\text{red}})| _{\mathcal{Y}%
_{S}}\rightarrow 0 
\]
\[
0\rightarrow \omega _{\mathcal{Y}/\Bbb{Z}}\rightarrow \omega _{\mathcal{Y}/%
\Bbb{Z}}( \Y_{S}^{\text{red}}) \rightarrow \omega _{%
\Y/\Z}( \Y_{S}^{\text{red}})|_{%
\Y_{S}^{\text{red}}}\rightarrow 0 
\]
and we apply the exact functor $\Hom_{\Oo_{\mathcal{Y}}}\left( \wedge ^{n}%
\E,-\right) $ to get exact sequences 
\[
0\rightarrow \Hom_{\Oo_{\mathcal{Y}}}( \wedge ^{n}\E%
,\omega _{\mathcal{Y}/\Bbb{Z}}( \Y_{S}^{\text{red}}-
\Y_{S}) ) \rightarrow \Hom_{\Oo_{\mathcal{Y}%
}}( \wedge ^{n}\E,\omega_{\mathcal{Y}/\Bbb{Z}}(
\Y_{S}^{\text{red}})) 
\]
\[
\rightarrow \Hom_{\Oo_{\mathcal{Y}}}( \wedge ^{n}\E,\omega
_{\mathcal{Y}/\Bbb{Z}}( \Y_{S}^{\text{red}})|_{%
\Y_{S}}) \rightarrow 0 
\]
and 
\[
0\rightarrow \Hom_{\Oo_{\mathcal{Y}}}(\wedge ^{n}\E,\omega _{%
\mathcal{Y}/\Bbb{Z}})\rightarrow \Hom_{\Oo_{\mathcal{Y}}}(\wedge ^{n}%
\E,\omega _{\mathcal{Y}/\Bbb{Z}}( \Y_{S}^{\text{red}%
}) ) 
\]
\[
\rightarrow \Hom_{\Oo_{\mathcal{Y}}}(\wedge ^{n}\E,\omega _{%
\mathcal{Y}/\Bbb{Z}}( \Y_{S}^{\text{red}})|_{\Y_{S}^{\text{red}}})\rightarrow 0. 
\]
Recall that$\;h$ denotes the structure map $h: \Y\rightarrow\Spec(\Z)$, $m_{i\text{ }}$
 denotes the multiplicity of the
component $\Y_{i}$  in $\Y_{S}$ and, as previously, for
each $i$ we let $p_{i}$ denote the prime which supports $\Y_{i}$. 
For brevity we shall write Hom$_{\Oo_{\mathcal{Y}}}( \wedge^\bullet 
\E, \omega _{\mathcal{Y}/\Bbb{Z}}) $ for $\tsum_{n}(
-1) ^{n}\Hom_{\Oo_{\mathcal{Y}}}( \wedge ^{n}\E,\omega _{%
\mathcal{Y}/\Bbb{Z}})$ etc. It then follows from the above and from
5.5 and 5.6 that $\Pi _{2}$ is equal to 
\[
\nu \circ h_{S\ast }\left( \Hom_{\Oo_{\mathcal{Y}}}(\wedge^\bullet
\E,\omega _{\mathcal{Y}/\Bbb{Z}}( \Y_{S}^{\text{red}%
})|_{\Y_{S}^{\text{red}}})-\Hom_{\Oo_{\mathcal{Y}%
}}(\wedge^\bullet\E,\omega _{\mathcal{Y}/\Bbb{Z}}( \Y
_{S}^{\text{red}})|_{\Y_{S}})\right) 
\]
\[
=\tprod_{i\in \mathcal{I}}p_{i}^{-\left( m_{i}-1\right) \left( -1\right)
^{d}\left( c^{d}\left( \E\right) \cdot \Y_{i}\right) }. 
\]
However, from 3.7 in [CPT2] (or see 5.1 in [CEPT2]), we know that 
\[
\left( -1\right) ^{d}c^{d}\left( \E\right) \cdot \Y_{i}=\chi
_{c}( \Y_{i}^{\ast }) 
\]
and so we have now shown 
\[
\chi _{L^{2}}(\Rr\Gamma \wedge ^{\bullet }\E)
^{2}=\tprod_{i\in \mathcal{I}}p_{i}^{-\left( m_{i}-1\right) \chi _{c}( 
\Y_{i}^{\ast }) }. 
\]

Finally we need to allow for the fact that in the above we have used the $%
L^{2}$-metric instead of the given K\"{a}hler metric. From the very
definition of the Quillen metric, we know that 
\[
\log \chi ( \Rr\Gamma \wedge ^{n}\E,\wedge ^{n}j^{D})
=\log \chi _{L^{2}}( \Rr\Gamma \wedge ^{n}\E) +\tau(
\wedge ^{n}\Omega _{Y},\wedge ^{n}j^{D}) 
\]
where $\tau \left( \wedge ^{n}\Omega _{Y},\wedge ^{n}j^{D}\right) $ denotes
the analytic torsion associated to $\wedge ^{n}\Omega _{Y}$ with respect to
the metric$\;\wedge ^{n}j^{D}.\;$But Theorem 3.1 in [RS] shows that 
\begin{equation}
\tsum_{n}\left( -1\right) ^{n}\tau \left( \wedge ^{n}\Omega_{Y},\wedge
^{n}j^{D}\right) =0\;\;
\end{equation}
and so we have now shown 
\[
\chi ( \Rr\Gamma\wedge ^{\bullet }\E, \wedge ^{\bullet
}j^{D}) =\tprod_{i\in \mathcal{I}}p_{i}^{-\left( m_{i}-1\right) \chi
_{c}( \Y_{i}^{\ast }) }. 
\]
This then completes the proof of Theorem 7.8.\ \ \ \ $\Box \medskip $

Observe that Theorems 7.8 and 7.9 show that 
\[
\gamma \circ \chi \left( \Rr\Gamma ( \wedge ^{\bullet }\omlogY )
,\wedge ^{\bullet }j_{Q}^{D}\right) ^{2}=\xi _{S}\left( 2\right) \varepsilon
_{0,S}( \Y, 2) ^{-1}. 
\]
We conclude this section by showing that the right hand factor $\xi
_{S}\left( 2\right) $ in the above can be removed by twisting the sheaf 
$\omlogY$ by $O_{\mathcal{Y}}( -\Y_{S})$.
\medskip
\medskip
 
\no{\bf Theorem 7.11} \ \ 
$
\gamma \circ \chi \left( \Rr\Gamma ( \wedge ^{\bullet }\omlogY (-\Y_S))
,\wedge ^{\bullet }j_{Q}^{D}\right) ^{2}=\varepsilon
_{0,S}( \Y, 2) ^{-1}.
$
\medskip
\medskip

\no{\sc Proof.} Since for each $i\geq 0\;$ 
\[
\wedge ^{i}(\omlogY(-\Y_S))=\wedge ^{i}\omlogY\otimes\Oo_{\Y}(-i\Y_S)
\] 
we obtain an exact sequence of complexes of sheaves 
\[
0\rightarrow \wedge ^{\bullet }(\omlogY(-\Y_S)) \rightarrow \wedge ^{\bullet }
\omlogY  \rightarrow \G^{\bullet}\rightarrow 0
\]
where for $0\leq i\leq d$
\[
\G^{i}=\wedge ^{i}\omlogY|_{i\Y_{S}}
\]
and so by (5.5), (5.6) and the equality displayed prior to (7.11) 
\[
\gamma \circ \chi \left( R\Gamma ( \wedge ^{\bullet }\omlogY(-\Y_S))  ,\wedge ^{\bullet }j_{Q}^{D}\right) ^{2}=\xi
_{S}( 2)\cdot \varepsilon_{0,S}(\Y,2)^{-1}\cdot \chi( \nu ( \G^{\bullet })) ^{2}
\]
and for $0\leq i\leq d$%
\[
\chi ( \nu ( \G^{i}))^{2}=\tprod_{p\in
S}p^{2f_{\ast }( \G^{i}) }=\tprod_{p\in S}p^{2i\chi( \Omega^i_{\Y_{\Q}})}
\]
since $\omlogY_{\Bbb{Q}}=\Omega_{\Y_{\Q}}.$ However, by
Serre duality we know that $\left( -1\right) ^{d-i}\chi ( \Omega _{%
\mathcal{Y}_{\Bbb{Q}}}^{d-i}) =( -1) ^{i}\chi ( \Omega
_{\mathcal{Y}_{\Bbb{Q}}}^{i}) $ and so we see that
\[
\tsum_{i=0}^{d}( -1) ^{i}i\cdot\chi ( \Omega _{\mathcal{Y}_{\Bbb{Q%
}}}^{i}) =\tsum_{i=0}^{d}( -1) ^{d-i}i\cdot\chi ( \Omega _{%
\mathcal{Y}_{\Bbb{Q}}}^{d-i}) =\tsum_{i=0}^{d}( -1)
^{i}( d-i)\cdot \chi ( \Omega _{\mathcal{Y}_{\Bbb{Q}}}^{i}) 
\]
hence 
\[
\tsum_{i=0}^{d}( -1) ^{i}2i\cdot\chi ( \Omega _{\mathcal{Y}_{\Bbb{%
Q}}}^{i}) =d\cdot\chi ( \mathcal{Y}_{\Bbb{Q}}) 
\]
which therefore shows that 
\[
\tprod_{i}\chi ( \nu ( \G^{i}))^{2(
-1) ^{i}}=\tprod_{p\in S}p^{2i( -1) ^{i}\chi( \Omega
_{\mathcal{Y}_{\Bbb{Q}}}^{i}) }=\tprod_{p\in S}p^{d\chi ( 
\Y_{\Bbb{Q}}) }=\xi _{S}( 2) ^{-1}
\]
as required.\ \ \ $\Box $

\bigskip

\section{Differentials}

 In this the final section of the article we construct arithmetic
classes associated to the sheaf of (regular) differentials $\Omega^1_{%
\mathcal{X}\text{/}\Bbb{Z}}$. Since $\Omega^1_{\mathcal{X}\text{/}\Bbb{Z}}$
is not in general locally free over $O_{\mathcal{X}},$ we resolve it by
locally free $G$--$\X$ sheaves as follows: we choose a $G$%
-equivariant embedding $i:\X\rightarrow {\cal P}$ of $\mathcal{X}$ into$%
\mathcal{\;}$a projective bundle $\mathcal{P\;}$over Spec$\left( \Bbb{Z}%
\right) $. The sheaf of differentials $\Omega^1_{\mathcal{X}\text{/}\Bbb{Z}}$
then has a resolution by locally free $G$--$\X$ sheaves 
\begin{equation}
\;0\rightarrow N^{\ast }\rightarrow P\stackrel{\pi }{\rightarrow }\Omega^1_{%
\mathcal{X}\text{/}\Bbb{Z}}\rightarrow 0
\end{equation}
where $P=i^{\ast }\Omega^1_{\mathcal{P}\text{/}\Bbb{Z}},$ and where $N^{\ast }
$ denotes the conormal bundle associated to the regular embedding $i.$ Let $%
\mathcal{F}^{\bullet }$ denote the length two complex 
\[
\F^{\bullet }\ :\ N^{\ast }\rightarrow P
\]
where the term $P$ is deemed to have degree zero. Thus we may view $\pi $ as
inducing a quasi-isomorphism of complexes, which we abusively also denote $%
\pi ,$ 
\[
\pi \ :\ \F^{\bullet }\rightarrow \Omega^1_{\mathcal{X}\text{/}\Bbb{Z}}.
\]
Here we further abuse notation and write $\Omega^1_{\mathcal{X}\text{/}\Bbb{Z}%
}$ for the complex which is $\Omega^1_{\mathcal{X}\text{/}\Bbb{Z}}$ in degree
zero and which is zero elsewhere.

For $j\geq 0,$ recall that we have the Dold-Puppe exterior power functors $%
\bigwedge^{j}\;$defined on bounded complexes of locally free $G$--$\X$
sheaves and which take quasi-isomorphisms to quasi-isomorphisms. (See [CPT3]
for an account of these functors which is particularly well-suited to their
use in this paper.)

We then endow the equivariant determinant of cohomology of the complex $%
\wedge ^{j}\left( \mathcal{F}^{\bullet }\right) $ with the metrics $\phi
_{j\bullet }$ induced, via $\wedge ^{j}( \pi _{\Bbb{C}}) ,$ from
the $\wedge ^{j}h_{Q\bullet }^{D}$ on the determinants of cohomology of $%
\Omega _{X\text{/}\Bbb{C}}^{j};\;$we then define arithmetic classes 
\begin{equation}
\chi \left( \Rr\Gamma \lambda^{j}\Omega^1_{\mathcal{X}\text{/}\Bbb{Z}%
},\wedge ^{j}h_{Q}^{D}\right) :=\chi \left( \Rr\Gamma \wedge^{j}( 
\F^{\bullet }) ,\phi _{j}\right)
\end{equation}

\begin{equation}
\chi \left( \Rr\Gamma \lambda ^{\bullet }\Omega^1_{\X/\Z},\wedge^{\bullet }h_{Q}^{D}\right) :=\tprod_{j=0}^{d}\chi \left( \Rr\Gamma \lambda ^{j}\Omega _{\X/\Z},\wedge
^{j}h_{Q}^{D}\right) ^{\left( -1\right) ^{j}}.
\end{equation}
Note that here the use of the symbols $\lambda^{j\;}\;$is entirely
symbolic; however, it is important to\ observe that the lefthand classes are
independent of the chosen embedding $i:\mathcal{X\rightarrow P}$: indeed,
for a further embedding $i^{\prime },$ with the obvious notation, $\wedge
^{j}( \mathcal{F}^{\prime \bullet }) $ is quasi-isomorphic to $%
\wedge ^{j}( \mathcal{F}^{\bullet }) ;$ furthermore their metrics
on the determinant of cohomology match under the corresponding
quasi-isomorphism; hence by 3.9 the arithmetic classes coincide.

The equivariant Arakelov-Euler characteristic $\chi ( \Rr\Gamma
\lambda ^{\bullet }\Omega^1_{\mathcal{X}\text{/}\Bbb{Z}},\wedge ^{\bullet
}h_{Q}^{D}) $ is the principal object of study in this section. Our
aim here is to relate it to the epsilon constant  $\varepsilon( 
\mathcal{Y}) ,$ whose definition we now briefly recall. Let $\Bbb{A}_{%
\Bbb{Q}}$ denote the ring of rational adeles; let $\psi =\tprod_v \psi _{v}$
denote a non-trivial additive character of $\Bbb{A}_{\Bbb{Q}}/\Bbb{Q}$;$\Bbb{%
\;}$let $dx$ denote the Haar measure on $\Bbb{A}_{\Bbb{Q}}/\Bbb{Q}$ such
that $\int_{\Bbb{A}_{\Bbb{Q}}/\Bbb{Q}}dx=1$ and let $dx=\tprod_v dx_{v}$ be a
factorisation of $dx\,$\ into local Haar measures $dx_{v}\;$with the
property that $\int_{\Z_{v}}dx_{v}=1$ for almost all $v$. Recall from
3.1.1 in [CEPT1] that for $\theta \in R_{G},$%
\[
\varepsilon_{v}\left( \Y, \theta ,\psi
_{v}, dx_{v}, l_{v}\right) =\varepsilon_{0,v}\left( \Y, \theta
, \psi _{v}, dx_{v}, l_{v}\right) \varepsilon \left( \Y_v,
\theta \right) .
\]
Here if $v<\infty$ then $\varepsilon \left( \Y_v, \theta
\right)$ is the epsilon constant associated to the special fibre $\Y_v$ and if $v=\infty$ then we take $\varepsilon \left( \Y_v,
\theta \right)=1$.  We then set 
\[
\varepsilon \left( \mathcal{Y}\text{,}\theta \right) =\tprod_{v}\varepsilon
_{v}\left(\Y_v, \theta ,\psi _{v},dx_{v},l_{v}\right) .
\]
Note that in this product almost all terms are 1 and moreover this product
is independent of choices of additive character and Haar measure.\ Thus in
the lefthand term we shall abuse notation and henceforth we shall not
overtly mention the choices of auxiliary primes $l_{v}.$

For future reference we now need to gather together some standard results
on fibral epsilon constants.

For this we require a minor variant on the notation introduced prior to
Theorem 5.10. As previously, given a prime number $p$ we fix a field
embedding $h:\overline{\Bbb{Q}}\rightarrow \overline{\Bbb{Q}}_{p},$ we put $%
( \overline{\Bbb{Q}})_{p}=\overline{\Bbb{Q}}\otimes \overline{%
\Bbb{Q}}_{p}$ and we let $J_{f}\rightarrow \left( \overline{\Bbb{Q}}\right)
_{p}^{\times }$ denote the map given by projection to the $p$th coordinate$%
;\;$given $x\in J_{f},$ we shall write $x_{p}$ for the $p$-component of $x$
in $( \overline{\Bbb{Q}}) _{p}^{\times }$.\ Let $\left| -\right|
_{p}:\overline{\Bbb{Q}}_{p}^{\times }\rightarrow p^{\Bbb{Q}}$ denote the $p$%
-adic absolute value which is normalised so that $\left| p\right|
_{p}=p^{-1}.$ We shall use the terminology of Definition 5.6 in [C] and for $%
f\in \Hom_{\Omega }( R_{G},( \overline{\Bbb{Q}})
_{p}^{\times }) $ we say that $\left| f\right| _{p}$ is \textit{%
well-defined} if $\left| h^{\ast }\left( f\right) \right| _{p}$ takes values
in $p^{\Bbb{Z}};$ in this case it follows that $\left| h^{\ast }\left(
f\right) \right| _{p}$ respects $\Omega _{p}$-action and we then write 
\[
\left| f\right| _{p}=h^{\ast -1}\left| h^{\ast }f\right| _{p}. 
\]
\smallskip

\no{\bf Theorem 8.1} \ {\sl For each prime number $p$, $\left|
\varepsilon ( \Y_p) \right|_{p}$ and 
$\left| \varepsilon ( b_{p}) \right|_{p}$ are
well-defined. Writing $U_{p}$ for the open sub-scheme 
$Y_{p}-b_{p}$, we have $\varepsilon ( \Y_{p}) =\varepsilon ( U_{p}) \varepsilon ( b_{p}) 
$ and
\[
\varepsilon  ( U_{p}) _{p}\left| \varepsilon  ( U_{p} )
\right| _{p}\in \Det( \Z_p[G]^{\times }) 
\]
 and for a prime number $q\neq p$ 
\[
\varepsilon  ( U_{p} ) _{q}\in \Det( \Z_{q}[G]^{\times }) . 
\]}
\medskip

\no{\sc Proof.} See [C] 5.7, 5.13 and 5.12.\ \ \ \ $\Box \medskip $

In order to make precise the fundamental relationship between 
$\chi ( \Rr\Gamma \lambda ^{\bullet }\Omega^1_{\mathcal{X}\text{/}%
\Bbb{Z}},\wedge ^{\bullet }h_{Q\bullet }^{D} ) $ and $\varepsilon
\left( \mathcal{Y}\right) ,$ we now need to introduce the \textit{arithmetic
ramification class}, which may be viewed as an arithmetic counterpart of the
ramification class occurring in Theorem 1.1 in [CPT1].\smallskip \smallskip
\medskip

\no{\bf Definition 8.2} \ Let ${\rm AR}(\X) \in A(\Z[G])$ 
be the arithmetic class which is represented by the idele valued
character function $\beta$, given by the rule that $\beta $ has trivial
archimedean coordinate and at a finite prime $q$%
\[
\beta _{q}=\varepsilon  ( b) \left| \varepsilon  (
b_{q} ) \right| _{q} 
\]
where $b_{q}$ denotes the union of the components of $b$ which are supported
by $q.\smallskip $

We are now in a position to be able to state the main result of this
article: 
\medskip

\no{\bf Theorem 8.3}  \ {\sl Let $\frak{d}$ be the arithmetic class 
$\chi (\Rr\Gamma \lambda ^{\bullet }\Omega^1_{\mathcal{X}\text{/}%
\Bbb{Z}},\wedge ^{\bullet }h_{Q}^{D})$. Then ${\frak{d}}^{s}\cdot {\rm AR}^{s} 
( \X) ^{-1}$  is a rational
class and 
\[
\theta \left( {\frak{d}}^{s}\cdot {\rm AR}^{s}(\X )
^{-1}\right) =\varepsilon ^{s}(\Y) ^{-1}. 
\]}
\medskip

As a first step towards the proof of this theorem, we use results from
[CPT2] to show that it will suffice to establish the corresponding result
after passage to degree zero: 
\medskip

\no{\bf Theorem 8.4} \  {\sl The class $\widetilde{\frak{d}}^{s}\cdot
\widetilde{{\rm AR}}^{s}(\X)$ is a rational
class and
\[
\theta \left( \widetilde{\frak{d}}^{s}\cdot \widetilde{\text{AR}}^{s}
(\X)^{-1}\right) =\widetilde{\varepsilon }^{s}(\Y)^{-1}. 
\]}
\smallskip

We begin by showing that Theorem 8.4 implies Theorem 8.3; we then conclude
the article by establishing Theorem 8.4.
\smallskip

Suppose then that Theorem 8.4 holds. From 4.E, with the notation of 4.2, we
know that 
\[
{\frak{d}}=\widetilde{\frak{d}}\cdot {\rm Ind}( {\frak{d}}_{0}) . 
\]
By 5.7 in [C] we know that $\varepsilon( b_{p},1_{G})$ is $\pm $
an integral power of $p$; hence we see that $\beta $ and $\widetilde{\beta }$
represent the same class in $A(\Z[G]) $ and so $\widetilde{%
\text{AR}}^{s}(\X) ={\rm AR}^{s}(\X).$

\medskip

\no{\bf Theorem 8.5} \ \ {\sl  ${\frak{d}}_{0}^{s}$ is a rational class and \
$
\theta \left( {\frak{d}}_{0}^{s}\right) =\varepsilon ( \Y ,2\cdot 1_{G}) ^{-1}.
$}
\medskip

\no{\bf Proof.} We endow $P$ and $N^{\ast }$  in (29) with $G$-invariant
hermitian metrics and denote the resulting hermitian bundles by $\widehat{P}$
and $\widehat{N^{\ast }}$. We then let $\eta _{1}$ denote the Bott-Chern
class assocated to the exact sequence (29), where $\Omega^1_{Y}$ is endowed
with the hermitian metric $h^{D},$ and we put $\widehat{\Omega}=\widehat{P}-%
\widehat{N^{\ast }}+\eta _{1}$ in the arithmetic Grothendieck group $%
\widehat{\rm K}_{0}(\Y) $ (see for instance II,
Sect. 6 in [GS1]); we recall \textit{loc. cit.} that $\widehat{\rm K}_{0}(\Y)$
 has a natural structure of a $\lambda $-ring
and we write $\widehat{f}_{\ast }$ for the push forward map from  $\widehat{\rm K}_{0}(\Y)$ to  
$\widehat{\rm K}_{0}(\Spec(\Z))$. Because $\mathcal{Y}$ is regular, we know that $\widehat{\rm K}_{0}(\Y)$
  is naturally isomorphic to $\widehat{\rm K}'_{0}(\Y)$ the
Grothendieck group of coherent hermitian sheaves (see Lemma 13 in [GS2]).
Thus we also have a natural map from the Grothendieck group of torsion $%
\mathcal{Y}$-sheaves supported on $S,$ denoted  $\widehat{\rm K}^S_{0}(\Y)$, to $
\widehat{\rm K}_{0}(\Y)$.
 Recall
that  $\widehat{\rm K}_{0}^S(\Y)$ is a module over the Grothendieck group of locally free $\mathcal{Y}$%
-sheaves $\widehat{\rm K}_{0}(\Y)$.

In Theorem 1.3 of [CPT1], with slightly different notation, it is shown that 
\[
\tsum_{i=0}^{d} \gamma \circ \chi \circ \widehat{f}_{\ast }(
( -1)^{i}\lambda ^{i}\widehat{\Omega }) =\left|
\varepsilon (\Y, 1_{G}) \right| ^{-1}
\]
whereas the class that we now wish to study is 
\[
\gamma \circ \chi ( \Rr\Gamma \lambda ^{\bullet }\Omega^1_{\mathcal{%
Y}\text{/}\Bbb{Z}},\wedge ^{\bullet }h_{Q,1}^{D}) =\tsum_{i=0}^{d}%
\gamma \circ \chi \circ \widehat{f}_{\ast }( (
-1) ^{i}( \lambda ^{i}( \widehat{P}-\widehat{N^{\ast }}%
) +\eta _{2}^{\left( i\right) })) 
\]
where the $i$th exterior power $\wedge ^{i}\Omega_{Y}=\Omega^i_Y$ carries the
hermitian metric $\wedge ^{i}h^{D}$, the terms of $\wedge ^{i}\mathcal{F}_{%
\Bbb{C}}^{\bullet }$ carry the metrics coming from $\widehat{P}$ and $%
\widehat{N^{\ast }}$, and where $\eta _{2}^{\left( i\right) }$ is the
Bott-Chern class associated to the exact sequence of hermitian bundles $%
\wedge ^{i}\widehat{\mathcal{F}_{\Bbb{C}}}^{\bullet }\rightarrow $ $\wedge
^{i}\widehat{\Omega }_{Y}.$ As our first step in proving the theorem, we
will show that in $\widehat{\text{K}}_{0}\left( \mathcal{Y}\right) $ 
\begin{equation}
\lambda ^{i}( \widehat{P}-\widehat{N^{\ast }}) +\eta _{2}^{\left(
i\right) }=\lambda ^{i}( \widehat{P}-\widehat{N^{\ast }}+\eta
_{1}) 
\end{equation}
which will then imply that 
\begin{equation}
\gamma \circ \chi ( \Rr\Gamma \lambda ^{\bullet }\Omega^1_{\mathcal{%
Y}\text{/}\Bbb{Z}},\wedge ^{\bullet }h_{Q,1}^{D}) =
\left|
\varepsilon (\Y, 1_{G}) \right| ^{-1}
\end{equation}

>From Lemma 7.10 we know that there is an exact sequence 
\[
0\rightarrow K\oplus N^{\ast }\rightarrow P\rightarrow  \omlogY \rightarrow
C\rightarrow 0
\]
where $K\,$\ and $C$ explicitly determined torsion $\Oo_{\mathcal{Y}}$-modules
supported on $S$. Thus for each $i,0\leq i\leq d,$ we have an equality in $%
\widehat{\text{K}}_{0}^{\prime }(\Y)$%
\[
\lambda ^{i}( \widehat{P}-\widehat{N^{\ast }}-\widehat{K}+\eta
_{1}) =\lambda ^{i}\left( \widehat{\Omega}^1_{\Y/\Z
}( \log \Y_{S}^{\text{red}}/\log S) -\widehat{C}\right) 
\]
which we can rewrite as 
\[
\lambda ^{i}( \widehat{\Omega }) +T_{1}^{\left( i\right)
}=\lambda ^{i}\left( \widehat{\Omega}^1_{\Y/\Z}( \log 
\Y_{S}^{\text{red}}/\log S))  \right) +T_{2}^{\left( i\right) }
\]
where $T_{1}^{\left( i\right) }$ and $T_{2}^{\left( i\right) }$ are the
following torsion classes 
\[
T_{1}^{\left( i\right) }=\tsum_{a+b=i,\; b>0}\lambda ^{a}\left( P-N^{\ast
}\right) \lambda ^{b}\left( -K\right) ,\;
\]
\[
T_{2}^{\left( i\right) }=\tsum_{a+b=i,\;b>0}\lambda ^{a}\left( \widehat{\Omega}^1_{\Y/\Z}( \log 
\Y_{S}^{\text{red}}/\log S))  \right) \lambda ^{b}( -C) .
\]
Next we consider the quasi-isomorphism of the Dold-Puppe exterior powers (where 
$P$ and $\omlogY$ are both deemed to have degree zero) 
\[
\wedge ^{i}\left( K^{\bullet }\oplus N^{\ast }\rightarrow P\right) \cong
\bigwedge ^{i}\left( \omlogY  \rightarrow C^{\bullet }\right) 
\]
and where $K^{\bullet }$ and $C^{\bullet }$ denote locally free resolutions
of $\ K$ and $C.$ Hence filtering the complex $\bigwedge ^{i}\left( K^{\bullet }\oplus
N^{\ast }\rightarrow P\right)$, by terms $\bigwedge ^{a}\left( N^{\ast
}\rightarrow P\right) \otimes \bigwedge ^{i-a}\left( K^{\bullet }\left[ 1\right]
\right)$, and filtering the complex $\bigwedge ^{i}(\omlogY \rightarrow
C^{\bullet }) $ by terms $\bigwedge ^{a}\left( \Omega^1_{\Y/\Z
}\left( \log \Y_{S}^{\text{red}}/\log S\right) \right) \otimes
\wedge ^{i-a}\left( C^{\bullet }\left[ -1\right] \right)$ (see p. 26 in
[S]), we obtain an equality in $\widehat{\text{K}}_{0}^{\prime }(\Y)$
\[
T_{1}^{\left( i\right) }+\lambda ^{i}( \widehat{P}-\widehat{N^{\ast }}) 
+\eta _{2}^{\left( i\right) }=\lambda ^{i}\left( \widehat{\Omega}^1_{\Y/\Z}( \log 
\Y_{S}^{\text{red}}/\log S)  \right)
+T_{2}^{\left( i\right) }
\]
which now establishes (32).

\bigskip

Now each 
$\varepsilon(\Y, 1_G)$ is a rational number 
and so by Theorem 7.9,  $\varepsilon(\Y, 2\cdot 1_G)
  \in \Q_{>0}.$ Thus by Lemma 7.7 we
see that if we can show 
\begin{equation}
{\frak{d}}_{0}=\chi  ( \Rr\Gamma \lambda ^{\bullet }\Omega^1_{\mathcal{%
Y}\text{/}\Bbb{Z}},\wedge ^{\bullet }h_{Q,1}^{D}) 
\end{equation}
then it will follow that ${\frak{d}}_{0}$ is a rational class and that 
\[
\theta \left( {\frak{d}}_{0}^{s}\right) =\varepsilon(\Y, 2\cdot 1_G)^{-1}.
\]
With the notation of 7.6 above and by the very definition of ${\frak{d}}_{0}$
(see 4.2), 
\[
{\frak{d}}_{0}=\chi 
( \Rr\Gamma \lambda ^{\bullet } (\F^\bullet)^G,\wedge ^{\bullet }h_{Q,1}^{D}) .
\]
Thus we are now required to show that for each $j,0\leq j\leq d,\;$\ the
natural map $\wedge ^{j}( \F^{\bullet G}) $ to $\left(\wedge
^{j}\F^{\bullet }\right) ^{G}$ is a quasi-isomorphism.  To
see this it will suffice to show the result after passing to a 
flat neighbourhood of each closed point of\ $y$ of $\mathcal{Y}$.
Writing $\X^{\prime }\rightarrow \Y^{\prime }$ for the
resulting base change to such a neighbourhood, we let $x^{\prime }$ resp.
$y^{\prime }$ denote a closed point of $\mathcal{X}^{\prime }$ resp. $%
\mathcal{Y}^{\prime }\;$above $y.$ From Theorem A.1 and Lemma A.2 in
[CEPT1], we know that, for a suitable choice of neighbourhood, $\mathcal{X}%
^{\prime }\;$contains $\left( G:I_{x^{\prime }}\right) $ disjoint
irreducible components which are permuted transitively by $G$ and where the
component which contains $x^{\prime }$ has stabiliser $I_{x^{\prime }},$ the
inertia group of $x^{\prime }.$ If $B_{1},..,B_{q}$ are the distinct
irreducible components of the inverse image $\pi ^{-1}(b)$ which contain the
image of $x^{\prime }$ on $\mathcal{X},$ then 
\[
I_{x^{\prime }}=I_{1}\oplus ..\oplus I_{q}
\]
where $I_{i}$ denotes the inertia group of the generic point of $B_{i};$
moreover, each $I_{i}$ carries a faithful abelian character $\phi _{i}$
given by the action of $I_{i}$ on the cotangent space of the generic point
of $B_{i}$. To be somewhat more precise, there are integers $%
n_{1},..,n_{d+1}$ coprime to the residual characteristic of $y$ so that,
after base extension by a suitable affine flat neighbourhood Spec$\left(
R\right) ,$ the connected open neighbourhood $V$ of $\X^{\prime }=\X\times \Spec(T)$
containing $x^{\prime }$
is the spectrum of
\[
\frac{R\left[ U_1, \ldots, U_{d+1}\right] }
{\left(
U_{1}^{n_{1}}-a_{1},\ldots ,U_{1}^{n_{d+1}}-a_{d+1}\right) }.
\]
Here $a_{1},\ldots ,a_{d+1}$ form a system of regular parameters of $\Y'$;
moreover there are integers $m_{i}$ for $1\leq i\leq d+1$ with
each $m_{i}$ coprime to the residual characteristic, $p$ say, of $y,$ and
with the property that $a_{1}^{m_{1}}\cdots a_{d+1}^{m_{d+1}}=p.$ Here, after
reordering if necessary, the characters $\phi _{i}$ are given by the action
of $I_{i}$ on $U_{i}.$ It now follows that $\Omega^1_{V/R}$ sits in an exact sequence 
\[
0\rightarrow K^{\bullet }\rightarrow \Omega^1_{V/R} \rightarrow 0
\]
with 
\[
K^{\bullet }\ :\ \Oo_V dr\rightarrow \bigoplus
_{i=1}^{d+1}\Oo_{V} dU_{i}
\]
and where $r=a_{1}^{m_{1}}..a_{d+1}^{m_{d+1}}-p.$ In the sequel for brevity
we shall write $K^{\bullet }=L\rightarrow E.$ Since the restriction $\F^{\bullet
}|V$ is quasi-isomorphic to $K^{\bullet }$ we are now
reduced to showing that $\wedge ^{m}( K^{\bullet I}) \simeq
\left( \wedge ^{m}K^{\bullet }\right) ^{I}$ for all $m\geq 0$
and for $I=I_{x'}$. This now
follows easily since we know (see for instance Sect. 3 in [CPT3]) that the
complex $\wedge ^{m}K^{\bullet }\,\;$is constituted entirely of terms which
are tensor products of modules of the form $\wedge ^{n}L\,$\ times either
one or no terms of the form $\wedge ^{n}E;$ the result then follows because $%
L\cong \Oo_{V}$, as $I$-modules, and
because, for any non-negative $n$, $\wedge ^{n}( E^{I}) \cong
\left( \wedge ^{n}E\right) ^{I}$ (using the fact that the $\phi _{i}$ come
from the distinct components in a direct sum decomposition).\ \ \ $\Box $

\medskip

\medskip

\no{\sc Proof of Theorem 8.4.}\ We write $\X_{S}^{\text{red%
}}=\tbigsqcup_{p\in S}\X_{p}^{\text{red}}$, let $i_{S}^{\text{red}}:\X_{S}^{\text{red}}\rightarrow \X$
 denote the associated
closed embedding and we let $U_{S}$ denote the complement of $\X_{S}^{\text{red}}$ in $\mathcal{X}$.

Composing the quasi-isomorphism $\pi :\F^{\bullet }\simeq \Omega^1_{%
\mathcal{X}\text{/}\Bbb{Z}}$ with the natural homomorphism $\omega :\Omega^1_{\mathcal{%
X}\text{/}\Bbb{Z}}\rightarrow \omlogX $, which is an isomorphism
over $U_{S},$ we get a chain map 
\[
\pi ^{\prime }:\F^{\bullet }\rightarrow \omlogX  
\]
which is a surjective quasi-isomorphism over $U_{S}$. Hence for $i\geq 0$
we obtain maps 
\[
\wedge ^{i}\pi ^{\prime }:\wedge ^{i}\F^{\bullet }\rightarrow
\wedge ^{i}\omlogX  
\]
which are surjective quasi-isomorphisms over $U_{S}$. Let 
\[
{\cal A}_{i}^{\bullet }=\ker (\wedge ^{i}\pi ^{\prime })\;\;\;\;\hbox{\rm and\
\ \ \ }{\cal B}_{i}^{\bullet }={\rm coker}(\wedge ^{i}\pi ^{\prime })
\]
so that ${\cal B}_{i}^{\bullet }$ and the cohomology sheaves ${\cal H}^{j}
( {\cal A}_{i}^{\bullet })$ of the complex ${\cal A}_{i}^{\bullet }$
 are all supported entirely over $S$. Let $\frak I$
denote the ideal sheaf of $O_{\mathcal{X}}$ associated to the closed
subscheme $\mathcal{X}_{S}^{\text{red}}$; we then write $ \left[{\cal H}^{j}
( {\cal A}_{i}^{\bullet })  \right] $ for the finite sum 
\[
\tsum_{n\geq 0}\left( {\frak{I}}^{n}{\cal H}^{j}
( {\cal A}_{i}^{\bullet })/{\frak{I}}^{n+1}{\cal H}^{j}
( {\cal A}_{i}^{\bullet })\right)
\]
 in G$_{0}( G,\X_{S}^{\text{red}}) $ and put 
\[
\left[ {\cal H}^{\bullet}
( {\cal A}_{i}^{\bullet }) \right]
=\tsum_{j}\left( -1\right) ^{j}\left[ {\cal H}^{j}
( {\cal A}_{i}^{\bullet })\right] \ \in \ {\rm G}_{0}( G,\X_{S}^{%
\text{red}}) . 
\]

We endow the equivariant determinants of cohomology of ${\cal A}
_{i}^{\bullet }$ and ${\cal B}_{i}^{\bullet }$ with the trivial metrics $%
\tau _{\bullet }$. Then from 6.4 and 5.6 we know that 
\[
\widetilde{\chi }\left( \Rr\Gamma \wedge ^{i}\F^{\bullet
},\wedge ^{i}h^{D}\right)\cdot \widetilde{\chi }\left( R\Gamma \wedge
^{i}\omlogX ,\wedge ^{i}h^{D}\right) ^{-1} 
\]
\[
=\widetilde{\chi }\left( \Rr\Gamma {\cal A}_{i}^{\bullet },\tau
_{\bullet }\right)\cdot  \widetilde{\chi }\left( R\Gamma {\cal B}
_{i}^{\bullet },\tau _{\bullet }\right) ^{-1} 
\]
\begin{equation}
=\widetilde{\nu \circ f_{S\ast }^{\text{red}}}\left( 
\left[ {\cal H}^{\bullet}
( {\cal A}_{i}^{\bullet }) \right]-\left[ {\cal H}^{\bullet}
( {\cal B}_{i}^{\bullet }) \right]\right).
\end{equation}
\smallskip

Since the class $F=\left( -1\right) ^{d}\tsum_{i}\left( -1\right) ^{i}
\left( 
\left[ {\cal H}^{\bullet}
( {\cal A}_{i}^{\bullet }) \right]-\left[ {\cal H}^{\bullet}
( {\cal B}_{i}^{\bullet }) \right]\right)$
 in ${\rm G}_{0} ( G,\X_{S}^{\text{red}} ) $ has the
property that its image in ${\rm G}_{0}( G,\X) ={\rm K}_{0}( G,\X) $ 
\begin{equation}
i_{S\ast }^{\text{red}}F=\left( -1\right) ^{d}\tsum_{i}\left( -1\right)
^{i}\left( [ \wedge ^{i}\F^{\bullet }] -\left[ \wedge
^{i}\omlogX \right] \right)
\end{equation}
\[
= c^{d} ( \Omega _{\mathcal{X}\text{/}\Bbb{Z}} ) -%
 c^{d} ( \omlogX) 
\]
we may take $F=\oplus _{p\in S}F_{p}$ to be the class $F$ in 6(a) of
[CPT1].
\medskip

\no{\bf Definition 8.6} \ For each $i\in \mathcal{I}$ we set $\X
_{i}=\pi ^{-1}( \Y_{i}) ;$ thus $\X_{i}$ is a
smooth projective variety over ${\bf F}_{p_{i}}$ of dimension $d$ which
carries a tame $G$-action. More generally for each non-empty subset $%
\mathcal{J}$ of $\mathcal{I}$ we define 
\[
\Y_{\cal J}=\cap _{j\in \mathcal{J}}\Y_{j},\qquad \X_{\mathcal{J}}=\cap _{j\in \mathcal{J}}\X_{j} 
\]
so that each $\mathcal{X}_{\mathcal{J}}$ is either empty or is a smooth
projective variety of dimension 
$d+1-\left| \mathcal{J}\right| .$ Again $\mathcal{X}_{\mathcal{J}}\;$carries
a tame $G$-action and the branch locus of the cover $\X_{\mathcal{J}%
}$/$\Y_{\mathcal{J}}$ is a divisor with strict normal crossings.
Let ${\cal I}_{p}$ denote the subset of those $i\in \mathcal{I}$ such
that $p_{i}=p$. For ${\cal J}\subset {\cal I}_p$ we write\ \ $f_{\mathcal{J}%
}$ for the structure map $f_{\mathcal{J}}:\X_{\mathcal{J}%
}\rightarrow \Spec({\bf F}_p)$ and as per 6.b in [CPT1] we
set 
\[
\Psi _{p} ( \X_{\mathcal{J}}/\Y_{\mathcal{J}} )
=\left( -1\right) ^{d-\left| {\cal J}\right| }f_{\mathcal{J\ast }} (
c^{d-\left| \mathcal{J}\right| } ( \Omega _{\X_{\mathcal{J}}/%
{\bf F}_{p}} )  ) \;\;\;\hbox{\rm in \ \ K}_{0}( {\bf F}_{p}[G]) 
\]
\[
\Psi _{p}=\tsum_{\phi \neq \mathcal{J\subset I}}\left( -1\right) ^{\left| 
\mathcal{J}\right| +1}\Psi _{p} ( \mathcal{X}_{\mathcal{J}}/\mathcal{Y}_{%
\mathcal{J}} ) 
\]
and 
\[
\Psi =\oplus_{p\in S}\Psi_{p}\ \in\ \oplus_{p\in S}{\rm K}_{0}({\bf F}_{p}[G]) . 
\]
\medskip

\no{\bf Theorem 8.7} \ {\sl (a) The classes $f_{p\ast }F_{p}$ and 
$\left( -1\right) ^{d}\Psi_{p}$  
differ by the class of a free ${\bf F}_{p}[G]$-module
in ${\rm K}_{0}({\bf F}_p[G]) $, and so 
$$
\widetilde{\nu }( f_{p\ast }F_{p}) =\;\widetilde{\nu } ( \Psi _{p} )
^{\left( -1\right) ^{d}}.
$$

(b) The class $\widetilde{\nu } ( \Psi _{p} ) $ 
is represented by the idele valued character homomorphism  $\delta _{p}$%
\[
\delta  ( \theta  ) _{v}=\left\{ 
\begin{array}{c}
\left| \varepsilon _{p} ( \Y_{p},\overline{\theta } )
\right| _{p}\;\;\;\;\hbox{\rm if \ }v=p; \\ 
\ \ \ \ 1\ \;\;\;\;\;\;\;\;\;\;\ \ \ \hbox{\rm if \ }v\neq p.
\end{array}
\right. 
\]}
\smallskip

\no{\sc Proof.} This is the content of (A) and (B) in 6(b) of [CPT1], but note that,
as explained in 4.B, we adopt the opposite convention on the representative
of a class in ${\rm K}_0{\rm T}(\Z[G])$ to that used in [CPT1].\ \ \
\ \ \ \ $\Box \medskip $

We are now in a position to complete the proof of Theorem 8.4. From (35),
(36) and part (a) of Theorem 8.7 we know that 
\[
\widetilde{{\frak{d}}}\cdot \widetilde{\frak{c}}^{-1}=\widetilde{\nu }\left(
i_{S\ast }\left( F\right) \right) ^{\left( -1\right) ^{d}}=\widetilde{\nu }%
\left( \Psi \right) 
\]
and by part (b) of Theorem 8.7 we know that $\widetilde{\nu }\left( \Psi
\right)$ is represented by the finite idele valued homomorphism on
characters $\delta =\tprod_{p}\delta _{p}$. From Theorem 7.5 we know that $%
\widetilde{\frak{c}}$ is a rational class and that $\theta \left( \widetilde{%
\frak{c}}^{s}\right) =\widetilde{\varepsilon }_{0}^{s}\left( \mathcal{Y}%
\right) ^{-1};\;$it therefore follows that $\widetilde{\frak{d}}^{s}$ is
represented by the character function with trivial Archimedean coordinate
and whose finite coordinate is 
\[
\widetilde{\varepsilon _{0}}^{s}( \Y ) ^{-1}\left[
\tprod_{p\in S}\left| \widetilde{\varepsilon }^{s}( \Y_{p}) \right|_{p}
\widetilde{\varepsilon }^{s}( \Y_{p}) \widetilde{\varepsilon }^{s}( \Y_{p}) ^{-1}%
\right] . 
\]

Therefore, to complete the proof of Theorem 8.4, we are now reduced to
showing:
\medskip

\no{\bf Proposition 8.8} \ {\sl The character function $ \tprod_{p\in S}\left| \varepsilon
( \Y_{p}) \right| _{p}\varepsilon ( \Y_{p} ) $  represents the arithmetic ramification class
${\rm AR}(\X) $.}
\medskip

\no{\sc Proof.} For $f,g\in \;\Hom_{\Omega }( R_{G},J_{f}) $ we write $%
f\sim g$ if $f$ and $g$ represent the same class in $A(\Z[G])$. From 8.2 we need to show that

\begin{equation}
\tprod_{p\in S}\left| \varepsilon
( \Y_{p}) \right| _{p}\varepsilon ( \Y_{p} )  \sim \tprod_{p\in
S}\varepsilon \left( b_{p}\right) \left| \varepsilon \left( b_{p}\right)
\right| _{p}.
\end{equation}
With the notation of 8.1 we know that\ for each prime $p\in S,$ 
\[
 \left| \varepsilon
( \Y_{p}) \right| _{p}\varepsilon ( \Y_{p} ) 
 =\left| \varepsilon  ( \Y_{p} ) \right| _{p}(\varepsilon  ( \Y_{p})
_{p}\times \tprod_{q\neq p}\varepsilon  ( \Y_{p} ) _{q}) 
\]
\[
=\left| \varepsilon \left( U_{p}\right) \right| _{p}\varepsilon \left(
U_{p}\right) _{p}\left| \varepsilon \left( b_{p}\right) \right|
_{p}\varepsilon \left( b_{p}\right) _{p}\times \tprod_{q\neq p}\varepsilon
\left( U_{p}\right) _{q}\varepsilon \left( b_{p}\right) _{q}. 
\]
But from Theorem \ 8.1 we know $\left| \varepsilon \left( U_{p}\right)
\right| _{p}\varepsilon \left( U_{p}\right) _{p}\sim 1$ and $\varepsilon
\left( U_{p}\right) _{q}$ $\sim 1$ whenever $q\neq p.$ This then establishes
(37), as required.$\;\;\;\Box $

\bigskip

\noindent{\sc T. Chinburg, \  University of Pennsylvania, Phila., PA
19104.}

 ted@math.upenn.edu

\medskip

\noindent{\sc G. Pappas, \ Michigan State University, E. Lansing, MI 48824.} 

 pappas@math.msu.edu
\medskip

\noindent{\sc M. J. Taylor, \ UMIST, Manchester, M60 1QD, UK.}

 Martin.Taylor@umist.ac.uk

\end{document}